\documentclass{article}

\usepackage{epsfig}
\usepackage{eepic}
\usepackage{pict2e}
\usepackage{latexsym}
\usepackage{amsmath}
\usepackage{amssymb}
\usepackage[english]{babel}
\usepackage{color}

\usepackage{geometry}
\newtheorem{thm}{Theorem}[subsection]
\newtheorem{cor}[thm]{Corollary}
\newtheorem{lem}[thm]{Lemma}
\newtheorem{prop}[thm]{Proposition}
\newtheorem{defn}[thm]{Definition}
\newtheorem{rem}[thm]{Remark}
\newtheorem{expl}[thm]{Example}

\newcommand{\M}{\mathcal{M}}

\newcommand{\Real}{\mathbb{R}}

\DeclareMathOperator*{\bigboxplus}{\text{\Huge $_{\boxplus}$}}

\DeclareMathOperator*{\bigsmileplus}{\text{\Large ${\stackrel{_{+}}{\smile}}$}}
\DeclareMathOperator*{\bigsmileminus}{\text{\Large ${\stackrel{-}{\smile}}$}}

\DeclareMathOperator*{\Card}{\mathrm{Card}}
\DeclareMathOperator*{\Ls}{\mathrm{Ls}}
\DeclareMathOperator*{\Linf}{\mathrm{Li}}
\DeclareMathOperator*{\Lim}{\mathrm{Lim}}
\DeclareMathOperator*{\sgn}{\mathrm{sgn}}

\begin{document}

\title{ Remarks on some     Limit Geometric Properties     related to an  Idempotent and Non-Associative Algebraic  Structure}

 \author{\thanks{University of Perpignan, LAMPS, 52 avenue Villeneuve,
66000 Perpignan, France. email{briec@univ-perp.fr} }Walter Briec}

 \maketitle
 \date

\begin{abstract}

This article analyzes the geometric properties of an idempotent, non-associative algebraic structure that extends the Max-Times semiring. This algebraic structure is useful for studying systems of Max-Times and Max-Plus equations, employing an appropriate notion of a non-associative determinant. We consider a connected ultrametric distance and demonstrate that it implies, among other properties, an analogue of the Pythagorean relation. To this end, we introduce a suitable notion of a right angle between two vectors and investigate a trigonometric concept associated with the Chebyshev unit ball. Following this approach, we explore the potential implications of these properties in the complex plane.

We provide an algebraic definition of a line passing through two points, which corresponds to the Painlev\'e-Peano-Kuratowski limit of a sequence of generalized lines. We establish that this definition leads to distinctive geometric properties; in particular, two distinct parallel lines may share an infinite number of points.

\end{abstract}

{\bf AMS:} 06D50, 32F17\\

{\bf Keywords:} Generalized mean, convexity,
convex hull, ultrametric, semilattice, $\mathbb{B}$-convexity.

\section*{Introduction}\label{SECMXVSP}

It is well known that algebras associated with idempotent semirings are related to particular geometric representations. This article investigates some of these properties through the non-associative symmetrization of Max-Times algebras. The concept of $\mathbb{B}$-convexity was introduced in \cite{bh} as a Painlev\'e-Kuratowski limit of generalized convexities.\footnote{{The letter $\mathbb{B}$, used to define a corresponding notion of convexity  may be confusing in light of the concept of convexity introduced by Anatole Beck in Banach spaces in 1962. It was chosen due to the origins of the names of several members of the mathematics department at the University of Perpignan at the time of the publication of the first paper \cite{bh}. Additionally, it can be linked to the work of Aaron Ben-Tal \cite{ben}, whose paper plays a crucial role in the context of these articles.  Another reason for this choice is that the operation $\boxplus$, pronounced ``$\mathbb{B}$oxplus," bears a graphical resemblance to the Chebyshev unit ball ($\ell_\infty$) when represented in two dimensions with Cartesian coordinates. This ball plays a central role in the proposed approach and is itself the Hausdorff limit of the sequence of $\ell_p$ balls as $p \to \infty$.  
}}
  This fundamental idea is related to earlier work by Ben-Tal \cite{ben} and Avriel \cite{avr1}. When restricted to the nonnegative orthant, $\mathbb{B}$-convexity coincides with Max-Times convexity. In this context, an algebraic formulation has been established in the idempotent Max-Times semiring $(\mathbb{R}_+, \max, \times)$ by applying the transformation $+ \rightarrow \max$ in the nonnegative Euclidean orthant. The idempotent Max-Times semiring is homeomorphic to the Max-Plus semiring $(\mathbb{R} \cup \{-\infty\}, \max, +)$ via a logarithmic transformation, and $\mathbb{B}$-convex sets emerge from a broad class of idempotent algebraic structures introduced in \cite{ms} and derived from the notion of dequantization (see also \cite{lms1}).

Hahn-Banach-type separation properties \cite{bh3} and fixed-point results have been established in \cite{bh2}. Other idempotent convex structures have been proposed in \cite{adilYe} and \cite{adilYeTi}. The algebraic formulation of $\mathbb{B}$-convexity was previously limited to the Max-Times semi-module $(\mathbb{R}_+^n, \vee, \times )$. To address this limitation, an idempotent and non-associative binary operation was introduced in \cite{b15} by considering a special class of idempotent and non-associative magmas. In this paper, we introduce an $n$-ary operation associated with a suitable class of regularized semicontinuous operators. It is worth noting that, when involving only two points, the idempotent and non-associative algebraic structures considered are also mentioned by Gaubert \cite{Gau98} and discussed Viro \cite{ov10}, who explored the complex case. 

{Note that there exist other approaches to extending idempotent algebraic structures. For example, Izhakian and Rowen \cite{ir10} introduced supertropical algebra, which involves a similar type of construction based on the notion of a valued monoid. Additionally, the concept of hyperfields, introduced in \cite{k83} and \cite{m06}, provides an alternative approach.}

Relaxing associativity preserves both symmetry and idempotency. It has been shown that this notion of convexity can be equivalently characterized by the Kuratowski-Painlev\'e limit of the generalized convex hull of two points as defined in \cite{bh}. Some separation results in $\mathbb{R}^n$ have been established in \cite{b17}, and the limiting structure of polytopes has been analyzed in \cite{b19}.

More recently, these structures have been examined from an algebraic perspective to compute the limit solutions of systems of Max-Times equations, particularly by introducing a special determinant formulation. These results also have geometric implications for determining the equation of a limit hyperplane passing through $n$ points, with possible extensions to Max-Plus algebra.

This paper continues prior research on this topic from a geometric perspective, introducing appropriate notions of distance. Throughout the paper, an idempotent analogue of subtraction plays a key role.

First, we define a space notion suited to the algebraic structure under consideration and generalize the definition of $\mathbb{B}$-space introduced in \cite{bh2}. Second, we construct a sequence of generalized metrics via an appropriate homeomorphic transformation of the vector space structure. We demonstrate that their limit properties yield an ultrametric distance intimately linked to the algebraic structure considered. To this end, we derive a notion of orthogonality from the limit inner product defined in \cite{b15}. As a result, we establish geometric properties of a special class of right-angle triples, specifically proving a Pythagorean-type relation within this framework. Using the limiting properties of the generalized inner product, we construct a pseudo-cosine limit function, where the pseudo-cosine is derived from the limiting inner product, and the pseudo-sine stems from the two-dimensional determinant formula established in \cite{b20}. Following this approach, we explore potential implications of these properties in the complex plane.

Finally, the paper investigates a notion of line corresponding to the non-associative algebraic structure under study. We show that the Painlev\'e-Kuratowski limit of a sequence of generalized lines admits an algebraic description extending a previous result in \cite{b15} for a line segment. We establish a condition ensuring that two boundary lines are parallel, analogous (in a certain sense) to the criterion for parallel lines in a vector space. The paper concludes with an analysis of the two-dimensional case, in which a special idempotent geometry is examined. Building on \cite{b20}, we demonstrate how to compute the equation of a two-dimensional line. The final result describes a plane geometry where a unique line passes through two distinct points, and where a line segment can be extended to a half-line. However, two distinct parallel lines may share an infinite number of points.

The paper is structured as follows. Section 1 introduces the framework of the non-associative idempotent algebraic structure defined in \cite{b15}, along with some associated convexity notions. Section 2 defines field and space concepts related to this algebraic structure, providing illustrative examples. Section 3 derives an ultrametric distance as the limit of a sequence of generalized metrics. In this context, a triangular equality is established for a special class of right-angle triples. Section 4 explores inner product properties established in \cite{b17} and introduces a related notion of pseudo-cosine. Section 5 presents an algebraic description of a line. The paper concludes with an analysis of the geometric implications of the considered algebraic structure.

\section{Idempotent and Non-Associative Algebraic and Convex Structure}

\subsection{Isomorphism of vector-Space Structures in Limit}
For all $p\in \mathbb{N}$,
let us consider a bijection $\varphi_p :\mathbb{R} \longrightarrow
\mathbb{R}$ defined by:
\begin{equation}\varphi_p:{\lambda} \longrightarrow {\lambda }^{2p+1}\end{equation}
 and $\phi_p (x_1,... ,x_n)=(\varphi_p(x_1),... ,\varphi_p(x_n))$. This is closely related to the approach proposed by
Ben-Tal   \cite{ben} and Avriel   \cite{avr1}. One can
induce a field structure on $\mathbb{R}$ for which $\varphi_p$
becomes a field isomorphism. {For $\lambda, \mu\in \mathbb{R}$, the indexed sum  
and the indexed product are simply given by  
\begin{equation}  
\lambda \stackrel{p}{+} \mu = \varphi_{p}^{-1} \big(\varphi_p(\lambda) + \varphi_p(\mu)\big) = (\lambda^{2p+1} + \mu^{2p+1})^\frac{1}{2p+1},  
\end{equation}  
and  
\begin{equation}  
\lambda \stackrel{p}{\cdot} \mu = \varphi_{p}^{-1} \big(\varphi_p(\lambda) \cdot \varphi_p(\mu)\big) = \lambda \cdot \mu.  
\end{equation}  }
These two operations define a scalar field $(\mathbb{R}, \stackrel{p}{+}, \cdot)$.    Given this change of notation via
$\varphi_p$ and $\phi_p$, we can define a $\Real$-vector space
structure on $\Real^n$ by: $\lambda \stackrel{\varphi_p}{\cdot}
x=\phi_p^{-1}\big(\varphi_p(\lambda)\cdot \phi_p(x)\big)=\lambda {\cdot} x$ and $x
\stackrel{\varphi_p}{+} y=\phi_p^{-1}\big(\phi_p(x)+\phi_p(y)\big)$.
We call these two operations the  {  {multiplication} by   scalar} and the
indexed sum (indexed by $\varphi_p$).

The $\varphi_p$-sum  denoted $\stackrel{\varphi_p}{\sum}$    of
$(x^{(1)},...,x^{(m)})\in \mathbb{R}^{n\times m}$ is defined by\footnote{For all positive natural numbers  $n$, $[n]=\{1,...,n\}$.}
{ \begin{equation}\stackrel{\varphi_p}{\sum_{{j}\in [m]}}x^{(j)}=
\phi_p^{-1}\Big(\sum_{j\in [m]}\phi_p (x^{(j)})\Big).\end{equation}} For the sake simplicity, throughout the paper we
denote for all $x,y\in \Real^n$:
\begin{equation}x\stackrel{p}{+}y=x\stackrel{\varphi_p}{+}y.\label{phip}\end{equation}
Recall that  Kuratowski-Painlev\'e lower limit of the sequence
of sets $\{A_n\}_{n\in \mathbb N}$, denoted $\Linf_{n\to\infty}A_n$, is the set of
points $x$ for which there exists a sequence $\{x^{(n)}\}_{n\in \mathbb N}$ of points
such that $x^{(n)}\in A_n$ for all $n$ and $x = \lim_{n\to\infty}x^{(n)}$. The Kuratowski-Painlev\'e upper limit of the sequence
of sets $\{A_n\}_{n\in \mathbb N}$, denoted $\Ls_{n\to\infty}A_n$, is the set of
points $x$ for which there exists a subsequence $\{x^{(n_k)}\}_{k\in \mathbb N}$ of points
such that $x^{(n_k)}\in A_{n_k}$ for all $k$ and $x = \lim_{k\to\infty}x^{(n_k)}$.
A sequence $\{A_n\}_{n\in \mathbb{N}}$ of subsets of $\Real^n$ is
said to converge, in the Kuratowski-Painlev\'e sense, to a set $A$
if $\Ls_{n\to\infty}A_n = A = \Linf_{n\to\infty}A_n$, in which case we
write $A = \Lim_{n\to\infty}A_n$.

 In \cite{b15} it was shown that for all $\lambda,\mu\in \Real$ we have:
{\begin{equation}\label{base}\lim_{p\longrightarrow
+\infty}\lambda \stackrel{p}{+}\mu=
\left\{\begin{matrix}{\lambda}\ &\hbox{ if } &|\lambda|&>&|\mu|\\
\frac{1}{2}(\lambda +\mu )&\hbox{ if }&|\lambda|&=&|\mu|\\
\mu& \hbox{ if }& |\lambda|&<&|\mu|.\end{matrix}\right.\end{equation}} Along
this line one can introduce the binary operation $\boxplus$
defined for all $\lambda,\mu\in \Real$ by:
{\begin{equation}
\lambda \boxplus \mu=\lim_{p\longrightarrow +\infty}\lambda \stackrel{p}{+}\mu.
\end{equation}}
Though the operation $\boxplus$ does not satisfy associativity, it
can be extended by constructing a non-associative algebraic
structure which returns to a given $n$-tuple a real value. For all
$x\in \mathbb R^n$ and all subsets $I$ of $[n]$, let us  consider
the map $\xi_I[x]:\Real \longrightarrow \mathbb Z$ defined for all
$\alpha \in \Real$ by
\begin{equation}\label{defxi}\xi_I[x](\alpha)= \Card
\{i\in I: x_i=\alpha\}- \Card \{i\in I: x_i=-\alpha\}.\end{equation}
This map measures the symmetry of the occurrences of a given value
$\alpha $ in the {coordinates} of a vector $x$.

For all $x\in \mathbb R^n$ let $\mathcal J_I(x)$ be a subset of $I$
defined by
\begin{equation}\mathcal J_I(x)=\Big\{j\in I: \xi_I[x](x_j)\not=0\Big\}=I\backslash \big \{i\in I: \xi_I[x](x_i)=0\big\}.\end{equation}
$\mathcal J_{I} (x)$ is called {\bf the residual  index set } of
$x$. It is obtained by dropping from $I$ all the $i$'s such that
$\Card \{j\in I: x_j=x_i\}= \Card \{j\in I: x_j=-x_i\}$.

For all positive natural numbers $n$ and
for all subsets $I$ of $[n]$, let $\digamma_I: \Real^n
\longrightarrow \Real $ be the map defined for all $x\in \Real^n$ by

\begin{equation}\digamma_{ I}(x)=\left\{\begin{matrix}\max_{i\in \mathcal
J_I (x)}x_{i} &\text{ if }&\mathcal J_I(x)\not=\emptyset&\text{ and } \quad \xi_I[x]\big (\max_{i\in
\mathcal J_I(x)}|x_i|\big)>0 \\
\min_{i\in \mathcal J_I(x)}x_i &\hbox{ if }&\mathcal J_I(x)\not=\emptyset&\text{ and }\quad \xi_I[x]\big(\max_{i\in \mathcal J_I(x)}|x_i|\big)<0\\
 0 &\text{ if }&{\mathcal J_I(x)=\emptyset {,}} &  
\end{matrix}\right.\end{equation}
where $\xi_I[x]$ is the map defined in \eqref{defxi} and $\mathcal
J_I(x)$ is the residual index set of $x$. The operation that takes an $n$-tuple $(x_1, \ldots, x_n)$ from $\mathbb{R}^n$ and returns a single real value $\digamma_I(x_1, \ldots, x_n)$ is called an  {\bf $n$-ary extension} of the binary operation $\boxplus$, for all natural numbers $n \geq 1$ and all $x \in \mathbb{R}^n$, where $I$ is a nonempty subset of $[n]$.

According to \cite{b15}, for any $n$-tuple $x = (x_1, \ldots, x_n)$, the operation can be defined as:
\begin{equation} \label{defrecOp}
\bigboxplus_{i \in I} x_i = \lim_{p \to \infty} \stackrel{\varphi_p}{\sum_{i \in I}} x_i = \digamma_I(x).
\end{equation}

Clearly, this operation encompasses as a special case the binary operation
{  mentioned in  \cite{Gau98}} and defined in equation \eqref{base}. For all
$(x_1,x_2)\in \Real^2$:
\begin{equation*}\bigboxplus_{i\in \{1,2\}}x_i=x_1\boxplus
x_2.\end{equation*} For example,  if $x=(-3,-2,3,3,1,-3)$,   we
have $F_{[6]}(-3,-2,3,3,1,-3)=F_{[2]}(-2,1)=-2=\bigboxplus_{i\in
[6]}x_i$. There are some basic properties that can be inherited
from the above algebraic structure: $(i)$ for all $\alpha\in
\Real$, one has: $ \alpha \Big(\bigboxplus_{i\in I}x_i\Big)=
\bigboxplus_{i\in I}(\alpha x_i)$;  $(ii)$ suppose that
$x\in\epsilon \Real_+^n$ where $\epsilon$ is $+1$ or $-1$. Then
$\bigboxplus_{i\in I}x_i=\epsilon \max_{i\in I} \{\epsilon x_i\}$; $(iii)$ we have $|\bigboxplus_{i\in I} x_i|\leq \bigboxplus_{i\in
I}|x_i|$; $(iv)$ for all $x\in \Real^n$:
\begin{equation}\label{Decompos}\Big[x_i \boxplus \big(\bigboxplus_{j\in I\backslash
\{i\}}x_j\big)\Big]\in \Big \{0, \bigboxplus_{j\in I}x_j\Big \}\quad
\text{ and }\quad \bigboxplus_{i\in I}x_i=\bigboxplus_{i\in
I}\Big[x_i \boxplus \big(\bigboxplus_{j\in I\backslash
\{i\}}x_j\big)\Big].\end{equation}

The algebraic structure $(\Real,\boxplus,\cdot)$  can be extended
to $\Real^n$. Suppose that $x,y\in \Real^n$, and let us denote $
x\boxplus y=(x_1\boxplus y_1,...,x_n\boxplus y_n).$  Moreover,
let us consider $m$ vectors $x^{{(1)}} ,...,x^{{(m)}}\in \Real^n$, and define
\begin{align}
\bigboxplus_{j\in [m]}x^{{(j)}}&=\Big(\bigboxplus_{j\in [m]}
x_{ 1}^{{(j)}},...,\bigboxplus_{j\in [m]} x_{ n}^{{(j)}}\Big).
\end{align}

The $n$-ary operation $(x_1,...,x_n)\rightarrow \bigboxplus_{i\in [n]}x_i$ is not associative.
To simplify the notations of the paper, for  all $z\in \{z_{i_1,...,i_m}: i_k\in I_k,k\in [m]\}$, where $I_1,...,I_m$ are $m$ index subsets of $\mathbb N$, we use the notation:
\begin{equation}
\bigboxplus_{ \substack{i_k\in I_k\\k\in [m]}}z_{i_1,...,i_m}=\bigboxplus_{\substack{(i_1,...,i_m)\in \prod_{k\in [m]} I_k}}z_{i_1,...,i_m}.
\end{equation}
 { Notice that since the operation $\boxplus$ is not associative, it can be ambiguous to apply without using the symbol $\bigboxplus$ indexed on a given finite subset $I$. In the remainder, for the sake of simplicity, and when it is more convenient, we will adopt the following notational convention.}
 For all $x\in \Real^n$:
\begin{equation}
\bigboxplus_{i\in [n]}x_i=x_1\boxplus\cdots\boxplus x_n.
\end{equation}
 
In the sequel, for all $x,y\in \Real^n$  we will often use the following notation to define an idempotent analogue of the subtraction operation:
\begin{equation}
x\boxminus y=x\boxplus(-y).
\end{equation}

 \subsection{Scalar Products and Determinants in Limit}

\noindent This section presents the algebraic properties induced
by an isomorphism between a scalar field and a set and their
implications on the scalar product.   If $\langle
\cdot ,\cdot \rangle$ is the canonical inner product over $\Real^{n}$, then
there exists a symmetric bilinear form
 $\langle \cdot,\cdot \rangle_{\varphi_p}:
 \Real^n\times \Real^n\longrightarrow \Real$ defined by:
\begin{equation}\langle x,y\rangle_{\varphi_p}=\varphi_p^{-1}\big(\langle \phi_p^{}(x),
\phi_p^{}(y)\rangle\big )=\big(\sum_{i\in [n]}x_i^{2p+1} y_i^{2p+1}\big)^{\frac{1}{2p+1}}.\end{equation}
Now, let us denote $\big [\langle y,\cdot\rangle_{\varphi_p}\leq
\lambda\big ]= \left\{x\in \Real^n: \langle
y,x\rangle_{\varphi_p}\leq \lambda\right\}$. For the sake of simplicity, let us denote $\langle \cdot,\cdot \rangle_p$ this scalar product.

In what follows, we state some results obtained  in \cite{b15,b17} by taking the limit. We first introduce the operation $ \langle \cdot,\cdot\rangle_\infty :\Real^n\times \Real^n\longrightarrow
 \Real$ defined for all $x,y\in \Real^n$ by $\langle x,y\rangle_\infty =\bigboxplus_{i\in
 [n]}x_iy_i$. Let $\|\cdot\|_\infty$ be the Tchebychev
 norm defined by $\|x\|_\infty=\max_{i\in [n]}|x_i|$. It is established in \cite{b15} that
for all $x,y\in \Real^n$, we have:
 $(i)$  $\sqrt{\langle x,x\rangle_\infty} =\|x\|_\infty$;
  $(ii)$  $|\langle x,y\rangle_\infty| \leq \|x\|_\infty \|y\|_\infty$;
   $(iii)$ for all $\alpha\in \Real$, $\alpha \langle x,y\rangle_\infty= \langle \alpha x,y\rangle_\infty=\langle  x,\alpha y\rangle_\infty$.
By definition, we have for all $x,y\in \Real^n$:

\begin{equation}\langle x,y\rangle_\infty=\lim_{p\longrightarrow \infty}\langle x,y\rangle_p.\end{equation}

A map $f:\Real^n\longrightarrow \Real $ defined as $f(x)=\langle a, x\rangle_\infty$ for some $a\in \Real^n$ is called an {\bf idempotent symmetric form}. In the following, for all maps $f:\Real^n\longrightarrow \Real$ and
all real numbers $c$, the notation $[\,f\leq c\,]$ stands for the
set $f^{-1}(\,]-\infty, c]\,)$.  Similarly, $[\,f< c\,]$ stands
for  $f^{-1}(\,]-\infty, c[\,)$ and $[\,f\geq c\,]=[\,-f\leq
-c\,]$.

  For all $u,v\in \Real$, let us define the binary operation
\begin{equation*}\label{Bform}u\stackrel{-}{\smile} v=
\left\{\begin{matrix}
u &\hbox{ if } &|u|& > &|v|\\
\min \{u,v\}&\hbox{ if }&|u|&=&|v|\\
v& \hbox{ if }& |u|&<&|v|.\end{matrix}\right.\end{equation*} 
{Note that in \cite{ir10}, a polynomial theory for supertropical algebra was proposed. To achieve this, a similar operation was constructed by extending the tropical semiring with a suitable valuation function. Both proposed constructions follow a similar principle. However, in our approach, this binary operation is derived from the idempotent semiring $(\Real_+, \max, \cdot)$ and the absolute value function.}

An
elementary calculus shows that $u\boxplus v=\frac{1}{2}\Big[u
\stackrel{-}{\smile} v-\big[(-u) \stackrel{-}{\smile}(-v)\big]
\Big]$. It was shown in \cite{bh3}, that this operation is associative. 

 Given $m$ elements $u_1, ..., u_m$ of $\Real$, not all of which are $0$, let $I_+$, respectively $I_-$,
 be the set of indices for which $0 < u_i$, respectively $u_i < 0$. We can then write
 $u_1\stackrel{-}{\smile}\cdots\stackrel{-}{\smile} u_m = (\stackrel{-}{\smile}_{i\in I_+}u_i) \stackrel{-}{\smile} (\stackrel{-}{\smile}_{i\in I_-}u_i)
 = (\max_{i\in I_+}u_i)\stackrel{-}{\smile}(\min_{i\in I_-}u_i)$ from which we have

 \begin{equation}\label{manysmiles}
 u_1\stackrel{-}{\smile}\cdots\stackrel{-}{\smile} u_m =
\left\{
\begin{array}{lcc}
\max_{i\in I_+}u_i  &\hbox{if}   & I_- = \emptyset \hbox{ or } \max_{i\in I_-}\vert u_i\vert <  \max_{i\in I_+}u_i    \\
\min_{i\in I_-}u_i  &  \hbox{if} &  I_+ = \emptyset \hbox{ or } \max_{i\in I_+}u_i <  \max_{i\in I_-}\vert u_i\vert \\
 \min_{i\in I_-}u_i  & \hbox{if}  &    \max_{i\in I_-}\vert u_i\vert  = \max_{i\in
 I_+}u_i.
\end{array}
\right.
 \end{equation}
We define an {\bf     idempotent lower symmetric form }on ${\mathbb R}^n$
as a map $g: {\mathbb R}^n\to\Real$ such that for all $(x_1, ...,
x_n)\in{\mathbb R}^n$,
 \begin{equation}\label{eqcarBform}
 g(x_1, ..., x_n) = {\langle a,x\rangle}_\infty^{-} =a_1x_1 \stackrel{-}{\smile}\cdots \stackrel{-}{\smile} a_nx_n.
 \end{equation}
It was established in \cite{b15} that for all $c \in \mathbb{R}$, 
the set $g^{-1}(\,]-\infty,c]) = \{x \in \mathbb{R}^n \mid g(x) \leq c\}$ is closed. It follows that  an     idempotent lower symmetric form  is
lower semi-continuous. It was also shown in \cite{b15} that for any  $a\in \Real^n$,   the    idempotent lower symmetrical   form  $g$ defined by $ g(x_1,
..., x_n) = a_1x_1\stackrel{-}{\smile}\cdots\stackrel{-}{\smile}
a_nx_n, $    is the lower
semi-continuous regularization of the map $x\mapsto \langle
a,x\rangle_{\infty}=\bigboxplus_{i\in [n]}a_ix_i$. 

Similarly, one can introduce a binary operation defined for all $u,v \in \mathbb{R}$ as  
\begin{equation}
u \stackrel{+}{\smile} v = -\big[(-u) \stackrel{-}{\smile} (-v)\big],
\end{equation}  
which is associative. Along these lines, one can define an {\bf idempotent upper symmetric form} as a map $h: {\mathbb R}^n\to\Real$  such that, for all $(x_1, ...,
x_n)\in{\mathbb R}^n$,
 \begin{equation}\label{eqcarBform}
 h(x_1, ..., x_n) = {\langle a,x\rangle}_\infty^{+} =a_1x_1 \stackrel{+}{\smile}\cdots \stackrel{+}{\smile} a_n x_n.
 \end{equation}
 For all $x\in \Real^n$, we clearly, have the following identities
\begin{equation}
 \langle a,x\rangle_\infty^{+}=-\langle a,-x\rangle_\infty^{-}
 \text{ and }
 \langle a,x\rangle_\infty^{-}=-\langle a,-x\rangle_\infty^{+}.
 \end{equation}
For  any  $a\in \Real^n$, the  idempotent upper symmetric    form $h$ defined by $ h(x_1,
..., x_n) = a_1x_1\stackrel{+}{\smile}\cdots\stackrel{+}{\smile}
a_nx_n, $,  is the upper
semi-continuous regularization of the map $x\mapsto \langle
a,x\rangle_{\infty}=\bigboxplus_{i\in [n]}a_ix_i$.

These forms can be related  to the construction of a ring involving a balance relation   symmetrizing the semi-ring $(\Real_+, \vee, \cdot)$ (see \cite{b20}, see also \cite{h85} and \cite{MP} in a Max-Plus context).

More, recently, the following notion of determinant was introduced in \cite{b20}. The set $\mathbb{R}$, equipped with the operation $\boxplus$, allows for the definition of a suitable notion of determinant (see \cite{b19,b20}).

\noindent We denote by $\mathcal{M}_n(\mathbb{R})$ the set of { square matrices of order $n$ with real entries.} We introduce over $\M_{n}(\Real)$ the
map $\Phi_p:\M_{n}(\Real)\longrightarrow \M_{n}(\Real)$ defined
for all matrices $A= \left(a_{i,j}\right)_{\substack
{i=1...n\\j=1..n}} \in \M_{n}(\Real)$ by:
\begin{align}
\Phi_p(A)&=\big(\varphi_p(a_{i,j})\big)_{\substack
{i=1...n\\j=1...n}}= \left({a_{i,j}}^{2p+1}\right)_{\substack
{i=1...n\\j=1...n}}.\end{align} Its reciprocal is the map
$\Phi_p^{-1}:\M_{n}(\Real)\longrightarrow \M_{n}(\Real)$ defined
by:\begin{align}
\Phi_p^{-1}(A)&=\left(\varphi_p^{-1}(a_{i,j})\right)_{\substack
{i=1...n\\j=1...n}}=\left({a_{i,j}}^{\frac{1}{2p+1}}\right)_{\substack
{i=1...n\\j=1...n}}.
\end{align}
$\Phi_p$ is a natural extension of the map $\phi_p$ from $\Real^n$
to $\mathcal M_n(\Real)$. $\Phi_p(A)$ is the $2p+1$ Hadamard power
of matrix $A$. In the following we introduce the matrix
product:
\begin{equation}
A \stackrel{p}{.}x=\phi_p^{-1}\big(\Phi_p(A).\phi_p(x)\big).
\end{equation}

It is easy to see that the map $x\mapsto A\stackrel{p}{\cdot}x$ is
$\varphi_p$-linear, that is linear with respect to the algebraic operation $\stackrel{p}{+}$ and the usual scalar multiplication.  Conversely,  if $g$ is a $\varphi_p$-linear
map then it can be represented by a matrix $A$ such that $g(x)=
A\stackrel{p}{\cdot}x$ for all $x\in \Real^n$.

 For any $n \times n$ matrix $A$, let $|A|$ denote its determinant. We now introduce the following definition of the $\varphi_p$-determinant:

\begin{equation}|A|_p=\varphi_p^{-1}|\Phi_p(A)|.\end{equation}
 Using the Leibnitz formula yields
\begin{equation}|A|_p=\Big(\sum_{\sigma\in \mathfrak S_n}\sgn(\sigma) \prod_{i\in [n]}a_{i,\sigma(i)}^{2p+1}\Big)^{\frac{1}{2p+1}}, \end{equation}
{where $\mathfrak S_n$ is the set of all the permutations defined on $[n]$. }

Along this line, one can define   a determinant in limit as follows:
\begin{equation}\lim_{p\longrightarrow \infty}|A|_p: =|A|_\infty=\bigboxplus_{\sigma\in \mathfrak S_n} \big(\sgn(\sigma)
 \prod_{i\in [n]}a_{i,\sigma(i)}\big).\end{equation}
This notion is compared in \cite{b20} to that of permanent originally defined in the framework of Max-Plus algebras (se \cite{butheged84} and \cite{b10}).
\subsection{An Idempotent and Non-associative Convex Structure }

In ~\cite{bh}, $\mathbb{B}$-convexity is introduced as a limit of
linear convexities. The
$\phi_p$-convex hull of  a finite set $A=\{x^{(1)},...,x^{(m)}\}\subset \mathbb{R}^n$ is
defined by:
\begin{equation}Co^{\phi_p}(A)=\Big\{\stackrel{\varphi_p}{\sum_{i\in [m]}}t_i\stackrel{\varphi_p}{\cdot}x^{(i)} :
\stackrel{\varphi_p}{\sum_{i\in [m]}}t_i=\varphi_p^{-1}(1),  \varphi_p(t_i)\geq
0, i\in [m]\Big\}\end{equation}
which can be rewritten {using the fact that $\varphi_p^{-1}(1)=1$ and $\varphi_p$ is increasing:}
\begin{equation}Co^{\phi_p}(A)=\Bigg\{\phi_p^{-1}\Big(\sum_{i\in [m]}t_i^{2p+1}{\cdot}\phi_p(x^{(i)})\Big) :
\big(\sum_{i\in [m]}t_i^{2p+1}\big)^{\frac{1}{2p+1}}=1 ,  t_i\geq
0, i\in [m]\Bigg\}.\end{equation}

  Moreover, for all $L\subset \mathbb{R}$, we simplify the notations
denoting $Co^{p}(L)=Co^{\phi_p}(L)$. The definition of $\mathbb{B}$-convex sets is based on
the definition of the Kuratowski-Painlev\'e limit of a sequence of
$\phi_p$-convex sets.
For a non-empty finite subset $A\subset \Real_{}^{n}$ we let
$Co^\infty(A)=\Ls_{p\longrightarrow \infty}Co^p(A)$.  A subset $C$ of $\Real^n$ is $\mathbb B$-convex if for all finite subsets $A$ of $C$, $Co^\infty(A)\subset C$.  The structure of $\mathbb{B}$-convex sets in $\Real_+^n$ will involve the order structure, with respect to the
positive cone of  $\Real^{n}$;  we denote by $\bigvee_{i\in [m]}
x^{(i)}$ the least upper bound of $x^{(1)},...,x^{(m)}\in \Real^{n}$,
that is:

\begin{equation*}
{\bigvee_{i\in [m]} x^{(i)}= \big(\max\{x_{ 1}^{(1)},...,x_{1}^{(m)}\},
\ldots, \max \{x_{ n}^{(1)},...,x_{n}^{(m)}\} \big)}\end{equation*}

\medskip

When $A={\{x^{(1)},...,x^{(m)}\}}\subset \Real_+^n$, it was shown in ~\cite{bh} that:
$$Co^\infty(A)=
\Big\{\bigvee_{i\in [m]}t_i x^{(i)}: t_i \in [0,1], \max_{i\in [m]} t_i =1\Big\}=\Lim_{p\longrightarrow \infty}Co^p(A).$$
 Moreover it is also proved
that the sequence is Hausdorff convergent. In addition, if the set $A$ is contained in any $n$-dimensional orthant, it was established in \cite{b15} that  $
Co^\infty(A) = \Big\{\bigboxplus_{i\in [m]}t_i x^{(i)}: t_i \in [0,1], \max_{i\in [m]} t_i =1\Big\}.
$

 If $L\subset \Real_+^n$ is
 $\mathbb{B}$-convex, then $L=\Lim_{p\longrightarrow \infty}Co^p(L)=\bigcap_{r=0}^{\infty}Co^p(L)$.
For all $S\subset \Real_+^n$, we denote $Co^{\infty}(S)=\Lim_{p\longrightarrow
\infty}Co^p(S)$. Properties of this operation are developed in ~\cite{bh}. In
particular, for all $S\subset \Real_+^n$, $Co^{\infty}(S)$ is $\mathbb{B}$-convex or  equivalently Max-Times convex.

An important limitation of this definition of $\mathbb B$-convexity is that its algebraic form is restricted to the non-negative orthant $\Real_+^n$, over which the idempotent Max-Times semiring is considered. To circumvent this difficulty an extended definition of $\mathbb B$-convexity other $\Real^n$ was suggested in \cite{b15}. A subset \( C \) of \( \mathbb{R}^n \) is said to be {\textbf{idempotent symmetric  convex}\footnote{Referred to as \(\mathbb{B}^\sharp\)-convexity in \cite{b15}, this concept is given a more intuitive terminology here.
}} if, for all \( x,y \in C \) and for all \( t \in [0,1] \),  
\begin{equation}  
x \boxplus ty \in C.  
\end{equation}

It was shown in \cite{b15}  that the sequence of subsets
$\{Co^{(p)}(x,y)\}_{p\in \mathbb N}$ of $\Real^n$ has a Painlev\'e-Kuratowski  limit denoted $Co^\infty(x,y)$. These points can be expressed trough an explicit algebraic formula in term of the idempotent and non-associative algebraic structure defined on $(\Real^n, \boxplus, \cdot)$.

\begin{center}
 { \scriptsize % This is a LaTeX picture output by TeXCAD.
% This is a LaTeX picture output by TeXCAD.
% File name: [Clipboard].
% Version of TeXCAD: 4.51
% Reference / build: 27-Nov-2018 (rev. a75)
% For new versions, check: http://texcad.sf.net/
% Options on the following lines.
%\grade{\on}
%\emlines{\off}
%\epic{\off}
%\beziermacro{\on}
%\reduce{\on}
%\snapping{\off}
%\pvinsert{% Your \input, \def, etc. here}
%\quality{8.000}
%\graddiff{0.005}
%\snapasp{1}
%\zoom{4.0000}
\unitlength 0.4mm % = 2.845pt
\linethickness{0.4pt}
\ifx\plotpoint\undefined\newsavebox{\plotpoint}\fi % GNUPLOT compatibility

\begin{picture}(182.25,157.5)(0,0)
{
\put(86.25,16.5){\vector(0,1){135}}
}
\put(146.25,142.5){$x^{(1)}$}
\put(116.25,106.5){\line(-1,0){60}}
\put(38.25,61.5){$y^{(1)}$}
\put(56.25,106.5){\line(0,-1){45}}
\put(62.25,142.5){$x^{(2)}$}
\put(32.25,106.5){$y^{(2)}$}
\put(56.25,136.5){\line(-1,0){30}}
\put(26.25,136.5){\line(0,-1){30}}
\put(146.25,112.5){$x^{(5)}$}
\put(146.25,52.5){$y^{(5)}$}
\put(127.25,95.75){$y^{(6)}$}
\put(65.75,58.5){$x^{(6)}$}
\put(146.25,106.5){\line(0,-1){45}}
\put(149.25,7.5){$x^{(4)}$}
\put(41.25,46.5){$y^{(4)}$}
\put(116.25,46.5){\line(-1,0){60}}
\put(11.25,16.5){$x^{(3)}$}
\put(62.25,28.5){$y^{(3)}$}
\put(41.25,31.5){\line(1,0){15}}
\put(91,70.75){$0 $}
\put(84.75,157.5){$x_2$}
\put(182.25,75){$x_1$}
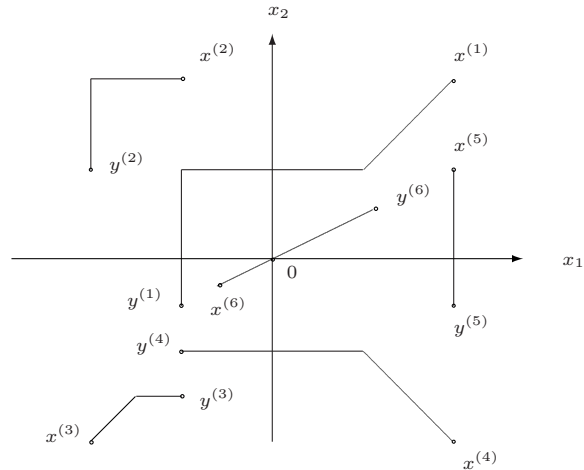
\put(36.5,0){{\bf Fig. 1:} Examples of sets $Co^\infty(x,y)$. }
\put(26.25,106.5){\circle*{1}}
\put(56.75,136.5){\circle*{1.118}}
\put(146.25,135.75){\circle*{1}}
\put(56.25,61.5){\circle*{1.118}}
\put(146.25,106.5){\circle*{1.118}}
\put(69,68.25){\circle*{1.118}}
\put(120.5,93.5){\circle*{1.118}}
\put(146.25,61.5){\circle*{1}}
\put(56.25,46.25){\circle*{1.118}}
\put(146.25,16.5){\circle*{1}}
\put(26.5,16.25){\circle*{1.118}}
\put(56.5,31.5){\circle*{1.118}}
\put(0,77){\vector(1,0){169.25}}
\put(86.5,76.75){\circle*{1}}
%\emline(145.75,136)(116.5,106.5)
\multiput(145.75,136)(-.03373702422,-.03402537486){867}{\line(0,-1){.03402537486}}
%\end
%\emline(145.75,16.75)(116.25,46.5)
\multiput(145.75,16.75)(-.03371428571,.034){875}{\line(0,1){.034}}
%\end
\put(26.25,16.75){\line(1,1){14.75}}
%\emline(68,68)(119.5,93.25)
\multiput(68,68)(.06875834446,.03371161549){749}{\line(1,0){.06875834446}}
%\end
\end{picture}

}
\end{center}

Moreover, it was established in \cite{b15} that this limit entirely characterizes an idempotent symmetric  convex set. A subset $C$ of $\mathbb R^n$ is idempotent symmetric  convex if and only if for all $x,y\in \mathbb R^n$
$Co^\infty(x,y)\subset C$. More recently, it has been established in \cite{b19}  that for all $x,y\in \Real^n$,
$Co^\infty(x,y)$ is {idempotent symmetric  convex} and thereby $Co^\infty(x,y) $ is the smallest {idempotent symmetric  convex} set that contains $\{x,y\}$.
\medskip

\begin{center}{\scriptsize % This is a LaTeX picture output by TeXCAD.
% File name: [Clipboard].
% Version of TeXCAD: 4.51
% Reference / build: 27-Nov-2018 (rev. a75)
% For new versions, check: http://texcad.sf.net/
% Options on the following lines.
%\grade{\on}
%\emlines{\off}
%\epic{\off}
%\beziermacro{\on}
%\reduce{\on}
%\snapping{\off}
%\pvinsert{% Your \input, \def, etc. here}
%\quality{8.000}
%\graddiff{0.005}
%\snapasp{1}
%\zoom{4.0000}
% This is a LaTeX picture output by TeXCAD.
% File name: [Clipboard].
% Version of TeXCAD: 4.51
% Reference / build: 27-Nov-2018 (rev. a75)
% For new versions, check: http://texcad.sf.net/
% Options on the following lines.
%\grade{\on}
%\emlines{\off}
%\epic{\off}
%\beziermacro{\on}
%\reduce{\on}
%\snapping{\off}
%\pvinsert{% Your \input, \def, etc. here}
%\quality{8.000}
%\graddiff{0.005}
%\snapasp{1}
%\zoom{4.0000}
\unitlength 0.4mm % = 2.845pt
\linethickness{0.4pt}
\ifx\plotpoint\undefined\newsavebox{\plotpoint}\fi % GNUPLOT compatibility
\begin{picture}(163.75,151)(0,0)
\put(0,70.75){\vector(1,0){153}}
%\emline(149.75,133.75)(115.5,107.25)
\multiput(149.75,133.75)(-.04357506361,-.03371501272){786}{\line(-1,0){.04357506361}}
%\end
\put(115.5,107.25){\line(-1,0){89.5}}
%\emline(26,107.25)(25.75,97.5)
\multiput(26,107.25)(-.03125,-1.21875){8}{\line(0,-1){1.21875}}
%\end
%\emline(25.75,97.5)(40.5,88.25)
\multiput(25.75,97.5)(.0536363636,-.0336363636){275}{\line(1,0){.0536363636}}
%\end
\put(40.5,88.25){\line(0,-1){33}}
%\emline(115.75,107.25)(102.75,98.75)
\multiput(115.75,107.25)(-.0515873016,-.0337301587){252}{\line(-1,0){.0515873016}}
%\end
\put(102.75,98.75){\line(0,-1){43.75}}
\put(41.25,55.75){\line(1,0){61.5}}
\put(22.75,37.75){\framebox(91.25,10)[cc]{}}
\put(163.75,71.5){\makebox(0,0)[cc]{$x_1$}}
\put(65.25,151){\makebox(0,0)[cc]{$x_2$}}
\put(59.25,65.5){\makebox(0,0)[cc]{$0$}}
\put(9.5,100.75){\line(0,1){47.75}}
\put(9.5,148.5){\line(1,0){24.75}}
\put(65.25,8){\vector(0,1){137.75}}
\put(38.75,28){\line(1,0){63.5}}
\put(107.5,24.5){\line(-1,0){71}}
\put(94.25,142.5){\line(1,0){24.25}}
\put(118.5,142.5){\line(0,-1){19.5}}
\put(65,0){\makebox(0,0)[cc]
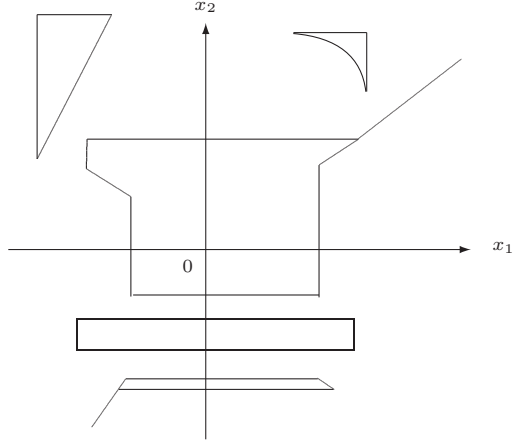
{{\bf Fig. 2:} Examples of {idempotent symmetric  convex  sets}}}
%\emline(107.75,24.5)(102.25,28.25)
\multiput(107.75,24.5)(-.049107143,.033482143){112}{\line(-1,0){.049107143}}
%\end
%\emline(34,148.25)(9.5,100.75)
\multiput(34,148.25)(-.03370013755,-.06533700138){727}{\line(0,-1){.06533700138}}
%\end
%\emline(38.75,28)(27.5,12)
\multiput(38.75,28)(-.0336826347,-.0479041916){334}{\line(0,-1){.0479041916}}
%\end
\qbezier(94.25,142.25)(116.25,140.375)(118.25,123)
\end{picture}

}
\end{center}

\section{Idempotent Symmetric and Ultrametric Spaces}\label{bspace}
The concept of $\mathbb{B}$-space, as introduced in \cite{bh3}, is restricted to the non-negative $n$-dimensional orthant $\mathbb{R}_+^n$. In this work, we extend it to the entire Euclidean vector space. This new definition will be useful for analyzing the algebraic structure of a line defined on $(\mathbb{R}^n, \boxplus, \cdot)$.

\subsection{Idempotent Symmetric Fields and Spaces}
  A finite-dimensional $\mathbb B$-space (Max-Times space) is, by definition, a subset $X$ of $\Real_+^n$
  such that $0\in X$, for all $t\geq 0$ and all $x\in X$, $tx\in X$ and for all $x,y\in X$, $x\vee ty\in X$.  

We first define a suitable notion of {\bf  idempotent pseudo-field}.   By comparing this definition with that of a field, note that we consider a weakened form of associativity  and idempotence is added. 

\begin{defn}A triple $(K,\boxplus_K,\cdot_K)$ is an idempotent pseudo-field 
  if:

\noindent $(a)$ There exits a binary operation $\boxplus_K:K\times
K\longrightarrow K$ that satisfies the following properties: $(i)$ for all $\lambda \in K$ we have $\lambda\boxplus_K \lambda=\lambda$ (idempotence); $(ii)$ for all $(\lambda,\mu)\in K\times K$ we have $\lambda \boxplus_K \mu=\mu\boxplus_K \lambda$ (commutativity); $(iii)$ there exists a neutral element $0_K\in K$ such that $\lambda \boxplus_K 0_K=0_K\boxplus_K \lambda=\lambda$ for all $\lambda\in K$; $(iv)$ every $\lambda\in K$ has a symmetric  element $-\lambda$ such  that $(-\lambda)\boxplus_K \lambda
 =\lambda\boxplus_K (-\lambda)=0_X$ (symmetry); $( v)$ for any mutually non-symmetric triple  of elements $(\lambda, \mu, \eta)\in K^3$, $(\lambda\boxplus_K \mu)\boxplus_K\eta=\lambda\boxplus_K (\mu \boxplus_K\eta)$ (weakened form of associativity).

\noindent $(b)$ There exists a scalar multiplication  $\cdot_K  :K\times
K\longrightarrow K$ that satisfies the following properties: $(i)$ for all $\lambda \in K$ and all $(\lambda,\mu)\in K\times K$, we
have $\lambda \cdot_K   \mu =  \mu\cdot_K  \lambda$ (commutativity);  $(ii)$ for all $\lambda,\mu, \eta \in K\times K\times K$ we have $(\lambda
\cdot_K   \mu)\cdot_K  \eta= \lambda \cdot_K   (\mu \cdot_K 
\eta)$ (associativity); $(iii)$ there is an element $1_K\in K$ such that for all $\lambda\in K$ we have $1_K\cdot_K  \lambda=\lambda$; $(iv)$ for all $\lambda\in K\backslash \{0_K\}$ there exists an element $\lambda^{-1}$ such that $\lambda \cdot_K \lambda^{-1}=1_K$.  
 
 \noindent $(c)$ The multiplication   $\cdot_K  $ is distributive over the operation $\boxplus_K$: for all $\lambda \in K$ and all $(\mu,\nu)\in K\times K$, we
have $$\lambda \cdot_K   (\mu \boxplus_K \nu)= (\lambda \cdot_K  \mu) \boxplus_K (\lambda \cdot_K  \nu).$$

 \end{defn}

\begin{expl}From \textnormal{\cite{b15}}, the triple $(\Real, \boxplus, \cdot)$, where $\boxplus$ is the binary operation defined in \eqref{base} is an idempotent pseudo-field.\end{expl}

In the following, when we use the $\boxplus$ symbol without further specification on $\Real^n$ this will mean that we are considering the binary operation defined in equation \eqref{base} and its related $n$-ary extension. 

\begin{expl}We consider the case of the non-associative symmetrisation of the Max-Plus semi-module analysed in \cite{b20}. Let $\mathbb M=\Real\cup \{-\infty\}$ and let us denote $(\mathbb M, \oplus, \otimes)$ the Maslov’s semi-module where we replace the operations $\vee$ with $\oplus$ and $+$ with $\otimes$. More precisely one can define on $\mathbb M$ the operations $\oplus$ and $\otimes$ respectively as 
$x\oplus y=\max\{x,y\}$ and $x\otimes y=x+y$, where  $-\infty$ is the neutral element of the operation $\oplus$.  Suppose now
that $x\in \Real_-$ and let us extend the logarithm function to
the whole set of real numbers. This we do by introducing the set
\begin{equation}\widetilde {\mathbb M}=\mathbb M \cup (\Real +\mathrm{i}\pi)\end{equation}
where $\mathrm{i}$ is the complex number such that $\mathrm{i}^2=-1$ and $\Real
+\mathrm{i}\pi=\{x+\mathrm{i}\pi:x\in \Real\}$. {Note that this approach is not an extension to the complex numbers. In fact, it designates a copy of the real numbers. Note also that in \cite{ir10}, the authors propose a closely related construction, though not a symmetric one, by creating a copy of the real numbers.  The formalism proposed in this paper is however convenient for introducing the following extended logarithmic function to $\widetilde {\mathbb M}$ and transferring the algebraic structure of the real set. }Let $\psi_{\ln}
:{\mathbb R} \longrightarrow \widetilde {\mathbb M}$ be the map defined by:
\begin{equation}
\psi_{ \ln}(x)=\left\{\begin{matrix}\ln(x)&\text{ if
}x>0\\-\infty&\text{ if }x=0\\\ln(-x)+\mathrm{i}\pi&\text{ if
}x<0.
\end{matrix}\right.
\end{equation}
The map $x\mapsto \psi_{ \ln}(x)$  is an isomorphism from $\mathbb
\Real$ to $\widetilde {\mathbb M}$. Let $\psi_{ \exp}(x): \widetilde
{\mathbb M}\longrightarrow \mathbb \Real$ be its inverse. Notice that
$\psi_{ \ln}(-1)=\mathrm{i}\pi$. By
definition we have $z\widetilde{\boxplus}
u=\psi_{\ln}\big(\psi_{\exp}(z)\boxplus \psi_{\exp}(u)\big )$. Moreover, we have $z\widetilde{\otimes}
u=\psi_{\ln}\big(\psi_{\exp}(z)\otimes \psi_{\exp}(u)\big )$. It follows that $(\widetilde{\mathbb M}, \tilde \oplus, \tilde \otimes)$ is isomorphic to $(\Real, \boxplus, \cdot)$ and is also an idempotent pseudo-field with the neutral elements $0_{\tilde {\mathbb M}}=-\infty$ and $1_{\tilde {\mathbb M}}=0$. 

{In Chapter 3 of \cite{koloma} it was shown that Max-Plus algebra can also be viewed as a limit algebraic structure via 
 the dequantization principle. }
For all $z\in \widetilde{\mathbb M}^n$, let  us  denote:
 \begin{equation}\widetilde{\bigboxplus _{i\in [n]}}z_i=
\psi_{\ln}\big(\bigboxplus_{i\in [n]}\psi_{\exp}(z_i)\big ).\end{equation}
From \cite{b20}, it follows that for all $z\in \widetilde{\mathbb M}^n$, we have
$$\widetilde{\bigboxplus_{i\in [n]}}z_i=\lim_{p\longrightarrow \infty}\frac{1}{2p+1}\psi_{\ln}\Big(\sum_{i\in [n]}\psi_{\exp}\big((2p+1)z_i\big)\Big). $$

\end{expl}

In the following, we give another example arising in a  Min -Times context (see  \cite{adilYe} for more details. ) 

\begin{expl} Let us define $\Real^{-1}:=\Real\backslash\{0\}\cup \{\infty\}$.  Let us consider the binary operation $\boxplus^{-1}: \mathbb \Real^{-1}\times \Real^{-1}\longrightarrow \Real^{-1}$ defined as:
\begin{equation}
x\boxplus^{-1}\! y=\left\{\begin{matrix}
x& \text{if} & |x|<|y|\\
y& \text{if} & |x|>|y|\\
+\infty& \text{if} & x+y=0\\
x& \text{if} &  x = y .  
\end{matrix}\right. 
\end{equation}
This operation is idempotent,  has $+\infty$ as   neutral element  and is symmetric. Moreover, {t}he standard product of non-zero numbers in $\Real^{-1}$ satisfies the desired condition. 
\end{expl}

\begin{rem}In \textnormal{\cite{ov10}}, an idempotent algebraic structure was derived from a suitable dequantization principle defined over the complex scalar field. Let us consider the binary operation $\boxplus_\mathbb C: \mathbb C\times \mathbb C\longrightarrow \mathbb C$ defined as:
\begin{equation}
z\boxplus_\mathbb C w=\left\{\begin{matrix}
z& \text{if} & |z|>|w|\\
w& \text{if} & |z|<|w|\\
0& \text{if} & z+w=0\\
\rho \frac{z+w}{|z+w|}& \text{if} & |z|=|w|=\rho>0.  
\end{matrix}\right. 
\end{equation}
It is easy to check that this operation is idempotent,  has $0$ as   neutral element  and is symmetric. Moreover, the standard product of complex numbers fulfils all the desired conditions. In addition note that this product is distributive over the operation $ \boxplus_\mathbb C$. In particular, if $|z|=|w|=\rho>0$ then we have for all $u\in \mathbb C$:
$$u\cdot (z\boxplus_\mathbb C w)=u \rho \frac{z+w}{|z+w|}=\rho \frac{uz+uw}{|z+w|}
=\rho |u|\frac{uz+uw}{|uz+uw|}=(uz)\boxplus_{\mathbb C} (uw).$$
However, it does not satisfy the weakened form of associativity of an idempotent pseudo-field. \end{rem}

\begin{defn}\label{SymIdemSpace}A triple $(X,\boxplus_X,\cdot_K)$ is an idempotent symmetric space defined over a  idempotent pseudo-field $K$
  if:

\noindent $(a)$ There exits a binary operation $\boxplus_X:X\times
X\longrightarrow X$ that satisfies the following properties: $(i)$ for all $x\in X$ we have $x\boxplus_X x=x$;  $(ii)$ for all $(x,y)\in X\times X$ we have $x\boxplus_X y=y\boxplus_X x$;  $(iii)$ there exists a neutral element $0_X\in X$ such that $x\boxplus_X 0_X=0_X\boxplus_X x$ for all $x\in X$; 
 $(iv)$ every $x\in X$ has a symmetric  element $-x$ such  that $(-x)\boxplus_X x=x\boxplus_X (-x)=0_X$.

\noindent $(b)$ There exists a scalar multiplication  $\cdot_K :\Real\times
X\longrightarrow X$ that satisfies the following properties: $(i)$ for all $\lambda \in \Real$ and all $(x,y)\in X\times X$, we
have $\lambda \cdot_K (x\boxplus_X y)=(\lambda \cdot_K x)\boxplus_X
(\lambda \cdot_K y)$; $(ii)$ for all $\lambda,\mu\in \Real\times \Real$ we have $(\lambda
\boxplus_K  \mu)\cdot_K x=(\lambda \cdot_K x)\boxplus_X (\mu \cdot_K
x)$, $(iii)$ for all $\lambda,\mu\in \Real\times \Real$ we have $(\lambda
\cdot_K\mu)\cdot_K x=\lambda \cdot_K(\mu \cdot_K x)$;  $(iv)$ for all $x\in X$ we have $1\cdot_K x=x$ and $-x=(-1)\cdot_K x$.

\end{defn}

{Suppose that $X$ is an idempotent symmetric space.} A  subset $Y$ of $X$ is an {\bf  idempotent symmetric subspace} if for all $x,y\in Y$ and all $t\in \Real$, $x\boxplus_X ty\in Y$.
Notice that, in this definition, 
associativity has been replaced with idempotence. 
Clearly, $(\mathbb{R}, \boxplus, \cdot)$ is an idempotent symmetric space.

\begin{lem} Suppose that $(X,\boxplus_X,\cdot_K)$ is an idempotent symmetric space defined over an idempotent pseudo-field $K$. Then, for all $x\in X$, one has:

$(a)$ $\lambda \cdot_K  0_X=0_X$

$(b)$ $ 0_K  \cdot_K x=0_X$

$(c)$ $\big((-y)\boxplus_X (-x)\big)\boxplus_X (x\boxplus_X
y)=0_X$.

\end{lem}
{\bf Proof:} $(a)$ By definition, for all $x\in X$ and all $\lambda\in {K}$, $\lambda \cdot_K 0_X=\lambda \cdot_K
\big(x\boxplus_X(-x)\big)$ for all $x\in X$. Using axioms $(b_1)$,
$(b_2)$ and  $(b_4)$ yields $\lambda \cdot_K 0_X=(\lambda \cdot_K
x)\boxplus_X (\lambda \cdot_K (-x))=(\lambda \cdot_K x)\boxplus_X
\big(\lambda \cdot_K (-1)\cdot_K x)\big)=(\lambda \cdot_K x)\boxplus_X
\big( (-1)\cdot_K \lambda\cdot_K x)\big)=(\lambda \cdot_K x)\boxplus_X
\big( -( \lambda\cdot_K x)\big)=0_X$. $(b)$ We have for all  $x\in X$, $ 0_K  \cdot_K
x=\big(1_K\boxplus (-1_K)\big)x=(1_K\cdot_K x)\boxplus_X  \big((-1_K)\cdot_K
x\big)= x\boxplus_X  (-x)=0_X$. $(c)$ For all $x,y\in X$, $\big((-y)\boxplus_X
(-x)\big)\boxplus_X (x\boxplus_X y)=\big((-1_K)\cdot_K y\boxplus_X
(-1_K)\cdot x\big)\boxplus_X (x\boxplus_X y)= \big((-1_K)\cdot(
y\boxplus_X  x)\big)\boxplus_X (x\boxplus_X y)=-(
x\boxplus_X  y)\boxplus_X (x\boxplus_X y)=0_X$.$\Box$\\

 \begin{prop} Let $X$ be an idempotent symmetric space. The intersection of two  idempotent symmetric subspaces of $X$ is an idempotent symmetric subspace.
\end{prop}

\begin{expl}\label{cartes} $(\Real^n, \boxplus, \cdot)$ is an idempotent symmetric space. $Y=\{(x_1,x_2,x_3,x_4)\in \Real^4: x_1=x_2,  x_3=x_4 \}$ is an idempotent symmetric subspace of $\Real^4$. 
\end{expl}

{In the following, if \( K \) is an idempotent pseudo-field, we define the operations for all \( x, y \in K^n \) and all \( \lambda \in K \) as follows:  
\begin{equation}  
x \boxplus_K y = (x_1, \dots, x_n) \boxplus_K (y_1, \dots, y_n) = (x_1 \boxplus_K y_1, \dots, x_n \boxplus_K y_n)  
\end{equation}  
and  
\begin{equation}  
\lambda \cdot_K x = \lambda \cdot_K (x_1, \dots, x_n) = (\lambda \cdot_K x_1, \dots, \lambda \cdot_K x_n).  
\end{equation}  }
The two following examples are based on this construction.
\begin{expl}\label{cartesMaslov} $(\widetilde{\mathbb M}^n, \tilde \boxplus, \tilde \otimes)$ is an idempotent symmetric space. 
\end{expl}

\begin{expl} \label{cartesAdil} $(\Real^{-n},  \boxplus^{-1}, \cdot)$ is an idempotent symmetric space, with   $\Real^{-n}=\big(\Real^{-1}\big)^n$.  
\end{expl} 
In the next example, it is shown that a similar construction can be applied to   certain   infinite-dimensional spaces.
 
\begin{expl}Suppose that $I$ is a subset of $\Real$ and let $\lambda$ be the Lebesgue measure.
Let $L_p(\lambda)$ be the set of the real valued maps $f:
I\longrightarrow \Real$ such that
$$\int_I |f|^p d\lambda<+\infty$$
for $p\in ]1,+\infty[$. Let us consider the binary
operation defined over $L_p(\lambda)$ by:
\begin{equation}f\boxplus g:x\mapsto f(x)\boxplus g(x)=
\left\{\begin{matrix}f(x)\ &\hbox{ if } &|f(x)|&>&|f(y)|\\
\frac{1}{2}\big(f(x)+f(y)\big)&\hbox{ if }&|f(x)|&=&|f(y)|\\
 f(y)& \hbox{ if }&
|f(x)|&<&|f(y)|.\end{matrix}\right.\end{equation} Clearly, one has
$|f(x)\boxplus g(x)|\leq |f(x)|+| g(x)|$. It follows that

$$\int_I |f\boxplus g|^p d\lambda<\int_I |f|^p d\lambda +\int_I | g|^p d\lambda<\infty.$$

Let $\cdot$ be the standard product of a real valued
 function by a real scalar, it follows that $(L_p(\lambda),\boxplus,\cdot)$ is an idempotent symmetric space.
\end{expl}

\begin{rem} Notice that that Examples \ref{cartesMaslov} and \ref{cartesAdil} are constructed from an homeomorphic transformation of $(\Real^n, \boxplus,\cdot)$. It is therefore easy to define two convex structures on   $(\widetilde{\mathbb M}^n, \tilde \boxplus, \tilde \otimes)$ and $(\Real^{-n},  \boxplus^{-1}, \cdot)$ respectively. Namely a subset $M$ of $\widetilde{\mathbb M} $ is $\mathbb M$-convex if for all $x,y\in M$ and all $t\in \mathbb M$, $x\widetilde{\boxplus}ty\in C$. Paralleling \cite{adilYe}, a subset $C$ of $\Real^{-n}$ is inverse-$\mathbb B$-convex if for all $x,y\in C$ and all $t\in [1,+\infty]$, $x\boxplus^{-1}ty\in C.$\end{rem}
The next result specifically concerns the algebraic structure considered in \cite{b15}. However it has some immediate implications for the idempotent symmetric  spaces defined in examples \ref{cartesMaslov} and \ref{cartesAdil}. 

\begin{prop} \label{comb} For all subsets $Y$ of  $ \Real^n $, the following claims are equivalent:

$(a)$ $Y$ is an idempotent symmetric subspace of $(\Real^n, \boxplus,\cdot)$.
\medskip

$(b)$  For all $(x_1,...,x_m)\in Y^m$ and all $t\in \Real^m$, we have  $ \bigboxplus\limits_{i\in [m]}t_ix_i\in Y.$ 
\end{prop}
{\bf Proof:} Let us prove that $(a)$ implies $(b)$. If $Y$ is
an idempotent symmetric subspace, this property is true for $m=2$. Suppose
it is true at rank $m$ and let us prove that it is true at rank
$m+1$. In other words, assume that for all $(x_1,...,x_m)\in Y^m$ we
have: $ \bigboxplus_{i\in [m]}t_ix_i \in Y$ for all $t\in \Real^m$, we need to prove that if
$(x_1,...,x_m, x_{m+1})\in Y^{m+1}$ then for all $t\in \Real^{m+1}$ we have $\bigboxplus_{i\in
[m+1]}t_ix_i\in Y$. To establish this property, we use equation
\eqref{Decompos} which implies that if $(x_1,...,x_m, x_{m+1})\in
Y^{m+1}$ then, for all $t\in \Real^{m+1}$, we have
$$\bigboxplus_{i\in [m+1]}t_ix_i=
\bigboxplus_{i\in [m+1]}\Big[t_ix_i\boxplus \Big(\bigboxplus_{j\in
[m+1]\backslash \{i\}}t_jx_j\Big)\Big].\quad (\star)$$

  For all $i$, the  vectors $z_i=t_ix_i\boxplus \Big(\bigboxplus_{j\in
[m+1]\backslash \{i\}}t_jx_j\Big)$ are copositive. Thus the associativity holds true for the $m$-tuple $(z_1,...,z_m)$. However, since we assume that the property holds true at rank $m$, $z_i\in Y$ for all $i$. Thus we deduce that $\bigboxplus_{i\in [m+1]}t_ix_i\in Y$. Therefore, the first part of the proof is {established}.
   Since $Y$ is a subset of $\Real^n$, that is an idempotent symmetric space, the converse inclusion is immediate. $\Box$\\

The following property is an immediate consequence.  

\begin{prop} Let $Y$ be an idempotent symmetric subspace of $(\Real^n, \boxplus,\cdot)$. Then  $Y$ is {idempotent symmetric  convex}.
\end{prop}

\subsection{ Generalized Norms and Ultrametric in Limit.}

Let $\|\cdot \|$ be the Euclidean norm defined for all $x\in \Real^n$ as $\|x \|=\big(\sum_{i\in [n]}{x_i}^2\big)^\frac{1}{2}$. Let $d$ be the Euclidean distance defined for all $x,y\in \Real^n$ as $d(x,y)=\|x-y\|$. For all positive natural numbers $p$, we introduce the map $\|\cdot \|_{{\varphi_p}}: \Real^n\longrightarrow \Real_+$
defined as
\begin{equation}
\|x \|_{{\varphi_p}}=\varphi_p^{-1}\big(\|\phi_p(x) \|\big).
\end{equation}
The following properties are immediate: $\|x \|_{{\varphi_p}}=0\iff x=0$; $\|\alpha x  \|_{{\varphi_p}}=|\alpha|\| x  \|_{{\varphi_p}}$;
$\|  x  \stackrel{p}{+}y \|_{{\varphi_p}}\leq \|  x  \|_{{\varphi_p}}\stackrel{p}{+} \|   y \|_{{\varphi_p}}$.
We deduce the following statements:
\begin{prop} For all $x,y\in \Real^n$ and all $\alpha \in \Real$, we have the the following properties:\\
$(a)$ $\|x\|_\infty=\lim_{p\longrightarrow \infty}\|x \|_{{\varphi_p}}$\\
$(b)$ $\| x\boxplus y\|_\infty\leq  \max\{\|x \|_{\infty},\|y \|_{\infty}\}$. \\
\end{prop}
{\bf Proof:} $(a)$ For all $x\in \Real^n$, we have
$$\|x \|_{{\varphi_p}}=\Big(\sum_{i\in [n]}|x_i|^{2(2p+1)}\Big)^{\frac{1}{2(2p+1)}}.$$
Taking the limit yields $(a)$. $(b)$ For all natural numbers $p$, we have:
$$\|  x  \stackrel{p}{+}y \|_{{\varphi_p}}\leq \|  x  \|_{{\varphi_p}}\stackrel{p}{+} \|   y \|_{{\varphi_p}}.\quad (A) $$
Suppose that $\{z^{(p)}\}_{p\in \mathbb N}$ is a sequence of $\Real_+^n$ which converges to some $z\in \Real_+^n$. It follows that we have $\lim_{p\longrightarrow \infty}\|z^{(p)} \|_{{\varphi_p}}=\|z \|_\infty$. For all $x\in \Real^n$,
let us denote $|x|=(|x_1|,...,|x_n|)$. By definition, we have
$$\|  x  \stackrel{p}{+}y \|_{{\varphi_p}}=\|  |x  \stackrel{p}{+}y| \|_{{\varphi_p}}.$$
Moreover, $\lim_{p\longrightarrow \infty}|x  \stackrel{p}{+}y|=|x\boxplus y|$. It follows that:
$$\lim_{p\longrightarrow \infty}\|  x  \stackrel{p}{+}y \|_{{\varphi_p}}=\||x\boxplus y|\|_\infty=\| x\boxplus y \|_\infty.$$

Since $\|  x  \|_{{\varphi_p}}\geq 0$ and
$\|  y  \|_{{\varphi_p}}\geq 0$, taking the limit on both sides in $(A)$ yields the result. $\Box$\\

Recall that, for all $\alpha, \beta\in \Real$, we adopt the notation $\alpha\boxminus \beta=\alpha\boxplus (-\beta)$. Let us consider the map $d_\boxplus: \Real^n\times \Real^n\rightarrow \Real_+$ defined as:
\begin{equation}d_\boxplus(x,y)=\max_{i\in [n]}|x_i\boxminus y_i|\end{equation}
In addition, let us consider the map $d_{{\varphi_p}}:\Real^n\times \Real^n\longrightarrow \Real_+$ defined by:
\begin{equation}d_{{\varphi_p}}(x,y)=\|x \stackrel{p}{-}y\|_{{\varphi_p}}.\end{equation}
This function is a distance relative to the algebraic structure obtained from the binary operation $\stackrel{p}{+}$. It defines a metric on any $\varphi_p$-vector space. Note that $d_{{\varphi_p}}$  is also a distance in the usual sense because as $p\geq 1$, we have for all $z\in \Real^n$ the relation:
\begin{equation}d_{{\varphi_p}}(x,y)\leq d_{{\varphi_p}}(x,z)\stackrel{p}{+}d_{{\varphi_p}}(z,y)\leq d_{{\varphi_p}}(x,z)+d_{{\varphi_p}}(z,y).
\end{equation}
\begin{lem}For all $x,y\in \Real^n$, we have
$$\lim_{p\longrightarrow\infty}d_{{\varphi_p}}(x,y)=d_\boxplus(x,y).$$
\end{lem}
{\bf Proof:}
For all $x\in \Real^n$, we have
$$\|x\stackrel{p}{-}y \|_{{\varphi_p}}=\Big(\sum_{i\in [n]}\Big(\Big|\big({x_i}^{2p+1}-{y_i}^{2p+1}\big)^{\frac{1}{2p+1}}\Big|^{2}\Big)^{2p+1}\Big)^{\frac{1}{2(2p+1)}}.$$
For all $i\in [n]$ and any natural number  $p$, set $z_i^{(p)}= \big|\big({x_i}^{2p+1}-{y_i}^{2p+1}\big)^{\frac{1}{2p+1}} \big|^2$. For all $i$, we have:
$\lim_{p\longrightarrow \infty}z_i^{(p)}=|x_i\boxminus y_i|^2 $. Since for all $p$, $z_i^{(p)}\geq 0$, it follows that
$$\lim_{p\longrightarrow \infty}\|x\stackrel{p}{-}y \|_{{\varphi_p}}=\Big(\max_{i\in [n]}  |x_i\boxminus y_i |^{2} \Big)^{\frac{1}{2}}=d_\boxplus (x,y). \quad \Box$$

\begin{prop}\label{ultramprop}The map $d_\boxplus: \Real^n\times \Real^n\longrightarrow \Real_+$ defines an {ultrametric} on $\Real^n$. For all $x,y,z\in \Real^n$:\\ $(1)$ $d_\boxplus(x,y)=d_\boxplus(y,x)$; \\$(2)$ $d_\boxplus(x,y)=0$ if and only if $x=y$;\\ $(3)$
 $d_\boxplus(x,y)\leq \max\{d_\boxplus(x,z),d_\boxplus(z,y)\}. $
\end{prop}
{\bf Proof:} $(1)$ is immediate. $(2)$ follows from the fact that for all $i$ we have $x_i\boxminus y_i=0\iff x_i=-y_i$. $(3)$ For all $x,y,z\in \Real^n$ and for all natural numbers $p$, we have
$$d_{{\varphi_p}}(x,y)\leq d_{{\varphi_p}}(x,z)\stackrel{p}{+}d_{{\varphi_p}}(z,y).$$
Since for all $x,y\in \Real^n$ we have $d_{{\varphi_p}}(x,y)\geq 0$ and $\lim_{p\longrightarrow \infty}d_{{\varphi_p}}(x,y)=d_{\boxplus}(x,y)$,  we deduce the result. $\Box$\\

Note that any ultrametric defines a metric since the inequality $d(x,y)\leq \max\{d(x,z), d(y,z)\}$ implies the triangular inequality $d(x,y)\leq  d(x,z)+ d(y,z)$. In particular  for all $x,y,z\in E$, at least one of the three equalities $  d(x,y)=d(y,z) $ or $  d(x,z)=d(y,z) $ or $  d(x,y)=d(z,x) $ holds.  From the above definition, one can conclude several typical properties of ultrametrics.

\begin{prop}\label{UMSpace}Let $Y$ be an idempotent symmetric subspace of $\Real^n$. $(Y,d_\boxplus)$ is an ultrametric space. 
\end{prop}
{\bf Proof:} Since \( Y \) is an idempotent symmetric space, the metric \( d_\boxplus \) is well-defined on \( Y \).
 From Proposition \ref{ultramprop}, the result is immediate. $\Box$\\

In the following we say that a sequence $\{x^{(k)}\}_{k\in \mathbb N}$ is {\bf $\digamma$-convergent }to $x\in \Real^n$ if
$\lim_{k\longrightarrow \infty} d_\boxplus (x^{(k)},x)=0$.

\begin{lem}Suppose that $\{x^{(k)}\}_{k\in \mathbb N}$ is $\digamma$-convergent to $x\in \Real^n$. Let $I(x)=\{i\in [n]: x_i\not=0\}$. Then, there exists some natural number $k_0$ such that for all $k\geq k_0$ and for all $i\in I(x)$ we have $x_i^{(k)}=x_i$. Moreover, for all
$i\notin I(x)$, $\lim_{k\longrightarrow \infty}x_i^{(k)}=0$.
\end{lem}
{\bf Proof:} Suppose that $\lim_{k\longrightarrow \infty}d_\boxplus (x_i^{(k)},x)=0$. Assume that there exits an increasing subsequence $\{k_{q}\}_{q\in \mathbb N}$ and some $i_0\in I(x)$ such that, for all $q$, $x_i^{(k_q)}\not=x_i$  and let us show a contradiction. In such a case, for all $q\in \mathbb N$, we have $|x_{i_0}^{(k_q)}\boxminus x_{i_0}|=\max \{|x_{i_0}^{(k_q)}|, |x_{i_0}|\}\not=0$. This implies that $$d_\boxplus (x^{(k_q)},x)= 
\max_{i\in I(x)}  \{|x_i^{(k_q)}\boxminus x_i| \}\geq |x_{i_0}^{(k_q)}\boxminus x_{i_0}|=\max \{|x_{i_0}^{(k_q)}|, |x_{i_0}|\}>0,$$ that is a contradiction. Hence, there is some $k_0$ such that for all $k\geq k_0$ and for all $i\in I(x)$ we have $x_i^{(k)}=x_i$. If $i\notin I(x)$ then $|x_i^{(k )}\boxminus x_i|=|x_i^{(k )}\boxminus 0|=|x_i^{(k )}|$ and $\lim_{k\longrightarrow \infty} d_\boxplus (x^{(k)},x)=0$. Therefore, $\lim_{k\longrightarrow \infty}x_i^{(k)}=0$. $\Box$ \\

It is now interesting to see that the idempotent symmetrized space $(\widetilde {\mathbb M}^n, \tilde \boxplus, \tilde \otimes)$ can also be endowed with a suitable form of ultrametric. Let $d_{\widetilde \boxplus}: \widetilde{\mathbb M}^n\times \widetilde{\mathbb M}^n\longrightarrow \Real \cup \{-\infty\}$ defined as:
\begin{equation}
d_{\widetilde \boxplus}(z,w)=\widetilde{\bigboxplus_{i\in [n]}}\big |z_i\widetilde{\boxminus}w_i\big |_{\widetilde{\mathbb M}},
\end{equation}
where for all $\gamma\in \widetilde {\mathbb M}$, $|\gamma \big |_{\widetilde{\mathbb M}}=\psi_{\ln}\big(|\psi_{\exp}(\gamma)|\big)$ and, for all $i$,
$z_i \widetilde{\boxminus}w_i =\psi_{\ln}\Big(\psi_{\exp}(z_i)\boxminus \psi_{\exp}(w_i)\Big)$. It is therefore easy to check that: $(i)$ for all $(z,w)\in \widetilde {\mathbb M}^n\times \widetilde {\mathbb M}^n$, $d_{\widetilde \boxplus}(z,w)=-\infty\iff z=w$; $(ii)$ for all $ z,t,w  \in \widetilde {\mathbb M}^n$,  $d_{\widetilde \boxplus}(z,w)\leq \max\{d_{\widetilde \boxplus}(z,t), d_{\widetilde \boxplus}(t,w)\}$. Therefore, paralleling Proposition \ref{UMSpace}, $(\widetilde {\mathbb M}, \tilde \boxplus, \tilde \otimes)$ can also be endowed with a topological structure. In particular, we can construct a standard ultrametric by defining for all $z,w\in \widetilde{\mathbb M}^n$ the distance $\bar d_{\widetilde{\boxplus}}(z,w)=\psi_{\exp}\Big(d_{\widetilde \boxplus}(z,w)\Big)$. 

\begin{prop}Let $Z$ be an idempotent symmetric subspace of $\widetilde{\mathbb M}^n$. Then $(Z,\bar d_{\widetilde{\boxplus}})$ is an ultrametric space. 

\end{prop}

\subsection{Geometry of the  Ultrametric  Ball}

In the following, let \( B_\boxplus(x, r) = \{z \in \mathbb{R}^n : d_\boxplus(x, z) {<} r\} \) and \( B_\boxplus(x, r] = \{z \in \mathbb{R}^n : d_\boxplus(x, z) \leq r\} \) respectively denote the open and closed ball{s} centered at \( x \) with radius \( r \).

  For all
${x}\in \mathbb R^n$ and all subsets $I$ of $[n]$, let us  consider
the map $\mathcal A_x:I \longrightarrow \mathbb R_+$ defined for all
$i \in I$ by $\mathcal A_x(i)= |x_i|.$
This map associates to each index $i$ the absolute value of the $i$-th coordinate of $x$. 
For all $\alpha \in \Real_+ $, $\mathcal A_x^{-1}(\alpha)=\{i\in I: |x_i|=\alpha\}$.

For all $x,y\in \mathbb R^n$, let $\mathcal L(x,y)$ be a subset of $I$
defined by
\begin{equation}\mathcal L(x,y)=\Big\{i\in [n]: x_i\not=y_i\Big\}.\end{equation}
$\mathcal L (x,y)$ is  obtained by dropping from $[n]$ all the $i$'s such that
$  x_i=y_i$.
It follows that {$d_\boxplus(x,y)=0$ if $x=y$ and }
\begin{equation}d_\boxplus(x,y)=\left\{\begin{matrix}\max\limits_{i\in \mathcal L(x,y)}|x_{i}| &\text{ if }& \max\limits_{i\in\mathcal L(x,y)}|x_i| \geq  \max\limits_{i\in\mathcal L(x,y)}|y_i|\\
\max\limits_{i\in \mathcal L(x,y)}|y_{i}| &\text{ if }& \max\limits_{i\in\mathcal L(x,y)}|x_i| \leq  \max\limits_{i\in\mathcal L(x,y)}|y_i|,
\end{matrix}\right.\end{equation}
if $x\not=y$.    In the next statement we describe the geometry of any ball centered at $x\in \Real^n$.

\begin{prop}\label{GeomBall} Let $x\in \Real^n$. Suppose that $ \mathcal A_x(\,[n]\,)=\bigcup_{i\in [n]} \mathcal A_x(\,i\,)=\{\alpha_j\}_{j\in [m]}$ with
 $\alpha_j<\alpha_{j+1}$ for all $j\in [m-1]$ and $\alpha_j\geq 0$ for all $j$.
 
 $(i)$ If $x=0$ and $ \mathcal A_x(\,[n]\,)=\{0\}$, then for all $\alpha\geq 0$, $B_\boxplus(0, \alpha]=B_\infty(0,\alpha]$.
 
 $(ii)$ If either $\alpha=0$ or $\alpha \in [0, \alpha_{1}[$ then $B_\boxplus(x, \alpha]=\{x\}$.
 
 $(iii)$ If $\alpha \in [\alpha_j, \alpha_{j+1}[$ then $z\in B_\boxplus(x, \alpha]$ if and only if:
 $$ |z_k|\leq \alpha, \text{ for all }k\in \mathcal A_x^{-1}\big([0, \alpha_j]\big)\text{ and } z_k=x_k \text{ for all }   k\in  \mathcal A_x^{-1}\big([  \alpha_{j+1}, \alpha_m]\big). 
$$

 $(iv)$ If $\alpha \geq \alpha_m= \|x\|_\infty$ then $B_\boxplus (x, \alpha]=B_\infty(0, \alpha].$

\end{prop}
{\bf Proof: } $(i)$ is immediate. $(ii)$ If $\alpha=0$, then $B_\boxplus(x, 0\,]=\{x\}$. If $\alpha_1=0$, the case $\alpha< \alpha_1$ is then excluded.
Suppose that $\alpha_1>0.$ For all $z\not=x$, $\max_{i\in [n]}|x_i\boxplus (-z_i)|\geq \max_{i\in [n]}|x_i|\geq \alpha_1$. Therefore $\alpha<\alpha_1$ implies that 
$B_\boxplus(x, \alpha\,]=\{x\}$. $(iii)$ Suppose that $\alpha \in [\alpha_1, \alpha_2]$. If there is some $i\in \mathcal A_x^{-1}\big([  \alpha_{2}, \alpha_m]\big)$ such that $z_i\not=x_i$, then $d_\boxplus(x,z)\geq \alpha_2>\alpha_1$. Therefore, if $z\in B_\boxplus(x, 0\,]$ we should have $z_i=x_i$ for all $i\in \mathcal A_x^{-1}\big([  \alpha_{2}, \alpha_m]\big)$. Moreover, for all $i\in  \mathcal A_x^{-1}\big(\{\alpha_1\}\big)$ we should have $|z_i|\leq \alpha.$ Let us extend this result {to some} arbitrary $j\leq m-1$. Suppose that $\alpha \in [\alpha_j, \alpha_{j+1}]$. If there is some $i\in \mathcal A_x^{-1}\big([  \alpha_{j+1}, \alpha_m]\big)$ such that $z_i\not=x_i$, then $d_\boxplus(x,z)\geq \alpha_{j+1}>\alpha_j$. Therefore, if $z\in B_\boxplus(x,\alpha]$, we should have $z_i=x_i$ for all $i\in \mathcal A_x^{-1}\big([  \alpha_{j+1}, \alpha_m]\big)$. Moreover, for all $i \in \mathcal A_x^{-1}\big([  \alpha_{1}, \alpha_j]\big)$ we should have $|z_i|\leq \alpha.$ $(iv)$ Suppose now that $\alpha\geq \|x\|_\infty$. By construction $\alpha_m=\|x\|_\infty$ and $z\in B_\boxplus(x,\alpha]$ if and only if $|z_i|\leq \|x\|_\infty$, which ends the proof. $\Box$ 

\begin{cor} For all $\alpha\geq 0$, and all $x\in \Real^n$, $B_\boxplus(x,\alpha]$ is closed with respect to the norm topology. 
\end{cor}
{\bf Proof:} Let $ \mathcal A_x(\,[n]\,)=\{\alpha_j\}_{j\in [m]}$ with
 $\alpha_j<\alpha_{j+1}$ for all $j\in [m-1]$ and $\alpha_j\geq 0$ for all $j$.
The result is immediate if either $\alpha=0$ or $\alpha \in [0, \alpha_{1}[$ since $B_\boxplus(x, \alpha]=\{x\}$. 
From Proposition \ref{GeomBall}, if $\alpha \in [\alpha_j, \alpha_{j+1}[$, then
 $$B_\boxplus(x, \alpha]=\Big\{z\in \Real^n\!\!\!: |z_k|\leq \alpha, k\in \mathcal A_x^{-1}\big([0, \alpha_j]\big)\Big\}\cap \Big\{ z\in \Real^n\!\!\!: z_k=x_k,   k\in  \mathcal A_x^{-1}\big([  \alpha_{j+1}, \alpha_m]\big)\Big\}.$$ 
 Since the intersection of two closed sets is a closed set, we deduce that $B_\boxplus(x, \alpha]$ is closed. 
 If $\alpha\geq \|x\|_\infty$, the result is immediate. $\Box$\\

\begin{center}{\scriptsize 
% This is a LaTeX picture output by TeXCAD.
% File name: [Clipboard].
% Version of TeXCAD: 4.51
% Reference / build: 27-Nov-2018 (rev. a75)
% For new versions, check: http://texcad.sf.net/
% Options on the following lines.
%\grade{\on}
%\emlines{\off}
%\epic{\off}
%\beziermacro{\on}
%\reduce{\on}
%\snapping{\off}
%\pvinsert{% Your \input, \def, etc. here}
%\quality{8.000}
%\graddiff{0.005}
%\snapasp{1}
%\zoom{4.0000}

% This is a LaTeX picture output by TeXCAD.
% File name: [Clipboard].
% Version of TeXCAD: 4.51
% Reference / build: 27-Nov-2018 (rev. a75)
% For new versions, check: http://texcad.sf.net/
% Options on the following lines.
%\grade{\on}
%\emlines{\off}
%\epic{\off}
%\beziermacro{\on}
%\reduce{\on}
%\snapping{\off}
%\pvinsert{% Your \input, \def, etc. here}
%\quality{8.000}
%\graddiff{0.005}
%\snapasp{1}
%\zoom{4.0000}
\unitlength 0.35mm % = 1.138pt
\linethickness{0.4pt}
\ifx\plotpoint\undefined\newsavebox{\plotpoint}\fi % GNUPLOT compatibility
\begin{picture}(325.01,213.35)(0,0)
\put(.25,168.65){\vector(1,0){145.2}}
\put(0,55.9){\vector(1,0){145.2}}
\put(172.25,169.05){\vector(1,0){145.2}}
\put(172,56.3){\vector(1,0){145.2}}
\put(61.11,14.8){\vector(0,1){81.3}}
\put(231.26,128.85){\vector(0,1){81.3}}
\put(58.01,128.85){\vector(0,1){81.3}}
\put(231.01,16.1){\vector(0,1){81.3}}
\put(107.11,204.35){\circle*{1.749}}
\put(107.61,102.85){\circle*{1.749}}
\put(282.431,108.611){\circle*{1.911}}
\put(15.61,102.85){\circle*{1.749}}
\put(181.921,108.611){\circle*{1.911}}
\put(15.36,11.35){\circle*{1.749}}
\put(181.648,7.897){\circle*{1.911}}
\put(107.36,10.85){\circle*{1.749}}
\put(282.158,8.101){\circle*{1.911}}
\put(107.11,88.35){\circle*{1.749}}
\put(274.61,88.1){\circle*{1.749}}
\put(277.01,205.65){\circle*{1.749}}
\put(277.01,136.65){\circle*{1.749}}
\put(55.41,49.45){\makebox(0,0)[cc]{$0$}}
\put(226.91,163.65){\makebox(0,0)[cc]{$0$}}
\put(226.66,50.9){\makebox(0,0)[cc]{$0$}}
\put(325.01,168.15){\makebox(0,0)[cc]{$x_1$}}
\put(324.76,55.4){\makebox(0,0)[cc]{$x_1$}}
\put(154.66,166.85){\makebox(0,0)[cc]{$x_1$}}
\put(154.41,54.1){\makebox(0,0)[cc]{$x_1$}}
\put(231.11,213){\makebox(0,0)[cc]{$x_2$}}
\put(57.86,213){\makebox(0,0)[cc]{$x_2$}}
\put(230.86,100.25){\makebox(0,0)[cc]{$x_2$}}
\put(60.51,98.95){\makebox(0,0)[cc]{$x_2$}}
\put(61.96,117){\makebox(0,0)[cc]{$0\leq \alpha<2$}}
\put(61.21,0){\makebox(0,0)[cc]{$ \alpha=3$}}
\put(238.61,118.3){\makebox(0,0)[cc]{$2\leq \alpha<3$}}
\put(237.86,1.3){\makebox(0,0)[cc]{$ \alpha>3$}}
\put(277.01,205.65){\line(0,-1){68.55}}
\put(107.26,102.369){\line(0,-1){91.172}}
\put(282.049,108.086){\line(0,-1){99.605}}
\put(15.51,101.869){\line(0,-1){91.172}}
\put(181.812,107.539){\line(0,-1){99.605}}
\put(104.16,164.6){\makebox(0,0)[cc]{$3$}}
\put(273.21,164.9){\makebox(0,0)[cc]{$3$}}
\put(273.96,52.25){\makebox(0,0)[cc]{$3$}}
\put(103.21,52.5){\makebox(0,0)[cc]{$3$}}
\put(53.248,204.2){\makebox(0,0)[cc]{$2$}}
\put(54.96,88.45){\makebox(0,0)[cc]{$2$}}
\put(225.26,205.5){\makebox(0,0)[cc]{$2$}}
\put(225.01,92.75){\makebox(0,0)[cc]{$2$}}
\put(116.76,213.35){\makebox(0,0)[cc]{$B_\boxplus\big((3,2), \alpha\big]$}}
\put(132.26,62.1){\makebox(0,0)[cc]{$B_\boxplus\big((3,2), \alpha\big]$}}
\put(310.81,177.5){\makebox(0,0)[cc]{$B_\boxplus\big((3,2), \alpha\big]$}}
\put(315.81,97.25){\makebox(0,0)[cc]{$B_\boxplus\big((3,2), \alpha\big]$}}
\put(107.61,102.867){\line(-1,0){92}}
\put(282.431,108.629){\line(-1,0){100.51}}
\put(107.11,10.432){\line(-1,0){92}}
\put(281.885,7.644){\line(-1,0){100.51}}
\put(160,0){\makebox(0,0)[cc]{{\bf Fig. 3:} Ultrametric Ball.}}
\end{picture}
}

\end{center}

\begin{cor} For all $\alpha\geq 0$, and all $x\in \Real^n$, $B_\boxplus(x,\alpha]$ is {idempotent symmetric  convex}.
\end{cor}
{\bf Proof:} From Proposition \ref{GeomBall}, $B_\boxplus(x,\alpha]$ is the cartesian product between a box and a singleton. Since any box is {idempotent symmetric  convex}, the result immediately follows. $\Box$\\

In \cite{b15} the following
formulation of a $\mathbb B$-polytope was established in the case of two points.\footnote{We use the convention that for all $(\alpha,\beta,\gamma,\delta)\in \mathbb R^4$, $ \alpha\boxplus \beta\boxplus \gamma\boxplus \delta=\digamma_{[4]}(\alpha,\beta,\gamma,\delta)$.} For all $x,y\in \Real^n$,
\begin{equation}\label{fond2}Co^{\infty}(x,y)=\Big \{t x\boxplus r x\boxplus s y\boxplus w y:
\max\{t,r,s,w\}=1,t,r,s,w\geq 0\Big\}.\end{equation}
In classical convexity, counting the same point twice in a finite set does not alter its convex hull. However, this is not the case here due to the specific properties of the non-associative operation we consider.
This formulation is useful to prove the next decomposition result  over $(\Real^n, \boxplus, \cdot)$.

\begin{prop}For all $x,y\in \Real^n$, and all $z\in Co^\infty(x,y)$ we have:
$$d_{\boxplus}(x,y)=\max\{d_{\boxplus}(x,z),d_{\boxplus}(z,y)\}.$$
\end{prop}
{\bf Proof:} If $x=y$, then the result is immediate. Let us assume that $x\not=y$. Suppose that $z\in  Co^\infty(x,y)$, in such a case there are some $r,s,t,w\in [0,1]$ with $\max\{r,s,t,w\}=1$ such that
$$z=rx\boxplus sx\boxplus ty\boxplus wy.$$
Suppose that $\alpha=d_\boxplus (x,y)$. Let $\mathcal L(x,y)=\{i\in [n]: x_i\not=y_i\}$. By construction
$$d_\boxplus (x,y)=\max\{|x_i|, |y_i|:i\in \mathcal L(x,y)\}.$$
Therefore, for all $i\notin \mathcal L(x,y)$, if $|x_i|>d_\boxplus (x,y)$ then $x_i=y_i$. Thus, we deduce that $z_i=rx_i\boxplus sx_i\boxplus ty_i\boxplus wy_i=x_i=y_i$. It follows that
$x_i\boxminus z_i=y_i\boxminus z_i=0$. Hence, we have:
$$d_\boxplus (x,z)=\max \{|x_i\boxminus z_i|:i\in \mathcal L(x,y) \}$$
and
$$d_\boxplus (y,z)=\max \{|y_i\boxminus z_i|:i\in \mathcal L(x,y) \}. $$
Now, notice that for all $i$
 $$|z_i|=|rx_i\boxplus sx_i\boxplus ty_i\boxplus wy_i|\leq   \max\{|x_i|, |y_i|\}$$
Since, for all $i$, $|x_i\boxminus z_i|\leq   \max\{|x_i|, |y_i|\}$ and $|y_i\boxminus z_i|\leq   \max\{|x_i|, |y_i|\}$, we deduce that
$$d_\boxplus (x,z)\leq \max \{|x_i|, |y_i|:i\in \mathcal L(x,y) \}=d_\boxplus (x,y)$$
and
$$d_\boxplus (y,z)\leq \max \{|x_i|, |y_i|:i\in \mathcal L(x,y) \}=d_\boxplus (x,y).$$
Thus
$$\max\{d_{\boxplus}(x,z),d_{\boxplus}(z,y)\}\leq d_\boxplus (x,y).$$
However, from the triangular inequality, we have
$$d_{\boxplus}(x,y)\leq \max\{d_{\boxplus}(x,z),d_{\boxplus}(z,y)\}.$$
Since the converse inequality holds, this   completes the proof. $\Box$\\

The next result is an immediate consequence.

\begin{prop}\label{multidec}Let $A=\{x^{(1)},...,x^{(m)}\}$ be a finite subset of $\Real^n$. Suppose moreover that  $(i)$ $ Co^\infty \big(x^{(1)},x^{(m)}\big)=\bigcup_{i\in [m-1]} Co^\infty \big(x^{(i)},x^{(i+1)}\big)$; $(ii)$ for all $i$ $Co^\infty \big(x^{(i-1)},x^{(i)}\big)\cap Co^\infty \big(x^{(i)},x^{(i+1)}\big)=\{x^{(i)}\}$. Then we have:
$$d_{\boxplus}\big(x^{(1)},x^{(m)}\big)=\max_{i\in [m-1]}d_\boxplus \big(x^{(i)},x^{(i+1)}\big).$$
\end{prop}
{\bf Proof:}
Since $x^{(2)}\in Co^\infty \big(x^{(1)},x^{(m)}\big)$, we have:
$$d_{\boxplus}\big(x^{(1)},x^{(m)}\big)=\max\big \{d_{\boxplus}\big(x^{(1)},x^{(2)}\big), d_{\boxplus}\big( x^{(2)},x^{(m)})\big )\big \}.$$
From the decomposition property established in \cite{b19}, {we have} $Co^\infty \big(x^{(k)},x^{(m)})=Co^\infty \big(x^{(k)},x^{(k+1)})\cup Co^\infty \big(x^{(k+1)},x^{(m)})$. It follows that for any $k\in [m-1]$, $x^{(k+1)}\in Co^\infty \big(x^{(k)},x^{(m)}\big)$. Hence:

$$d_\boxplus \big(x^{(k)},x^{(m)}\big)=\max \Big \{ d_\boxplus\big(x^{(k)},x^{(k+1)}\big) , d_\boxplus\big (x^{(k+1)},x^{(m)}\big)\Big\}.$$
By recurrence the proof is then immediate. $\Box$\\

\begin{prop}Let $x,y\in \Real^n$. If $x=-y$, then for all  $z\in Co^\infty (x,y)\backslash \{x,y\}$, $d_{\boxplus}(x,y)= d_{\boxplus}(x,z)=d_{\boxplus}(z,y).$

\end{prop}
{\bf Proof:} Suppose that $z\in Co^\infty \big(x,y)$, in such a case there are some $r,s,t,w\in [0,1]$ with $\max\{r,s,t,w\}=1$ such that
 $z=rx\boxplus sx\boxplus ty\boxplus wy.$ Since $z\notin\{x,y\}$ we have  $\max\{r,s\}=\max\{t,w\}=1$, and  $\min\{r,s\}=\max\{t,w\}=1$, otherwise either $z=x$ or $z=y$. If $r=s=t=w=1$, then $z=0$ and the result is immediate. Consequently $\max\{r,s\}=\max\{t,w\}=1$, $\min\{r,s\}<1$ and
 $\min\{t,w\}<1$. Dropping the symmetric  terms, it follows that
 $$z=rx\boxplus sx\boxplus ty\boxplus wy=\min\{r,s\}x\boxplus \min\{t,w\}y.$$
Since $\min\{r,s\}<1$ and $\min\{t,w\}<1$, this implies that $d_{\boxplus}(x,y)= d_{\boxplus}(x,z)=d_{\boxplus}(z,y).$$\Box$\\

\section{ Some Trigonometric Properties in Limit}\label{trigo}

In this section, we show that certain properties related to the notion of orthogonality can be derived within the framework of the idempotent and symmetric algebraic structure we consider. 

\subsection{Orthogonality and Pythagorean Properties in Limit}

In the following, for all $p\in \mathbb N$, we say that $x,y\in \Real^n$ are {\bf $\varphi_p$-orthogonal} if:
\begin{equation}
\langle x,y\rangle_p=0.
\end{equation}
In addition, we say that $x,y\in \Real^n$ are {\bf $\digamma$-orthogonal} if:
\begin{equation}
\langle x,y\rangle_\infty=0.
\end{equation}

In the sequel, the earlier results are extended to the general case of a triple $[x,y,z]$ of $\Real^n$.
Notice that we have for all natural numbers $p$
\begin{equation}
  \big \langle x\stackrel{p}{-}z, y\stackrel{p}{-}z\big \rangle_p=\stackrel{\varphi_p}{\sum_{i\in [n]}}x_iy_i
  \stackrel{p}{+}\stackrel{\varphi_p}{\sum_{i\in [n]}}(-x_iz_i)
  \stackrel{p}{+}\stackrel{\varphi_p}{\sum_{i\in [n]}}(-y_iz_i)\stackrel{p}{+}\stackrel{\varphi_p}{\sum_{i\in [n]}}{z_i}^2.
  \end{equation}
Therefore:
{  \begin{align}\lim&_{p\longrightarrow \infty}\big \langle x\stackrel{p}{-}z, y\stackrel{p}{-}z\big \rangle_p=\digamma_{[4n]}\big(x_1y_1,...,x_ny_n,-x_1z_1,...,-x_nz_n, -y_1z_1,...,-y_nz_n, {z_1}^2,...,{z_n}^2\big) . \end{align}}
In the following, we introduce the operation $\langle\! \langle \cdot,\cdot,\cdot\rangle\! \rangle_\boxplus: \Real^n\times \Real^n\times \Real^n
\longrightarrow \Real$ defined as:

\begin{equation}
 \langle\! \langle x,y,z\rangle \!\rangle_\infty= \lim_{p\longrightarrow \infty}\big \langle x\stackrel{p}{-}z, y\stackrel{p}{-}z\big \rangle_p .
\end{equation}

{In the following}, for all $x,y,z\in \Real^n$, {we say that} $[x,y,z]$ is $\digamma$-{\bf right-angled in $z$} if \begin{equation}
   \langle\! \langle x,y,z\rangle \!\rangle_\infty=0.
\end{equation}

In general, a right-angled triangle is not a $\digamma$-right-angled triple.    This is shown in the next example.

\begin{expl}Let $x=(1,-2,-1)$,$y=(2,3,-2)$ and $z=(3, 2,-3)$. We have
$x-z=(-2,-4,2)$ and $y-z=(-1,1,1)$. We have $\langle x-z,y-z\rangle=0$. Therefore $[x,y,z]$ is   a right-angled triangle.However
$\digamma_{{[12]}} (2,-6,2,3, 2,-3,-6,-6,-6,9,4,9)=9\not=0$. Hence $[x,y,z]$ is not a $\digamma$-right-angled triple.
\end{expl}
The next result was established in \cite{b19}. It plays an important role in the remainder. It means that if the limit sum is zero then the $\phi_p$-generalized sum is zero for any $p$.

\begin{lem} \label{induc}Suppose that there is some $x=(x_1,...,x_n)\in \Real^n$ such that $\bigboxplus_{i\in [n]}x_i=0$. Then for all $p\in \mathbb N$, $\stackrel{\varphi_p}{\sum_{i\in [n]}}x_i=0$.
\end{lem}

\begin{prop} Let $ x,y,z\in \Real^n$, and let $[x,y,z]$ be a triple.
  The triple $[x,y,z]$ is $\digamma$-right-angled in $z$ all if and only if for all $p\in \mathbb N$, $[x ,y,z]$ is $\varphi_p$-right-angled in $z$, i.e. $\big \langle x\stackrel{p}{-}z, y\stackrel{p}{-}z\big \rangle_p=0$.

 \end{prop}
{\bf Proof:}  If $[x,y,z]$ is a $\digamma$-right-angled triple in $z$ then:

$$ \langle\! \langle x,y,z\rangle \!\rangle_\infty=0.
$$It follows that
$$  \big \langle x\stackrel{p}{-}z, y\stackrel{p}{-}z\big \rangle_p=0
$$for all $p\in \mathbb N$ and the first implication follows.

Conversely if for all $p$,  $[x,y,z]$ is a  $\varphi_p$-right-angled triple in $z$, then for all $p$
$$
  \big \langle x\stackrel{p}{-}z, y\stackrel{p}{-}z\big \rangle_p=0.
$$
Hence
$$
  \lim_{p\longrightarrow \infty}\big \langle x\stackrel{p}{-}z, y\stackrel{p}{-}z\big \rangle_p=\langle\! \langle x,y,z\rangle \!\rangle_\infty=0,
$$
which proves the reciprocal.  $\Box$\\

\begin{center}
{\scriptsize

\unitlength 0.37mm % = 1.053pt
\linethickness{0.4pt}
\ifx\plotpoint\undefined\newsavebox{\plotpoint}\fi % GNUPLOT compatibility
\begin{picture}(344.5,122.5)(0,0)
\put(64.25,8){\vector(0,1){108.5}}
\put(0,64.25){\vector(1,0){134.75}}
%\emline(26,110.5)(52.25,85.25)
\multiput(26,110.5)(.094765343,-.0911552347){277}{\line(1,0){.094765343}}
%\end
%\emline(52.25,85.25)(14,43)
\multiput(52.25,85.25)(-.0910714286,-.1005952381){420}{\line(0,-1){.1005952381}}
%\end
\put(14,43){\line(0,1){0}}
%\emline(14,43)(25.75,111.25)
\multiput(14,43)(.091085271,.529069767){129}{\line(0,1){.529069767}}
%\end
%\emline(42.5,47)(64,64)
\multiput(42.5,47)(.114973262,.090909091){187}{\line(1,0){.114973262}}
%\end
%\emline(64,64)(96,25)
\multiput(64,64)(.0911680912,-.1111111111){351}{\line(0,-1){.1111111111}}
%\end
%\emline(43,46.5)(95.75,25.5)
\multiput(43,46.5)(.228354978,-.090909091){231}{\line(1,0){.228354978}}
%\end
\put(93,99){\line(0,-1){47.5}}
\put(93,51.5){\line(1,0){34}}
%\emline(126.75,51.5)(92.75,99)
\multiput(126.75,51.5)(-.091152815,.1273458445){373}{\line(0,1){.1273458445}}
%\end
\put(141.75,64.5){\makebox(0,0)[cc]{$x_1$}}
\put(62.75,122){\makebox(0,0)[cc]{$x_2$}}
\put(68,67.5){\makebox(0,0)[cc]{$0$}}
\put(56.25,88.5){\makebox(0,0)[cc]{$a$}}
\put(24.5,114.75){\makebox(0,0)[cc]{$b$}}
\put(10,40.25){\makebox(0,0)[cc]{$c$}}
\put(92,102.5){\makebox(0,0)[cc]{$e$}}
\put(131.25,50){\makebox(0,0)[cc]{$d$}}
\put(92.25,45.75){\makebox(0,0)[cc]{$f$}}
\put(99,22.5){\makebox(0,0)[cc]{$h$}}
\put(38,45.5){\makebox(0,0)[cc]{$g$}}
\put(64.75,0){\makebox(0,0)[cc]
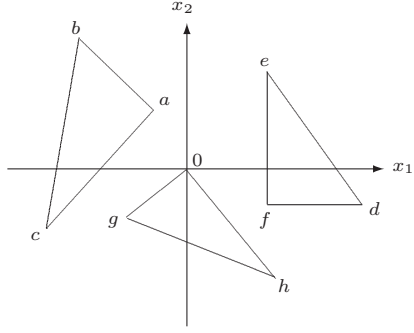
{{\bf Fig. 4:}   $\digamma$-right-angled triples with convex lines. }}
\put(267,8.5){\vector(0,1){108.5}}
\put(202.75,64.75){\vector(1,0){134.75}}
%\emline(245.25,47.5)(266.75,64.5)
\multiput(245.25,47.5)(.114973262,.090909091){187}{\line(1,0){.114973262}}
%\end
%\emline(266.75,64.5)(298.75,25.5)
\multiput(266.75,64.5)(.0911680912,-.1111111111){351}{\line(0,-1){.1111111111}}
%\end
\put(295.75,99.5){\line(0,-1){47.5}}
\put(329.5,52){\line(0,1){47.5}}
\put(295.75,52){\line(1,0){34}}
\put(329.5,99.5){\line(-1,0){34}}
\put(344.5,65){\makebox(0,0)[cc]{$x_1$}}
\put(265.5,122.5){\makebox(0,0)[cc]{$x_2$}}
\put(270.75,68){\makebox(0,0)[cc]{$0$}}
\put(259,89){\makebox(0,0)[cc]{$a$}}
\put(227.25,115.25){\makebox(0,0)[cc]{$b$}}
\put(212.75,40.75){\makebox(0,0)[cc]{$c$}}
\put(294.75,103){\makebox(0,0)[cc]{$e$}}
\put(334,50.5){\makebox(0,0)[cc]{$d$}}
\put(295,46.25){\makebox(0,0)[cc]{$f$}}
\put(301.75,23){\makebox(0,0)[cc]{$h$}}
\put(240.75,46){\makebox(0,0)[cc]{$g$}}
\put(267.5,.5){\makebox(0,0)[cc]
{{\bf Fig. 5:} $\digamma$-right-angled triple{s} }}
\put(267.5,-7){\makebox(0,0)[cc]
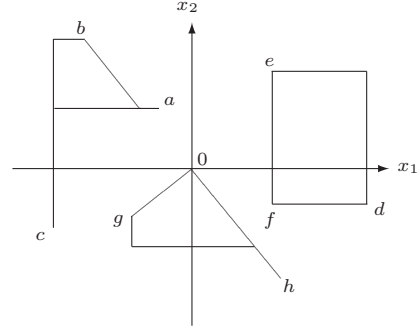
{  with {idempotent symmetric  convex} lines.}}
\put(217.5,43.5){\line(0,1){67.5}}
\put(217.5,111){\line(1,0){11.25}}
%\emline(228.75,110.75)(248.25,86.25)
\multiput(228.75,110.75)(.091121495,-.114485981){214}{\line(0,-1){.114485981}}
%\end
\put(248.25,86.25){\line(1,0){7}}
\put(217.75,86.25){\line(1,0){31}}
\put(289.25,36.75){\line(-1,0){43.75}}
\put(245.5,36.75){\line(0,1){11}}
\end{picture}

}

\end{center}
\medskip

\medskip
\begin{rem}A $\digamma$-right-angled triple is a right-angled triangle. However, the converse is not true. In Figure 4, $[a,b,c]$ is right-angled but not $\digamma$-right-angled. The triples $[d,f,e]$ and $[g,0,h]$ are $\digamma$-right-angled.  

In Figure 5, the same triples are represented using {idempotent symmetric convex} lines. Clearly, the right angle is preserved only for the triples $[d,f,e]$ and $[g,0,h]$.

\end{rem}

In the following, for all subsets $I$ of $[n]$ and all $x\in \Real^n$, let us denote $x_{[I]}=\sum_{i\in I}x_ie_i$, where $\{e_1,...,e_n\}$ is the canonic basis of $\Real^n$.

\begin{prop} Let $ x,y,z\in \Real^n$, and let $[x,y,0]$ is a $\digamma$-right-angled triple in $0$. Then there exists a partition  $\{I_1,...,I_{m}, J\}$ of   $[n]$ with $ m= \lfloor\frac{n}{2}\rfloor $, $\Card\,( J)\in \{0,1\}$  and such that for all $k$:
$$\Big[x_{[I_k]}, y_{[I_k]}, 0_{[I_k]}\Big]$$ is a   right-angled triangle in $0_{[I_k]}$ of $\Real^{I_k}$ {with} $\Card \,(I_k)=2$. Moreover, if $J=\{j\}$ is nonempty then $ x_j   y_j  =0$.
\end{prop}
{\bf Proof:} If $[x,y,0]$ is a $\digamma$-right-angled triple in $0$ then:

$$\bigboxplus_{i\in [n]} x_iy_i=0.$$
 $\bigboxplus_{i\in [n]}x_iy_i=0 $ implies that there exists a partition  $\{I_1,...,I_{m}, J\}$ of   $[n]$ with $ m= \lfloor \frac{n}{2}\rfloor$, $\Card \,(I_k)=2$ for all $k$ and $\Card\,( J)\in \{0,1\}$ such that:
$$\sum_{i\in [I_k]}x_iy_i=\bigboxplus_{i\in [I_k]}x_iy_i=0,$$
and $x_jy_j=0$ if $J=\{j\}$ is a nonempty set.
Consequently $$\sum_{i\in [I_k]} x_i   y_i  =0.  $$
Suppose that $I_k=\{i',i''\}$. Then
$$ x_{i'}  y_{i' }\boxminus   x_{i''} y_{i''}=0\iff  x_{i'}  y_{i' }-  x_{i''} y_{i''}=0.$$
It follows that $\langle x_{[I_k]},y_{[I_k]}\rangle_p=0$ for all $p.$ $\Box$

\begin{cor}Suppose that $n=2$. If $[x,y,0]$ is a  right-angled triangle in $z$ of $\Real^2$, then $[x,y,0]$ is $\digamma$-right-angled in $0$ and it is $\varphi_p$-right-angled in $0$ for all natural numbers $p$.

\end{cor}

The following result is an idempotent and non-associative analogue of the Pythagorean theorem.

\begin{prop}Let   $ x,y,z\in \Real^n$, and let $[x,y,z]$ is a $\digamma$-right-angled triple in $z$. Then for all $p\in \mathbb N$:
$$ d^2_{{\varphi_p}}(x,y) = d^2_{{\varphi_p}}(x,z) \stackrel{p}{+}d^2_{{\varphi_p}}(y,z),$$
and
$$d_\boxplus(x,y)=\max\{d_\boxplus(x, z),d_\boxplus(y,z)\}.$$

\end{prop}
{\bf Proof:} If $[x,y,z]$ is a $\digamma$-right-angled triple in $z$ then $$\langle\! \langle x,y,z\rangle \!\rangle_\infty=0.$$ Thus, for all $p$, $\langle x\stackrel{p}{-}z,y\stackrel{p}{-}z\rangle_p=0$. It follows that
$$\langle \phi_p(x)-\phi_p(z),\phi_p(y)-\phi_p(z)\rangle=0.$$
Hence $$\|\phi_p(x)-\phi_p(y)\|^2=\|\phi_p(x)-\phi_p(z)\|^2+\|\phi_p(y)-\phi_p(z)\|^2. $$
We deduce that
{\small \begin{align*}&\Big\|\phi_p\Big(\phi_p^{-1}\big(\phi_p(x)-\phi_p(y)\big)\Big)\Big\|^2=\varphi_p \Big( \varphi_p^{-1}(\|\phi_p\Big(\phi_p^{-1}\big(\phi_p(x)-\phi_p(z)\big)\Big)\|^2)\Big)+\varphi_p \Big( \varphi_p^{-1}( \|\phi_p\Big(\phi_p^{-1}\big(\phi_p(y)-\phi_p(z)\big)\Big)\|^2 \big)\Big). \end{align*}}
Thus:
$$d^2_{{\varphi_p}}(x,y)=\|x\stackrel{p}{-}y\|_{{\varphi_p}}^2=\|x\stackrel{p}{-}z\|_{{\varphi_p}} \stackrel{p}{+}\|y\stackrel{p}{-}z\|_{{\varphi_p}}=d^2_{{\varphi_p}}(x,z) \stackrel{p}{+}d^2_{{\varphi_p}}(y,z). $$
Taking the limit, yields:
$$d^2_\boxplus(x,y)=\lim_{p\longrightarrow \infty}d^2_{{\varphi_p}}(x,y)=\lim_{p\longrightarrow \infty}d^2_{{\varphi_p}}(x,z) \stackrel{p}{+}d^2_{{\varphi_p}}(y,z)=\max\{d^2_\boxplus(x, z),d^2_\boxplus(y,z)\}. $$
The final statement immediately follows. $\Box$
\subsection{Cosine, Sine in Limit and some Parametrization of the Unit Square}

The purpose of this section is to demonstrate that we can use the limiting properties of standard trigonometric concepts to derive analogues of well-known geometric properties within the framework of the non-associative algebra we are studying. 
Paralleling the usual definition of a cosine, let us define the $\varphi_p$-pseudo-cosine of the angle between two vectors as:
\begin{equation}\cos_p(x,y)=\frac{\langle x,y\rangle_p}{\|x\|_{{\varphi_p}}\|y\|_{{\varphi_p}}}.\end{equation}
Taking the limit, we define the pseudo $\digamma$-cosine as:
\begin{equation}\cos_\infty(x,y)=\frac{\langle x,y\rangle_\infty}{\|x\|_{\infty}\|y\|_{\infty}}.\end{equation}

It is now natural to define a  {sine function in limit}. In \cite{b20} a $\varphi_p$ exterior product is defined for all $v_1,v_2, ...,v_n\in \Real^n$ 
as:
\begin{equation}
\big(v_1\stackrel{p}{\wedge} v_2\stackrel{p}{\wedge}\cdots \stackrel{p}{\wedge}v_n\big )=|v_1,v_2,...,v_n|_p\big(e_1{\wedge} e_2 {\wedge}\cdots  {\wedge}e_n\big ){,}
\end{equation}
where:
\begin{align}
\big(v_1\stackrel{p}{\wedge} v_2\stackrel{p}{\wedge}\cdots \stackrel{p}{\wedge}v_n\big )&=\varphi_p^{-1}\big(
 \phi_p(v_1) {\wedge}\phi_p( v_2) {\wedge}\cdots  {\wedge}\phi_p(v_n)\ \big)\\
 &=\varphi_p^{-1}\Big(|\phi_p(v_1),\phi_p(v_2),...,\phi_p(v_n)|\Big) \big( e_1  {\wedge}  e_2  {\wedge}\cdots  {\wedge} e_n \big ).
\end{align}
Taking the limit yields:
\begin{equation}
\big(v_1\stackrel{\infty}{\wedge} v_2\stackrel{\infty}{\wedge}\cdots \stackrel{\infty}{\wedge}v_n\big )=|v_1,v_2,...,v_n|_\infty\big(e_1\stackrel{\infty}{\wedge} e_2\stackrel{\infty}{\wedge}\cdots \stackrel{\infty}{\wedge}e_n\big ).
\end{equation}
In the case of two vectors the exterior product can be identified to the exterior product and the oriented area of the cross product can be interpreted as the positive area of the parallelogram having $x$ and $y$ as sides:
\begin{equation}
{ x\wedge y=\left |x,y \right |=\left\|x \right\|\left\|y \right\| \sin (x,y) .}
\end{equation}
The $\varphi_p$-sine function for two vectors $x,y\in \Real^2$ is defined as
\begin{equation}
\sin_p(x,y)=\frac{|x,y|_p}{{\|x\|_{{\varphi_p}}\|y\|_{{\varphi_p}}}}=\frac{\varphi_p^{-1}\big(\varphi_p(x_1)\varphi_p(y_2)-\varphi_p(x_2)
\varphi_p(y_1)\big)}{{\|x\|_{{\varphi_p}}\|y\|_{{\varphi_p}}}}.
\end{equation}
Along this line, let us define
\begin{equation}
\sin_\infty(x,y)=\frac{|x,y|_\infty}{ \|x\|_{\infty}\|y\|_{\infty}}=\lim_{p\longrightarrow \infty}\sin_p(x,y)=\frac{x_1y_2\boxminus x_2y_1}{\max\{|x_1,x_2|\}\max\{|y_1,y_2|\}}.
\end{equation}

\begin{lem}\label{trigo}For all $x,y\in \Real^n$, we have the following properties. \\
$(a)$ $ \cos_p(x,y)=\varphi_p^{-1}\big ( \cos \big (\phi_p(x),\phi_p(y)) \big)$ and $ \sin_p(x,y)=\varphi_p^{-1} \big(\sin (\phi_p(x),\phi_p(y) )\big).$\\
$(b)$ ${\cos_p(x,y)}^2\stackrel{p}{+}{\sin_p(x,y)}^2=1$.\\
$(c)$ $\max\{|\cos_\infty(x,y)|, |\sin_\infty(x,y)|\}=1$.\\
\end{lem}
{\bf Proof:} $(a)$ The first statement is immediate. To prove the second part note that:
\begin{equation*}
\sin_p(x,y)=\frac{|x,y|_p}{{\|x\|_{{\varphi_p}}\|y\|_{{\varphi_p}}}}=\frac{\varphi_p^{-1}\big(|\phi_p(x), \phi_p(y)|)}{{\varphi_p^{-1}(\|\phi_p(x)\|)\varphi_p^{-1}(\|\phi_p(y)\|)}}.
\end{equation*}
Therefore
\begin{equation*}
\sin_p(x,y)= \varphi_p^{-1}\Big(\frac{|\phi_p(x), \phi_p(y)|}{{ \|\phi_p(x)\| \|\phi_p(y)\| }}\Big)=\varphi_p^{-1} \big(\sin (\phi_p(x),\phi_p(y)) \big).
\end{equation*}
$(b)$ We have 
\begin{align*}  {\cos_p (x,y)}^2  \stackrel{p}{+}{\sin_p(x,y)}^2 &=\varphi_p \Big(\big[\varphi_p^{-1}\big ( \cos \big (\phi_p(x),\phi_p(y)) \big)\big]^2\Big)+
\varphi_p \Big(\big[\varphi_p^{-1}\big ( \sin \big (\phi_p(x),\phi_p(y)) \big) 
\big]^2\Big)\\
&= \varphi_p \Big(\varphi_p^{-1}\big [\big ( \cos \big (\phi_p(x),\phi_p(y))\big]^2 \big)\Big)+
\varphi_p \Big(\varphi_p^{-1}\big ( \big[\sin \big (\phi_p(x),\phi_p(y)) \big]^2\big) 
\Big)\\&= \big ( \cos \big (\phi_p(x),\phi_p(y)) )\big)^2 +
 \big (  \sin \big (\phi_p(x),\phi_p(y)) \big)^2=1.  \end{align*}
$(c)$ From $(a)$ and $(b)$, we have
$$\max\{|\cos_\infty(x,y)|^2, |\sin_\infty(x,y)|^2\}=\lim_{p\longrightarrow \infty} {\cos_p (x,y)}^2  \stackrel{p}{+}{\sin_p(x,y)}^2=1.  $$
Therefore, 
$$\max\{|\cos_\infty(x,y)| , |\sin_\infty(x,y)| \}= 1. \quad \Box  $$ \\

Let   $C_\infty(0,1]=\big\{(x_1,x_2)\in \Real^2: \max\{|x_1|,|x_2|\}=1\big\}$ be the circle defined with respect to the $\|\,\|_\infty$
metric in $\Real^2$. Let us consider the points $
a =(1,0)$, $b= (1,1)$, $c=(0,1)$, $d=(-1,1)$, $e=(-1,0)$, $f=(-1,-1)$, $g=(0,-1)$, $h=(1,-1)$.  By construction:

\begin{equation}
C_\infty(0,1]=C_\boxplus(0,1] =[a,b]\cup [ b,d]\cup [  d,f]\cup [   f,h]\cup [h,a].
\end{equation}

% This is a LaTeX picture output by TeXCAD.
% File name: [Clipboard].
% Version of TeXCAD: 4.51
% Reference / build: 27-Nov-2018 (rev. a75)
% For new versions, check: http://texcad.sf.net/
% Options on the following lines.
%\grade{\on}
%\emlines{\off}
%\epic{\off}
%\beziermacro{\on}
%\reduce{\on}
%\snapping{\off}
%\pvinsert{% Your \input, \def, etc. here}
%\quality{8.000}
%\graddiff{0.005}
%\snapasp{1}
%\zoom{4.0000}
\begin{center}{\scriptsize 
% This is a LaTeX picture output by TeXCAD.
% File name: [Clipboard].
% Version of TeXCAD: 4.51
% Reference / build: 27-Nov-2018 (rev. a75)
% For new versions, check: http://texcad.sf.net/
% Options on the following lines.
%\grade{\on}
%\emlines{\off}
%\epic{\off}
%\beziermacro{\on}
%\reduce{\on}
%\snapping{\off}
%\pvinsert{% Your \input, \def, etc. here}
%\quality{8.000}
%\graddiff{0.005}
%\snapasp{1}
%\zoom{4.0000}
\unitlength 0.37mm % = 1.053pt
\linethickness{0.4pt}
\ifx\plotpoint\undefined\newsavebox{\plotpoint}\fi % GNUPLOT compatibility
\begin{picture}(137.5,150.75)(0,0)
\put(63,10){\vector(0,1){134}}
\put(0,72){\vector(1,0){130}}
\put(58.5,76.25){\makebox(0,0)[cc]{$0$}}
\put(137.5,71.5){\makebox(0,0)[cc]{$x_1$}}
\put(62.75,150.75){\makebox(0,0)[cc]{$x_2$}}
\put(63.25,71.75){\circle*{1.581}}
\put(63.25,109.75){\circle*{1.581}}
\put(63,37){\circle*{1.581}}
\put(101,109.5){\circle*{1.581}}
\put(100.75,36.75){\circle*{1.581}}
\put(100.75,71.75){\circle*{1.581}}
\put(28,109.5){\circle*{1.581}}
\put(27.75,36.75){\circle*{1.581}}
\put(27.75,71.75){\circle*{1.581}}
%\emline(63.75,72.5)(90.25,110.25)
\multiput(63.75,72.5)(.0910652921,.1297250859){291}{\line(0,1){.1297250859}}
%\end
\thicklines
%\dashline{1}(90.5,109.5)(90.5,72.25)
\put(90.31,109.31){\line(0,-1){.9803}}
\put(90.31,107.349){\line(0,-1){.9803}}
\put(90.31,105.389){\line(0,-1){.9803}}
\put(90.31,103.428){\line(0,-1){.9803}}
\put(90.31,101.468){\line(0,-1){.9803}}
\put(90.31,99.507){\line(0,-1){.9803}}
\put(90.31,97.547){\line(0,-1){.9803}}
\put(90.31,95.586){\line(0,-1){.9803}}
\put(90.31,93.626){\line(0,-1){.9803}}
\put(90.31,91.665){\line(0,-1){.9803}}
\put(90.31,89.705){\line(0,-1){.9803}}
\put(90.31,87.744){\line(0,-1){.9803}}
\put(90.31,85.784){\line(0,-1){.9803}}
\put(90.31,83.823){\line(0,-1){.9803}}
\put(90.31,81.863){\line(0,-1){.9803}}
\put(90.31,79.902){\line(0,-1){.9803}}
\put(90.31,77.942){\line(0,-1){.9803}}
\put(90.31,75.981){\line(0,-1){.9803}}
\put(90.31,74.021){\line(0,-1){.9803}}
%\end
\put(27.75,37.25){\framebox(72.75,72.5)[cc]{}}
\put(84,66.75){\makebox(0,0)[cc]{$\mathrm{pcos} (\theta)$}}
\put(60,0){\makebox(0,0)[cc]{{\bf Fig. 6:} Pseudo Cosine.}}
\put(105.25,77.75){\makebox(0,0)[cc]{$a$}}
\put(105.25,112.25){\makebox(0,0)[cc]{$b$}}
\put(65.75,113.25){\makebox(0,0)[cc]{$c$}}
\put(24,112.25){\makebox(0,0)[cc]{$d$}}
\put(24,65.25){\makebox(0,0)[cc]{$e$}}
\put(24,32.5){\makebox(0,0)[cc]{$f$}}
\put(60.25,32.5){\makebox(0,0)[cc]{$g$}}
\put(104.25,32.5){\makebox(0,0)[cc]{$h$}}
\put(83.75,22.75){\makebox(0,0)[cc]{$C_\infty(0,1)$}}
\end{picture}

}
\end{center}
\medskip\medskip
\medskip\medskip

 In the following we will identify   the angles of the square with the length of the segments of its circumference. Let us consider the  map  $\alpha: C_\infty(0,1]\longrightarrow [0,8]$ defined by:

\begin{equation}
\alpha(x_1,x_2)=\left\{ \begin{matrix}x_2&\text{ if }& x\in [a,b]\\2-x_1&\text{ if }& x\in [b,d]\\4-x_2&\text{ if }& x\in [d,f]\\
6+x_1&\text{ if }& x\in [f,h]\\8+x_2&\text{ if }& x\in [h,a]\end{matrix}\right..
\end{equation}
Clearly, we have $\alpha(a)=0$, $\alpha(b)=1$,  $\alpha(c)=2$,  $\alpha(d)=3$,  $\alpha(e)=4$,  $\alpha(f)=5$,  $\alpha(g)=6$, $\alpha(h)=7.$

 It is easy to prove that this map is a bijection from $  C_\infty(0,1]$ to $ [0,8]$. Notice that $8$ is the perimeter of $C_\infty(0,1]$. For all $\theta\in [0,8]$ we define the pseudo-cosine $\mathrm{pcos} : \Real\longrightarrow [-1,1]$ as a periodic function with period 8:
 \begin{equation}
 \mathrm{pcos} (\theta)=\cos_\infty\big (\alpha^{-1}(\theta), e_1\big )\quad \mod 8.
 \end{equation}
 It follows that for any integer $k\in\mathbb Z$
 \begin{equation}
\mathrm{pcos} (\theta+ 8k)=\left\{ \begin{matrix}1&\text{ if }& \theta\in [0,1]\\2-\theta&\text{ if }&  \theta\in [1,3]\\-1&\text{ if }& \theta\in [3,5]\\
-6+\theta&\text{ if }& \theta\in [5,7]\\1&\text{ if }& \theta \in [7,8]\end{matrix}\right..
\end{equation}

\begin{center}{\scriptsize
\unitlength 0.4mm % = 2.845pt
\linethickness{0.4pt}
\ifx\plotpoint\undefined\newsavebox{\plotpoint}\fi % GNUPLOT compatibility
\begin{picture}(246.25,140.5)(0,0)
\put(22,72.75){\vector(1,0){213}}
\put(121.25,9.75){\vector(0,1){130.75}}
\put(206.75,73){\circle*{1.118}}
\put(34.25,73){\circle*{1.803}}
\put(121.25,132.25){\circle*{1.118}}
\put(121.25,22.5){\circle*{1.118}}
\put(164.25,72.5){\circle*{1.118}}
\put(142.5,72.75){\circle*{1.5}}
\put(187.25,72.5){\circle*{1.414}}
\put(76.25,73){\circle*{1.118}}
\put(99.75,72.5){\circle*{1.118}}
\put(54,72.75){\circle*{1}}
\put(102,132.25){\line(1,0){36.75}}
%\emline(138.75,132.25)(187.25,22.25)
\multiput(138.75,132.25)(.03372739917,-.07649513213){1438}{\line(0,-1){.07649513213}}
%\end
%\emline(101.75,132.5)(53.75,22.75)
\multiput(101.75,132.5)(-.03373155306,-.07712579058){1423}{\line(0,-1){.07712579058}}
%\end
%\emline(53.75,22.75)(15,22.5)
\multiput(53.75,22.75)(-4.84375,-.03125){8}{\line(-1,0){4.84375}}
%\end
%\emline(15,22.5)(0,52)
\multiput(15,22.5)(-.0337078652,.0662921348){445}{\line(0,1){.0662921348}}
%\end
\put(186.75,22.75){\line(1,0){39.5}}
%\emline(226.25,22.75)(246.25,62.5)
\multiput(226.25,22.75)(.0337268128,.0670320405){593}{\line(0,1){.0670320405}}
%\end
\put(109.5,22.5){\makebox(0,0)[cc]{$-1$}}
\put(113.25,138){\makebox(0,0)[cc]{$+1$}}
\put(113.25,76){\makebox(0,0)[cc]{$0$}}
\put(142.25,66.25){\makebox(0,0)[cc]{$ 1$}}
\put(164.25,67.5){\makebox(0,0)[cc]{$2$}}
\put(187.25,68.25){\makebox(0,0)[cc]{$3$}}
\put(207,69.25){\makebox(0,0)[cc]{$4$}}
\put(99.5,67.5){\makebox(0,0)[cc]{$-1$}}
\put(76,68.5){\makebox(0,0)[cc]{$-2$}}
\put(53.75,68.75){\makebox(0,0)[cc]{$-3$}}
\put(33,67.25){\makebox(0,0)[cc]{$-4$}}
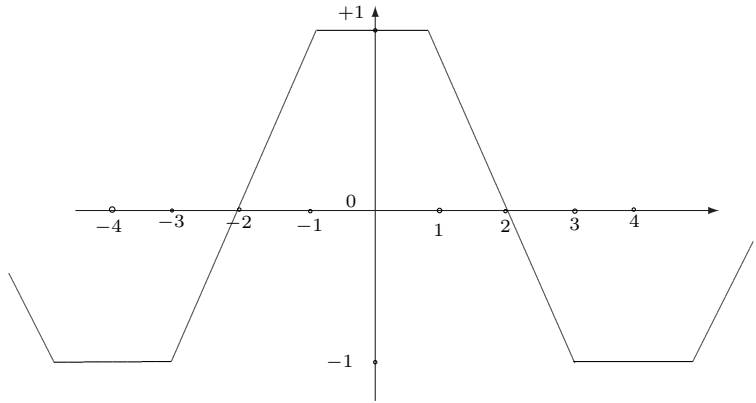
\put(121.25,0){\makebox(0,0)[cc]{{\bf Fig. 7:} Pseudo Cosine Function }}
\end{picture}}
\end{center}

Parallel to the previous definition, for all \( \theta \in [0,8] \), we define the pseudo-sine function \( \mathrm{psin}: \mathbb{R} \longrightarrow [-1,1] \) as a periodic function with period 8:

 \begin{equation}
 \mathrm{psin}(\theta)=\sin_\infty\big (\alpha^{-1}(\theta), e_1\big )\quad \mod 8.
 \end{equation}
 It follows that for any integer $z\in\mathbb Z$
 \begin{equation}
\mathrm{psin}(\theta+ 8k)=\left\{ \begin{matrix}\theta&\text{ if }& \theta\in [0,1]\\1&\text{ if }&  \theta\in [1,3]\\4-\theta&\text{ if }& \theta\in [3,5]\\
-1&\text{ if }& \theta\in [5,7]\\8-\theta&\text{ if }& \theta \in [7,8]\end{matrix}\right..
\end{equation}
{Hence, for all} $\theta \in \Real$, we have from Lemma \ref{trigo}
\begin{equation}
\max \{\big|\mathrm{pcos}(\theta)\big|,\big|\mathrm{psin}(\theta)\big|\}=\max \Big\{\big|\cos_\infty \big(\alpha^{-1}(\theta), e_1\big)\big|,\big|\sin_\infty \big(\alpha^{-1}(\theta), e_1\big)\big|\Big\}=1.
\end{equation}

\subsection{Some Formalism on the Complex Scalar Field}

We can continue the previous analogies, this time using complex numbers and studying the relative properties of the unit square.
Let us denote $C_{{\varphi_p}}(0,1]=\{x\in \Real^n: \|x\|_{{\varphi_p}}=1\}$. We have the relations $\phi_p\Big(C_{{\varphi_p}}(0,1]\Big)=C(0,1]$. 
Now, let us consider the subsets of $\mathbb C$ defined as $ C^\natural_{{\varphi_p}}(0,1]=\{z=a+\mathrm{i}b\in \mathbb C: \|(a,b)\|_{{\varphi_p}}=1\}$ and 
$ C^\natural_{\infty}(0,1]=\{z=a+\mathrm{i}b\in \mathbb C: \|(a,b)\|_{\infty}=1\}$. Let us consider the bijective  map $\varphi^\natural: \mathbb C\longrightarrow \mathbb C$ defined as:
\begin{equation}
\varphi_p^\natural(a+\mathrm ib)=\varphi_p(a)+\mathrm i\varphi_p(b). 
\end{equation}
Its reciprocal is the map ${\varphi^\natural}^{-1}: \mathbb C\longrightarrow \mathbb C$ defined as 
${\varphi_p^\natural}^{-1}(a+\mathrm ib)=\varphi_p^{-1}(a)+\mathrm i\varphi_p^{-1}(b). $ The map $\varphi_p^\natural$ is {homogeneous} {with respect to} the real scalar field, i.e. for all $\rho\in \Real $, $ 
\varphi_p^\natural\big(\rho(a+\mathrm ib)\big )=\rho\big(\varphi_p(a)+\mathrm i\varphi_p(b)\big)=\rho\varphi_p^\sharp \big(a+\mathrm i b\big). $ Therefore, for all $z=\rho e^{\mathrm i\theta}\in \mathbb C$, with $\rho=|z|$, we have 
\begin{equation} {\varphi_p^\natural}^{-1}(z)=\rho\big (\varphi_p^{-1}(\cos \theta)+\mathrm i\varphi_p^{-1}(\sin \theta)\big). \end{equation}
For all $z,w\in  \mathbb C^2 $, let us define:
\begin{equation}
z\stackrel{p}{\cdot}w={\varphi_p^\natural}^{-1}\Big( \varphi_p^\natural(z) \varphi_p^\natural(w)\Big)
\end{equation}
It follows that if $z=a+\mathrm i\, b$ and $w=c+\mathrm i d$, then:
\begin{equation}z\stackrel{p}{\cdot}w= \varphi_p^{-1}\Big( \varphi_p (a) \varphi_p (b)-  \varphi_p (b) \varphi_p (d)\Big)+\mathrm i  \,\varphi_p ^{-1} \Big(\varphi_p (a) \varphi_p (c)+  \varphi_p (b) \varphi_p (c)\Big).\end{equation}
Equivalently
\begin{equation}z\stackrel{p}{\cdot}w=  (ab\stackrel{p}{-}bd)+\mathrm i \, ( ad\stackrel{p}{+}bc).\end{equation}
Along this line, taking the limit with respect to the product topology of $  \mathbb R^2 $, we define the product: 
\begin{equation}z\boxtimes w= \lim_{p\longrightarrow \infty} z\stackrel{p}{\cdot}w=
(ab\boxminus bd)+\mathrm i \, ( ad\boxplus bc).\end{equation}

Note that this product is equivalent to the matrix product defined in \cite{b20} for matrices representing complex numbers. We have for all complex numbers of the form $z=a +\mathrm i \, b$ and  $w=c +\mathrm i \, d$:
\begin{equation}
\begin{pmatrix} a&-b\\b&a \end{pmatrix}\boxtimes \begin{pmatrix} c&-d\\d&c \end{pmatrix}=\begin{pmatrix} ac\boxminus bd&-(ad\boxplus bc)\\  ad\boxplus bc & ac\boxminus bd \end{pmatrix}.
\end{equation}
Moreover, for all $p\in \mathbb N$, and all $z=a+\mathrm i b\in \mathbb C$ we define the $\varphi_p$-module of $z$ as $|z|_p=\|(a,b)\|_{{\varphi_p}}$. Similarly, 
 $|z|_\infty=\|(a,b)\|_{\infty}=\max\{|a|,|b|\}$. In the following, we also extend the operation $\boxplus$ to the complex numbers. For all $z=a+\mathrm i b$ and $w=c+\mathrm i d$, we have :\begin{equation}z\boxplus w=(a\boxplus c)+\mathrm i (b\boxplus d). \end{equation}
We can then establish the following relations.

\begin{lem} We have the following properties:\\
$(a)$ For all $z\in {\mathbb C}$, there exists some $\theta \in [0,8]$ such that
$$z=\|z\|_\infty (\mathrm{pcos}(\theta)+\mathrm i \,\mathrm{psin}(\theta))$$
$(b)$ If $z=a+\mathrm i\, b$ then $\sqrt{z\boxtimes \bar z}=|z|_\infty=\max\{|a|,|b|\}$.\\
$(c)$ If $z\not=0$, then $z\cdot \frac{1}{|z|_\infty^2} \bar z=1.$\\
$(d)$ Suppose that $z,w\in  C^\natural_{\infty}(0,1] $, then $z\boxtimes w\in C^\natural_{\infty}(0,1]$. \\
$(e)$ Suppose that $z,w\in  \mathbb C $, then $|z\boxtimes w|_\infty =|z|_\infty|w|_\infty$. \\

\end{lem}
{\bf Proof:} $(a)$ Suppose that $z=a+\mathrm i\, b$. If $z=0$, the result is trivial.  If $z\not=0$ then
$$z=\max\{|a|,|b|\}\frac{a+\mathrm i\, b}{ \max\{|a|,|b|\}}=|z|_\infty\Big(\frac{a}{ \max\{|a|,|b|\}}+\mathrm i\frac{ b}{ \max\{|a|,|b|\}}\Big). $$
However  $\Big(\frac{a}{ \max\{|a|,|b|\}}+\mathrm i\frac{ b}{ \max\{|a|,|b|\}}\Big)\in C_\infty^\natural(0,1].$ Hence there is some $\theta\in [0,8]$ such that 
$$ \Big(\frac{a}{ \max\{|a|,|b|\}}+\mathrm i\frac{\, b}{ \max\{|a|,|b|\}}\Big)=\mathrm{pcos}\theta +\mathrm i\,\mathrm{psin}\theta.  $$
Hence 
$$z=\|z\|_\infty \big(\mathrm{pcos}(\theta)+\mathrm i\, \mathrm{psin}(\theta)\big). $$

$(b)$ We have $z\boxtimes \bar z=(a+\mathrm i \, b)\boxtimes (a-\mathrm i \, b)=(a^2\boxplus b^2) + \mathrm i\, (ab\boxminus ab)=a^2\boxplus b^2=\max\{|a|,|b|\}^2$ which yields the result. $(c)$ is immediate from $(b)$. 

$(d)$ Suppose that 
$z=a+\mathrm i\, b$ and $w=c+\mathrm i\, d$. We first assume that $z,w\in  C^\natural_{\infty}(0,1].$ This implies that $\max\{|a|,|b|\}=\max\{|c|,|d|\}=1$. By definition we have:
$$z\boxtimes w=\big(a   c \boxminus b   d\big) +\mathrm i\, \big(a  d \boxplus b  c \big). $$
We need to prove that $\max\big\{|a   c \boxminus b   d|  , |a  d \boxplus b  c|\big\}=1.$ We consider the following cases:\\
\noindent $(i)$ if $|a|<1$ and $|c|<1$ then $|bd|=1$ and consequently  $|a   c \boxminus b d|=1$; $(ii)$ if $|a|<1$ and $|d|<1$ then $|bc|=1$, therefore $|a  d \boxplus b  c |=1$; $(iii)$ if $|c|<1$ and $|b|<1$ then $|ad|=1$ and we deduce that $|a  d \boxplus b  c |=1$; $(iv)$ if $|b|<1$ and $|d|<1$ then $|ac|=1$, hence we have  $|a  c \boxminus b  d |=1$.  

  Suppose now that  $|a|=|b|=|c|=|d|=1$. All we need to show is that we cannot have $a   c \boxminus b   d=a  d \boxplus b  c =0$. In such a case we have 
  $a   c = b   d$ and $a  d =- b  c $. However,  since $|c|=|d|=1 $, this implies that $ { d^2} =-    c^2$ that is a contradiction.  Consequently, we deduce that that if $z,w\in  C^\natural_{\infty}(0,1]$ then $z\boxtimes w\in  C^\natural_{\infty}(0,1]$. 
  
$(e)$ Suppose that $z=|z|_\infty \big(\mathrm{pcos}(\theta)  +\mathrm i\, \mathrm{psin} (\theta)\big) $ and $w=|w|_\infty\big(\mathrm{pcos} (\alpha )+\mathrm i\, \mathrm{psin}( \alpha)\big)$. 
We have: {\small $$z\boxtimes w=|z|_\infty|w|_\infty\big (\mathrm{pcos}(\theta)   \mathrm{pcos}(\alpha) \boxminus \mathrm{psin}(\theta)   \mathrm{psin}(\alpha) +\mathrm i\,  (\mathrm{pcos}(\theta)   \mathrm{psin}(\alpha) \boxplus \mathrm{psin}(\theta)   \mathrm{pcos }(\alpha) )\big). $$}
Hence, since 
{\small $$ \left|\mathrm{pcos}(\theta)   \mathrm{pcos}(\alpha) \boxminus \mathrm{psin}(\theta)   \mathrm{psin}(\alpha) +\mathrm i\,  (\mathrm{pcos}(\theta)   \mathrm{psin}(\alpha) \boxplus \mathrm{psin}(\theta)   \mathrm{pcos }(\alpha) ) \right|_\infty=1 $$}
we deduce the result. $\Box$\\

\section{Algebraic Formulation of the Kuratowski-Peano Limit of a Line}\label{fund}

\subsection{Algebraic Formulation of a Line}

In this section, we demonstrate that the non-associative algebraic structure we consider allows for establishing an algebraic formulation of any limiting line passing through two points. In the standard convex case it is easy to give an algebraic formulation of a line passing trough two different points. In such a case one can define the line passing trough $x$ and $y$ ($x\not=y$) as:
\begin{equation}
\mathcal D(x,y)=\{t x +s y:s+t=1, s,t\in \Real \}.
\end{equation}
Note that there is no restriction imposing the non-negativity of $s$ and $t$ which limits a line to a line-segment. Paralleling this   definition one can define a $\varphi_p$-line as:
\begin{equation}
\mathcal D^{(p)}(x,y)=\big \{t x \stackrel{p}{+}s y:s\stackrel{p}{+}t=1, s,t\in \Real \big\}.
\end{equation}
It was pointed out in \cite{b15} that the subset \( \{ tx \boxplus sy : \max\{s,t\} = 1, s,t \geq 0 \} \) is not path-connected and does not describe, over \( \mathbb{R}^n \), the convex hull in the limit of \( \{x, y\} \).
 To that end it was shown in \cite{b15} that $Co^{\infty}(x,y)=\Big \{t x\boxplus r x\boxplus s y\boxplus w y:
\max\{t,r,s,w\}=1,t,r,s,w\geq 0\Big\}$. In what follows, we take a parallel approach to give the algebraic description of a line. {Let us consider} the following subset:
\begin{equation}\mathcal D^\infty(x,y):=\Big \{t x\boxplus r x\boxplus s y\boxplus w y:
t\boxplus r\boxplus s\boxplus w=1,t,r,s,w\in \Real\Big\}.\end{equation}

\begin{lem} \label{LiInclus} For all $x,y\in \Real^n$ with $x\not=y$,
 $$\mathcal D^\infty(x,y)\subset
\Linf_{p\longrightarrow \infty}\mathcal D^{(p)}(x,y).$$
\end{lem}
{\bf Proof:} Suppose that $z\in \Big \{t x\boxplus r x\boxplus s y\boxplus w y:
 t\boxplus r\boxplus s\boxplus w=1,t,r,s,w\in \Real\Big\}$. By definition,
there {exist} $\alpha_1,\alpha_2,\beta_1,\beta_2\in \Real$ with
$\alpha_1\boxplus \alpha_2 \boxplus \beta_1\boxplus \beta_2=1$ and such that
$$z=\alpha_1 x\boxplus \alpha_2 x \boxplus \beta_1 y\boxplus \beta_2 y.$$
Since $\alpha_1\boxplus \alpha_2 \boxplus \beta_1\boxplus \beta_2=1$, there exists some $p_0\in \mathbb N$ such that, for all $p\geq p_0$, $\alpha_1
\stackrel{p}{+}\alpha_2 \stackrel{p}{+} \beta_1 \stackrel{p}{+}
\beta_2\not=0$. For all $p\geq p_0$, define
$$z^{(p)}=\frac{1}{\alpha_1
\stackrel{p}{+}\alpha_2 \stackrel{p}{+} \beta_1 \stackrel{p}{+}
\beta_2 }\big(\alpha_1 x\stackrel{p}{+}\alpha_2 x\stackrel{p}{+}
\beta_1 y\stackrel{p}{+} \beta_2 y\big).$$
By construction, for all $p\geq p_0$,
$z^{(p)}\in \mathcal D^{(p)}(x,y)$. Taking the limit on both sides yields: \begin{align*}\lim_{p\longrightarrow
\infty}z^{(p)}&= \frac{1}{\alpha_1\boxplus \alpha_2 \boxplus \beta_1\boxplus \beta_2}\big(\alpha_1 x\boxplus\alpha_2 x\boxplus \beta_1 y\boxplus \beta_2
y\big)\\&=\alpha_1 x\boxplus\alpha_2 x\boxplus \beta_1 y\boxplus
\beta_2 y=z.\end{align*} Consequently, $z\in \Linf_{p\longrightarrow
+\infty}\mathcal D^{(p)}(x,y)$.
 $\Box$

\begin{lem}\label{limcond} Let $\big\{(s^{(p)},t^{(p)})\big\}_{p\in \mathbb N}$ be  sequence of $\Real^2$. Suppose there is an increasing subsequence $\{p_q\}_{q\in \mathbb N}$ of natural numbers  such that
$s^{(p_q)}\stackrel{p_q}{+}t^{(p_q)}=1$ and $\lim\limits_{q\longrightarrow \infty}(s^{(p_q)},t^{(p_q)})=(s,t)$. We have the following properties:\\
$(a)$ If $\max\{\,|s|, |t|\, \}>1$ then $s=-t$.\\
$(b)$ If $s=-t$ then $ |s|=|t|\geq 1$.\\
$(c)$ If $\min\{\,|s|, |t|\,\}\leq 1$ then $\max\{ s ,  t \}=1$.\\
$(d)$ If $s\not=-t$ then $\max\{ s ,  t \}=1$.
\end{lem}
{\bf Proof:} $(a)$ In such a case if $|s|>1$ then $t=\lim\limits_{q\longrightarrow \infty}{t_{p_{q }}}=\lim\limits_{q\longrightarrow \infty}1\stackrel{p_{q }}{-}{s_{p_{q }}}=-s$. The proof is symmetrical if
 $|t|>1$. $(b)$ Let us prove that in such a case we have $|s|\geq 1$ and $|t|\geq 1$. If $|s|<1$ then  $t=\lim\limits_{q\longrightarrow \infty}{t_{p_{q }}}=\lim\limits_{q\longrightarrow \infty}1\stackrel{p_{q }}{-}{s_{p_{q }}}=1$, which contradicts $s=-t$. Similarly, if $|t|<1$ we have a contradiction. Therefore $|s|\geq 1$ and $|t|\geq 1$. $(c)$ If $s=1$, then the conclusion is obvious. If $|s|<1$, then  $t=\lim\limits_{q\longrightarrow \infty}{t_{p_{q }}}=\lim\limits_{q\longrightarrow \infty}1\stackrel{p_{q }}{-}{s_{p_{q }}}=1$ and $\max\{s,t\}=1$. To conclude if $s=-1$, then
 $t=\lim\limits_{q\longrightarrow \infty}{t_{p_{q}}}=\lim\limits_{q\longrightarrow \infty}1\stackrel{p_{q }}{+}{1}=1$, which ends the proof of $(c)$. $(d)$ From $(a)$, $s\not=-t$ implies that $\max\{\,|s|, |t|\, \}\leq 1$. Hence  $\min\{\,|s|, |t|\,\}\leq 1$ and the result is then immediate from $(c)$. $\Box$\\
 
 The next lemma is a preliminary result to show the upper limit of the sequence $\{\mathcal D^{{(p)}}\}_{p\in \mathbb N}$ is identical to $\mathcal D^\infty$.
 
 \begin{lem}  Let $x,y\in \Real^n$ with $x\not=y$ and suppose that $z\in \Ls_{p\longrightarrow \infty} \mathcal D^{{(p)}}(x,y)$. Then there exists a sequence  $\big\{(s^{(p)},t^{(p)})\big\}_{p\in \mathbb N}\subset \Real^2$ and an increasing subsequence $\{p_{{q}}\}_{q\in \mathbb N}$ of natural numbers such that
$s^{(p_{{q}})}\stackrel{p_{{q}}}{+}t^{(p_q)}=1$, $\lim\limits_{q\longrightarrow \infty}(s^{(p_q)},t^{(p_q)})=(s,t)$ and 
$$z_i=sx_i\boxplus ty_i,$$
for all $i\in [n]$ such that $sx_i\not=-ty_i$.
 
 \end{lem}
 {\bf Proof:} Suppose that $z\in \Ls_{p\longrightarrow \infty}\mathcal D^{(p)}(x,y)$. In such a case there exists an increasing sequence of natural numbers $\{p_k\}_{k\in \mathbb N}$ and sequence $\{z^{(p_k)}\}_{k\in \mathbb N}$ with $z^{(p_k)}\in \mathcal D^{(p_k)}(x,y)$ for all $k$ and $\lim_{q\longrightarrow \infty}z^{(p_k)}=z$. This implies that there exists a sequence $\{ s^{(p_k)} ,t^{(p_k)}\}_{k\in \mathbb N}\subset\Real\times\Real$ such that
\begin{equation}
z^{(p_k)}=t^{(p_k)} x \stackrel{p_k}{+}s^{(p_k)} y,
\end{equation}
with $s^{(p_k)}\stackrel{p_k}{+}t^{(p_k)}=1$ for all $k$. Let us prove that   $\big\{(s^{(p_k)},t^{(p_k)})\big\}_{q\in \mathbb N}$ is a bounded sequence. Suppose that this is not the case and let us show a contradiction. Suppose for example that $\{s^{(p_k)}\}_{k\in \mathbb N}$ is not bounded. An elementary calculus shows that since $t^{(p_k)}=1\stackrel{p_k}{-}s^{(p_k)}$. Thus:
\begin{equation}
z^{(p_k)}= x \stackrel{p_k}{+}s^{(p_k)} (y \stackrel{p_k}{-}x).
\end{equation}
However   $\lim_{k\longrightarrow\infty }y \stackrel{p_k}{-}x=y\boxminus x$. Since $x\not=y$, $ y\boxminus x\not=0$ and this implies  that $\lim_{k\longrightarrow \infty}\|z^{(p_q)}\|=\infty$ that is a contradiction. Consequently $ \{s^{(p_k)}\}_{k\in \mathbb N}$ is bounded.  Since $s^{(p_k)}=1\stackrel{p_k}{-}t^{(p_k)}$ this also implies that {$ \{t^{(p_k)}\}_{k\in \mathbb N}$} is bounded. Therefore, we deduce that    the sequence $\{(s^{(p_k)},t^{(p_k)})\}_{k\in \mathbb N}$ is   bounded. Consequently one can extract an increasing subsequence $\{p_{k_m}\}_{m\in  {\mathbb N} }$ such that 
$$\lim_{m\longrightarrow \infty}(s^{(p_{k_m})},t^{(p_{k_m})})=(s,p).$$
Let us denote $I=\{i\in [n]: sx_i\boxplus ty_i=0\}$. For all $i\notin I$, since  $sx_i\not=-ty_i$ we have
$$\lim_{m\longrightarrow \infty}s^{(p_{k_m})}x_i \stackrel{p_{k_m}}{+} t^{(p_{k_m})}y_i=sx_i\boxplus ty_i, $$
{which} yields the result. $\Box$\\

\begin{prop} \label{LiInclus}   For all $x,y\in \Real^n$ with $x\not=y$,
 $$\Lim_{p\longrightarrow \infty}\mathcal D^{(p)}(x,y)=\mathcal D^\infty(x,y).$$
\end{prop}
{\bf Proof:} To prove this result, {we need} to establish that $\Ls_{p\longrightarrow \infty}\mathcal D^{(p)}(x,y)\subset\Big \{t x\boxplus r x\boxplus s y\boxplus w y:
 t\boxplus r\boxplus s\boxplus w=1,t,r,s,w\in \Real\Big\}$.
Suppose that $z\in \Ls_{p\longrightarrow \infty}\mathcal D^{(p)}(x,y)$. By hypothesis, there exists an increasing sequence of natural numbers $\{p_q\}_{q\in \mathbb N}$ and sequence $\{z^{(p_q)}\}_{q\in \mathbb N}$ with $z^{(p_q)}=t^{(p_q)} x \stackrel{p_q}{+}s^{(p_q)} y\in \mathcal D^{(p_q)}(x,y)$ for all $q$ {such that} $\lim_{q\longrightarrow \infty}z^{(p_q)}=z$. Moreover, $\lim_{p\longrightarrow \infty}(s^{(p_{q })},t^{(p_{q })})=(s,p)$ and
$$\lim_{p\longrightarrow \infty}s^{(p_{q })}x_i \stackrel{p_q}{+} t^{(p_{q })}y_i=sx_i\boxplus ty_i, \quad (1)$$
{for all $i$ such that $sx_i\not=- ty_i$}. Let  $I=\{i\in [n]: sx_i\boxplus ty_i= 0\}$. 
We consider three cases.

$(a)$ $I=\emptyset$. If $\min \{|s|,|t|\}>1$ then, from Lemma \ref{limcond}.(a), $s=-t$.
Let us consider the vector
$$sx \boxplus x\boxplus (1-\epsilon )y\boxplus ty, $$
where $\epsilon>0$ is sufficiently small. Clearly, we have $s  \boxplus 1 \boxplus (1-\epsilon ) \boxplus t =   1 \boxplus (1-\epsilon ) =1$. Moreover,  since $|s|=|t|>1$, for all $i\in [n]$ we have $sx_i\boxplus ty_i=sx_i \boxplus x_i\boxplus (1-\epsilon )y_i\boxplus ty_i$. From equation $(1)$

$$z=sx\boxplus ty=sx \boxplus x\boxplus (1-\epsilon )y\boxplus ty. $$
Therefore $z\in \mathcal D^\infty(x,y)$. 
If $\min\{|s|,|t|\}\leq 1$ then, from Lemma \ref{limcond}, $\max\{s,t\}=1$. Suppose
that $t=1$. Then, $sx\boxplus ty=sx\boxplus 0x\boxplus y\boxplus y$ and
$s\boxplus 0\boxplus 1\boxplus 1=1$.  One can proceeds symmetrically if $s=1$. \\

$(b)$ $I\not=\emptyset$ and for all $i\in I$, $x_i=y_i$. In such a case for all $i\in I$, $z_i^{(p_{q})}=x_i=z_i$. For all $i\in I$, $sx_i=-ty_i=-tx_i$ and this implies that $s=-t$. From Lemma Lemma \ref{limcond}.(b) we have $|s|\geq 1$ and $|t|\geq 1$. Let us consider the vector
$$sx \boxplus x\boxplus (1-\epsilon )y\boxplus ty, $$
where $\epsilon>0$ is sufficiently small. Clearly, since $s=-t$, we have $s  \boxplus 1 \boxplus (1-\epsilon ) \boxplus t =   1 \boxplus (1-\epsilon ) =1$.

If $i\in I$  then, since $s =-t $ and $x_i=y_i$, one has  $sx_i \boxplus x_i\boxplus (1-\epsilon)y_i\boxplus ty_i=x_i\boxplus (1-\epsilon)y_i=x_i=y_i=z_i. $\\

 If $i\notin I$ then, since $sx_i\not=-tx_i$, $|s|\geq 1$ and $|t|\geq 1$, we have  $sx_i \boxplus x_i\boxplus (1-\epsilon)y_i\boxplus ty_i=s x_i\boxplus ty_i=z_i. $\\
We conclude that $z\in \mathcal D^\infty(x,y)$.

$(c)$   $I\not=\emptyset$ and there is some $i_0\in I$, such that $x_{i_0}\not=y_{i_0}$. We first prove that in such a  case we have
$s\not=-t$. Suppose that this is not the case. We have $sx_{i_0}\boxplus ty_{i_0}=sx_{i_0}\boxplus (-s)y_{i_0}=0$. However this implies $x_{i_0}= y_{i_0}$  that is a contradiction. Moreover, from Lemma \ref{limcond}.(c) we have $\max \{|s|,|t|\}=1$ and $\max\{s,t\}=1$. Assume for example that $t=1$.

 Suppose that $i\in I$.  We have   $t y_i=y_i=-s  x_i$ and therefore
$$s^{(p_{q })}x_i \stackrel{p_q}{+} t^{(p_{q })}y_i=s^{(p_{q })}x_i \stackrel{p_q}{-}  {s } t^{(p_{q })}x_i=
\Big(s^{(p_{q })}  \stackrel{p_q}{-}   {s }  t^{(p_{q })}\Big)x_i. $$
By hypothesis $\lim_{k\longrightarrow \infty}z_i^{(p_{q })}=z_i$. Hence the sequence $\big\{s^{(p_{q })}  \stackrel{p_q}{-} st^{(p_{q })}\Big\}_{q\in \mathbb N}$ has a limit. Let us denote
$$s'=\lim_{q\longrightarrow \infty}s^{(p_{q})}  \stackrel{p_q}{-}   {s } t^{(p_{q })}.$$
By hypothesis we have $z_i=s'x_i.$ Note that, since
$$|s^{(p_{q })}  \stackrel{p_q}{-}   {s } t^{(p_{q })}|\leq |s^{(p_{q })} | \stackrel{p_q}{+}  | {s } t^{(p_{q })}|,$$
we have taking the limit $|s'|\leq \max\{|s|, t|s|\}\leq |s|$. We consider two cases.

$(i)$ Suppose that $ s'\not=-s$. Let us consider the vector
$$ sx \boxplus s'x \boxplus 0 y\boxplus   y. $$

 If $i\in I$  then, since $sx_i =-    y_i$, we have $w_i=sx_i \boxplus s'x_i \boxplus 0 y_i\boxplus   y_i=s'x_i =z_i$.\\

 If $i\notin I$, then since $ s \not=-s'$, we have either $s=s'$ or $|s'|<|s|$. Hence, $ sx_i \boxplus s'x_i \boxplus 0 y_i\boxplus   y_i=sx_i  \boxplus   y_i=z_i$.

It follows that $z=sx \boxplus s'x \boxplus 0 y\boxplus   y.$ Moreover, since $-s\not=t=1$ and $|s'|\leq |s|$, we have $s\boxplus s'\boxplus 0\boxplus 1=1$. Hence $z\in \mathcal D^\infty(x,y)$. \\

$(ii)$ Suppose now that $ s' =-s$. Let us consider the vector
$$ sx \boxplus 0 x \boxplus  y\boxplus   y. $$

If $i\in I$  then, since $sx_i =-    y_i$, we have $sx_i \boxplus 0x_i \boxplus   y_i\boxplus   y_i=y_i$. However since $y_i=-sx_i$, we deduce that $y_i=s'x_i=z_i$.\\

If
 $i\notin I$ then,  $sx_i\not=-y_i$. Thus, we have $ sx_i \boxplus 0 x_i \boxplus   y_i\boxplus   y_i=sx_i  \boxplus   y_i=z_i$. Hence, we conclude that $z=sx \boxplus 0 x \boxplus   y\boxplus   y.$ Moreover, {since $\max\{|s|,|t|\}=1$},
 ${s} \boxplus 0  \boxplus   1\boxplus   1=1$.

 Consequently, we have shown that $z\in \Big \{t x\boxplus r x\boxplus s y\boxplus w y:
 t\boxplus r\boxplus s\boxplus w=1,t,r,s,w\in \Real\Big\}$. If $s=1$ the proof is symmetrical. Hence
 $$\Ls_{p\longrightarrow \infty}\mathcal D^{(p)}(x,y)\subset\Big \{t x\boxplus r x\boxplus s y\boxplus w y:
 t\boxplus r\boxplus s\boxplus w=1,t,r,s,w\in \Real\Big\}. \quad \Box$$

\begin{center}
{\scriptsize % This is a LaTeX picture output by TeXCAD.
% File name: [Clipboard].
% Version of TeXCAD: 4.51
% Reference / build: 27-Nov-2018 (rev. a75)
% For new versions, check: http://texcad.sf.net/
% Options on the following lines.
%\grade{\on}
%\emlines{\off}
%\epic{\off}
%\beziermacro{\on}
%\reduce{\on}
%\snapping{\off}
%\pvinsert{% Your \input, \def, etc. here}
%\quality{8.000}
%\graddiff{0.005}
%\snapasp{1}
%\zoom{4.0000}
\unitlength 0.37mm % = 2.845pt
\linethickness{0.37pt}
\ifx\plotpoint\undefined\newsavebox{\plotpoint}\fi % GNUPLOT compatibility
\begin{picture}(248.25,147)(0,0)
\put(0,11.5){\makebox(0,0)[cc]{$p\in) 0,\infty[$}}
\put(118.5,46.5){\makebox(0,0)[cc]{$p=0$}}
\put(248.25,67.75){\makebox(0,0)[cc]{$x_1$}}
\put(136.25,147){\makebox(0,0)[cc]{$x_2$}}
\put(133.5,73){\makebox(0,0)[cc]{$0$}}
\put(137,0){\makebox(0,0)[cc]
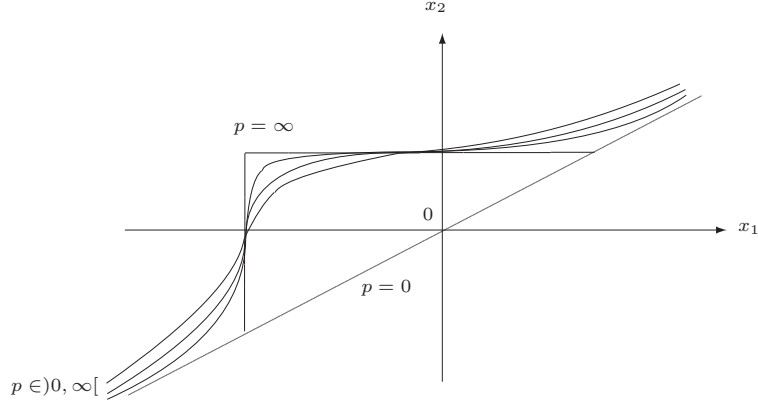
{{\bf Fig. 8:} Limit of a sequence of $\varphi_p$-lines.}}
\qbezier(18.75,8.75)(68.125,38.625)(68,64)
\qbezier(225.5,117.5)(182.5,95.25)(123.5,95)
\qbezier(123.5,95)(68.125,94.75)(68.25,64.5)
\qbezier(18.75,6.75)(68,29.125)(68.25,64)
\qbezier(68.25,64)(69.375,84.75)(74,88.5)
\qbezier(74.25,88.5)(74.625,94.625)(123.5,95.25)
\qbezier(123.5,95.25)(196.875,95.125)(225.75,115.5)
\qbezier(18.5,12.5)(64.375,44.625)(67.75,64.25)
\qbezier(67.75,64.25)(73.375,76)(78.5,80.75)
\qbezier(78.75,80.75)(82.125,86.75)(123,94.75)
\qbezier(123,94.75)(175.75,97.875)(223.5,119.5)
%\emline(26.25,8.25)(231.25,115.25)
\multiput(26.25,8.25)(.064627994956,.033732660782){3172}{\line(1,0){.064627994956}}
%\end
\put(138.5,13){\vector(0,1){124.75}}
%\emline(193,95)(68,94.75)
\multiput(193,95)(-15.625,-.03125){8}{\line(-1,0){15.625}}
%\end
%\emline(68,94.75)(67.75,31)
\multiput(68,94.75)(-.03125,-7.96875){8}{\line(0,-1){7.96875}}
%\end
\put(74.75,103.5){\makebox(0,0)[cc]{$p=\infty$}}
\put(25,67.25){\vector(1,0){215.25}}
\end{picture}
}
\end{center}
   In the next example, it is shown that the formalism described above allows us to define subsets of $\mathbb{R}^n$ that are unbounded.

   \begin{expl}\label{construct}Suppose that $n=3$ and let us consider the points $x=(3,-2, 1 )$ and $y=(1,-1, 1 )$. Suppose that $t=1$, $r=\delta=-w$ and $s=0$ where $\delta$ is a real number. Clearly, $t\boxplus r \boxplus s \boxplus w=1\boxplus \delta \boxplus(-\delta)\boxplus 0=1$. It follows that
   $$z_\delta:=tx\boxplus rx \boxplus s y\boxplus wy= t\begin{pmatrix}
   3\\-2\\ 1
   \end{pmatrix}\boxplus \delta\begin{pmatrix}
   3\\-2\\ 1
   \end{pmatrix}\boxplus (-\delta)\begin{pmatrix}
   1\\-1\\ 1
   \end{pmatrix}\boxplus 0 \begin{pmatrix}
   1\\-1\\ 1
   \end{pmatrix}= \begin{pmatrix}
   3\delta\\-2\delta\\ 1
   \end{pmatrix}\in \mathcal D^\infty(x,y).$$
   It follows that 
   $$\lim_{|\delta|\longrightarrow \infty}\|z_\delta\|=+\infty.$$
   \end{expl}
 
 In the following it is shown that $\mathcal D^\infty(x,y)$ is path connected.

 \begin{prop}Let $x,y\in \Real^n$ and suppose that $x\not=y$. Then $\mathcal D^\infty(x,y)$ is path-connected.

 \end{prop}
 {\bf Proof:} Suppose that $u,v\in \mathcal D^\infty(x,y)$ with $u\not=v$. Then, there exists a sequence $\{(u^{(p)}, v^{(p)})\}_{p\in \mathbb N}$ such that for all $p\in \mathbb N$ $(u^{(p)}, v^{(p)})\in \mathcal D^{(p)}(x,y)\times \mathcal D^{(p)}(x,y)$ and $\lim_{p\longrightarrow \infty}(u^{(p)}, v^{(p)})=(u,v)$. Let $\mathcal I(u,v)=\{i \in [n]: u_iv_i<0\}$. Since $\lim_{p\longrightarrow \infty}(u^{(p)}, v^{(p)})=(u,v)$, there is some $p_0\in \mathbb N$ such that for all $p\geq p_0$,  $\mathcal I(u,v)=\mathcal I(u^{(p)}, v^{(p)})$. For all $i\in \mathcal I(u,v)$ and all $p\geq p_0$, there is a unique point  $\gamma_i^{(p)}\in Co^p(x,y) $ such that $\{\gamma_i^{(p)}\}=H_i\cap \mathcal D^{(p)}(x,y)$  where, for all $i\in [n]$, $H_i=\{x\in \mathbb R^n:x_i=0\}$.
 For all $i\in \mathcal I(u,v)$, $$\Ls_{p\longrightarrow \infty}\big(H_i\cap \mathcal D^{(p)}(x,y)\big )\subset H_i\cap \mathcal D^\infty(x,y)$$ that is a  nonempty subset of $\mathcal D^\infty(x,y)$. Let $\gamma_i\in \Ls_{p\longrightarrow \infty}\big (H_i\cap \mathcal D^{(p)}(x,y)\big )$. This implies that there is an increasing sequence $\{p_q\}$ such that $$\lim_{q\longrightarrow \infty} \gamma_i^{(p_q)} =\gamma_i \in H_i\cap \mathcal D^\infty(x,y).$$ Let us denote $n(u,v)=\mathrm{Card} \,\mathcal I(u,v) $. For all $p\geq p_0$ one can find a sequence $\{i_k\}_{k\in [n(u,v)-1]}$ such that $$\gamma_{i_k}^{(p_q)} \cdot \gamma_{i_{k+1}}^{(p_q)} \geq 0.$$ It follows that
 $$\Lim_{q\longrightarrow \infty}Co^{p_q}\Big(\gamma_{i_k}^{(p_q)} ,\gamma_{i_{k+1}}^{(p_q)} \Big) =Co^\infty\Big (\gamma_{i_k}^{ } ,\gamma_{i_{k+1}}^{ } \Big).$$
 However, \begin{align*}\Ls_{p\longrightarrow \infty}Co^p(u^{(p)}, v^{(p)})&=\bigcup_{k\in [{n(u,v)-1}]} \Ls_{q\longrightarrow \infty}Co^{p_q}\Big(\gamma_{i_k}^{(p_q)} ,\gamma_{i_{k+1}}^{(p_q)} \Big)\\&=\bigcup_{k\in [{n(u,v)-1}]} \Lim_{q\longrightarrow \infty}Co^{p_q}\Big(\gamma_{i_k}^{(p_q)} ,\gamma_{i_{k+1}}^{(p_q)} \Big)\\&=\bigcup_{k\in [{n(u,v)-1}]}Co^\infty\Big (\gamma_{i_k}^{ } ,\gamma_{i_{k+1}}^{ } \Big).\end{align*}
 Now, since any set $Co^\infty\Big (\gamma_{i_k}^{ } ,\gamma_{i_{k+1}}^{ } \Big) $ is path-connected, it follows that
 there exists a continuous map $\psi: [0,1]\longrightarrow \mathcal D^\infty(x,y)$ such that $\psi(0)=u$ and $\psi(1)=v$. $\Box$\\

\begin{center}{\scriptsize

\unitlength 0.4mm % = 1.138pt
\linethickness{0.4pt}
\ifx\plotpoint\undefined\newsavebox{\plotpoint}\fi % GNUPLOT compatibility
\begin{picture}(223.25,147)(0,0)
\put(223.25,67.75){\makebox(0,0)[cc]{$x_1$}}
\put(111.25,147){\makebox(0,0)[cc]{$x_2$}}
\put(103,71){\makebox(0,0)[cc]{$0$}}
\put(112,0){\makebox(0,0)[cc]
{{\bf Fig. 9:} $\digamma$-line spanned from $x$ and $y$.}}
\put(107,11.75){\vector(0,1){124.75}}
\put(0,67.25){\vector(1,0){215.25}}
%\emline(168,94.75)(197.5,107.5)
\multiput(168,94.75)(.194078947,.083881579){152}{\line(1,0){.194078947}}
%\end
\put(186.5,103){\circle*{1.581}}
\put(42.25,77){\circle*{1.5}}
%\emline(41.75,39.75)(15.25,24.25)
\multiput(41.75,39.75)(-.144021739,-.08423913){184}{\line(-1,0){.144021739}}
%\end
\put(168,94.5){\line(-1,0){126}}
\put(42,94.5){\line(0,-1){55}}
\put(180.5,109.25){\makebox(0,0)[cc]{$x$}}
\put(31,79.25){\makebox(0,0)[cc]{$y$}}
\put(2.25,29.5){\makebox(0,0)[cc]{$\mathcal D^\infty(x,y)$}}
\end{picture}

}

\end{center}

\begin{center}
{\scriptsize% This is a LaTeX picture output by TeXCAD.
% File name: [Clipboard].
% Version of TeXCAD: 4.51
% Reference / build: 27-Nov-2018 (rev. a75)
% For new versions, check: http://texcad.sf.net/
% Options on the following lines.
%\grade{\on}
%\emlines{\off}
%\epic{\off}
%\beziermacro{\on}
%\reduce{\on}
%\snapping{\off}
%\pvinsert{% Your \input, \def, etc. here}
%\quality{8.000}
%\graddiff{0.005}
%\snapasp{1}
%\zoom{4.0000}
\unitlength 0.3mm % = 1.138pt
\linethickness{0.4pt}
\ifx\plotpoint\undefined\newsavebox{\plotpoint}\fi % GNUPLOT compatibility
\begin{picture}(215.25,207.75)(0,0)
\put(109.25,12.75){\vector(0,1){183.25}}
%\vector(30,139.75)(200,53.75)
\put(200,53.75){\vector(2,-1){.07}}\multiput(30,139.75)(.066692820714,-.033738721067){2549}{\line(1,0){.066692820714}}
%\end
%\vector(38,54)(199,158.25)
\put(199,158.25){\vector(3,2){.07}}\multiput(38,54)(.052103559871,.033737864078){3090}{\line(1,0){.052103559871}}
%\end
\put(25,131.5){\line(0,-1){116.5}}
%\emline(25,15)(14.25,2.75)
\multiput(25,15)(-.0336990596,-.0384012539){319}{\line(0,-1){.0384012539}}
%\end
\put(215.25,207){\line(-1,-1){19.25}}
%\emline(25.25,131.75)(153.75,206.75)
\multiput(25.25,131.75)(.057804768331,.033738191633){2223}{\line(1,0){.057804768331}}
%\end
%\emline(153.75,206.75)(196,187.75)
\multiput(153.75,206.75)(.0749113475,-.0336879433){564}{\line(1,0){.0749113475}}
%\end
%\dashline{1}(195.75,188.25)(25.5,15.75)
\multiput(195.68,188.18)(-.0324286,-.0328571){21}{\line(0,-1){.0328571}}
\multiput(194.318,186.8)(-.0324286,-.0328571){21}{\line(0,-1){.0328571}}
\multiput(192.956,185.42)(-.0324286,-.0328571){21}{\line(0,-1){.0328571}}
\multiput(191.594,184.04)(-.0324286,-.0328571){21}{\line(0,-1){.0328571}}
\multiput(190.232,182.66)(-.0324286,-.0328571){21}{\line(0,-1){.0328571}}
\multiput(188.87,181.28)(-.0324286,-.0328571){21}{\line(0,-1){.0328571}}
\multiput(187.508,179.9)(-.0324286,-.0328571){21}{\line(0,-1){.0328571}}
\multiput(186.146,178.52)(-.0324286,-.0328571){21}{\line(0,-1){.0328571}}
\multiput(184.784,177.14)(-.0324286,-.0328571){21}{\line(0,-1){.0328571}}
\multiput(183.422,175.76)(-.0324286,-.0328571){21}{\line(0,-1){.0328571}}
\multiput(182.06,174.38)(-.0324286,-.0328571){21}{\line(0,-1){.0328571}}
\multiput(180.698,173)(-.0324286,-.0328571){21}{\line(0,-1){.0328571}}
\multiput(179.336,171.62)(-.0324286,-.0328571){21}{\line(0,-1){.0328571}}
\multiput(177.974,170.24)(-.0324286,-.0328571){21}{\line(0,-1){.0328571}}
\multiput(176.612,168.86)(-.0324286,-.0328571){21}{\line(0,-1){.0328571}}
\multiput(175.25,167.48)(-.0324286,-.0328571){21}{\line(0,-1){.0328571}}
\multiput(173.888,166.1)(-.0324286,-.0328571){21}{\line(0,-1){.0328571}}
\multiput(172.526,164.72)(-.0324286,-.0328571){21}{\line(0,-1){.0328571}}
\multiput(171.164,163.34)(-.0324286,-.0328571){21}{\line(0,-1){.0328571}}
\multiput(169.802,161.96)(-.0324286,-.0328571){21}{\line(0,-1){.0328571}}
\multiput(168.44,160.58)(-.0324286,-.0328571){21}{\line(0,-1){.0328571}}
\multiput(167.078,159.2)(-.0324286,-.0328571){21}{\line(0,-1){.0328571}}
\multiput(165.716,157.82)(-.0324286,-.0328571){21}{\line(0,-1){.0328571}}
\multiput(164.354,156.44)(-.0324286,-.0328571){21}{\line(0,-1){.0328571}}
\multiput(162.992,155.06)(-.0324286,-.0328571){21}{\line(0,-1){.0328571}}
\multiput(161.63,153.68)(-.0324286,-.0328571){21}{\line(0,-1){.0328571}}
\multiput(160.268,152.3)(-.0324286,-.0328571){21}{\line(0,-1){.0328571}}
\multiput(158.906,150.92)(-.0324286,-.0328571){21}{\line(0,-1){.0328571}}
\multiput(157.544,149.54)(-.0324286,-.0328571){21}{\line(0,-1){.0328571}}
\multiput(156.182,148.16)(-.0324286,-.0328571){21}{\line(0,-1){.0328571}}
\multiput(154.82,146.78)(-.0324286,-.0328571){21}{\line(0,-1){.0328571}}
\multiput(153.458,145.4)(-.0324286,-.0328571){21}{\line(0,-1){.0328571}}
\multiput(152.096,144.02)(-.0324286,-.0328571){21}{\line(0,-1){.0328571}}
\multiput(150.734,142.64)(-.0324286,-.0328571){21}{\line(0,-1){.0328571}}
\multiput(149.372,141.26)(-.0324286,-.0328571){21}{\line(0,-1){.0328571}}
\multiput(148.01,139.88)(-.0324286,-.0328571){21}{\line(0,-1){.0328571}}
\multiput(146.648,138.5)(-.0324286,-.0328571){21}{\line(0,-1){.0328571}}
\multiput(145.286,137.12)(-.0324286,-.0328571){21}{\line(0,-1){.0328571}}
\multiput(143.924,135.74)(-.0324286,-.0328571){21}{\line(0,-1){.0328571}}
\multiput(142.562,134.36)(-.0324286,-.0328571){21}{\line(0,-1){.0328571}}
\multiput(141.2,132.98)(-.0324286,-.0328571){21}{\line(0,-1){.0328571}}
\multiput(139.838,131.6)(-.0324286,-.0328571){21}{\line(0,-1){.0328571}}
\multiput(138.476,130.22)(-.0324286,-.0328571){21}{\line(0,-1){.0328571}}
\multiput(137.114,128.84)(-.0324286,-.0328571){21}{\line(0,-1){.0328571}}
\multiput(135.752,127.46)(-.0324286,-.0328571){21}{\line(0,-1){.0328571}}
\multiput(134.39,126.08)(-.0324286,-.0328571){21}{\line(0,-1){.0328571}}
\multiput(133.028,124.7)(-.0324286,-.0328571){21}{\line(0,-1){.0328571}}
\multiput(131.666,123.32)(-.0324286,-.0328571){21}{\line(0,-1){.0328571}}
\multiput(130.304,121.94)(-.0324286,-.0328571){21}{\line(0,-1){.0328571}}
\multiput(128.942,120.56)(-.0324286,-.0328571){21}{\line(0,-1){.0328571}}
\multiput(127.58,119.18)(-.0324286,-.0328571){21}{\line(0,-1){.0328571}}
\multiput(126.218,117.8)(-.0324286,-.0328571){21}{\line(0,-1){.0328571}}
\multiput(124.856,116.42)(-.0324286,-.0328571){21}{\line(0,-1){.0328571}}
\multiput(123.494,115.04)(-.0324286,-.0328571){21}{\line(0,-1){.0328571}}
\multiput(122.132,113.66)(-.0324286,-.0328571){21}{\line(0,-1){.0328571}}
\multiput(120.77,112.28)(-.0324286,-.0328571){21}{\line(0,-1){.0328571}}
\multiput(119.408,110.9)(-.0324286,-.0328571){21}{\line(0,-1){.0328571}}
\multiput(118.046,109.52)(-.0324286,-.0328571){21}{\line(0,-1){.0328571}}
\multiput(116.684,108.14)(-.0324286,-.0328571){21}{\line(0,-1){.0328571}}
\multiput(115.322,106.76)(-.0324286,-.0328571){21}{\line(0,-1){.0328571}}
\multiput(113.96,105.38)(-.0324286,-.0328571){21}{\line(0,-1){.0328571}}
\multiput(112.598,104)(-.0324286,-.0328571){21}{\line(0,-1){.0328571}}
\multiput(111.236,102.62)(-.0324286,-.0328571){21}{\line(0,-1){.0328571}}
\multiput(109.874,101.24)(-.0324286,-.0328571){21}{\line(0,-1){.0328571}}
\multiput(108.512,99.86)(-.0324286,-.0328571){21}{\line(0,-1){.0328571}}
\multiput(107.15,98.48)(-.0324286,-.0328571){21}{\line(0,-1){.0328571}}
\multiput(105.788,97.1)(-.0324286,-.0328571){21}{\line(0,-1){.0328571}}
\multiput(104.426,95.72)(-.0324286,-.0328571){21}{\line(0,-1){.0328571}}
\multiput(103.064,94.34)(-.0324286,-.0328571){21}{\line(0,-1){.0328571}}
\multiput(101.702,92.96)(-.0324286,-.0328571){21}{\line(0,-1){.0328571}}
\multiput(100.34,91.58)(-.0324286,-.0328571){21}{\line(0,-1){.0328571}}
\multiput(98.978,90.2)(-.0324286,-.0328571){21}{\line(0,-1){.0328571}}
\multiput(97.616,88.82)(-.0324286,-.0328571){21}{\line(0,-1){.0328571}}
\multiput(96.254,87.44)(-.0324286,-.0328571){21}{\line(0,-1){.0328571}}
\multiput(94.892,86.06)(-.0324286,-.0328571){21}{\line(0,-1){.0328571}}
\multiput(93.53,84.68)(-.0324286,-.0328571){21}{\line(0,-1){.0328571}}
\multiput(92.168,83.3)(-.0324286,-.0328571){21}{\line(0,-1){.0328571}}
\multiput(90.806,81.92)(-.0324286,-.0328571){21}{\line(0,-1){.0328571}}
\multiput(89.444,80.54)(-.0324286,-.0328571){21}{\line(0,-1){.0328571}}
\multiput(88.082,79.16)(-.0324286,-.0328571){21}{\line(0,-1){.0328571}}
\multiput(86.72,77.78)(-.0324286,-.0328571){21}{\line(0,-1){.0328571}}
\multiput(85.358,76.4)(-.0324286,-.0328571){21}{\line(0,-1){.0328571}}
\multiput(83.996,75.02)(-.0324286,-.0328571){21}{\line(0,-1){.0328571}}
\multiput(82.634,73.64)(-.0324286,-.0328571){21}{\line(0,-1){.0328571}}
\multiput(81.272,72.26)(-.0324286,-.0328571){21}{\line(0,-1){.0328571}}
\multiput(79.91,70.88)(-.0324286,-.0328571){21}{\line(0,-1){.0328571}}
\multiput(78.548,69.5)(-.0324286,-.0328571){21}{\line(0,-1){.0328571}}
\multiput(77.186,68.12)(-.0324286,-.0328571){21}{\line(0,-1){.0328571}}
\multiput(75.824,66.74)(-.0324286,-.0328571){21}{\line(0,-1){.0328571}}
\multiput(74.462,65.36)(-.0324286,-.0328571){21}{\line(0,-1){.0328571}}
\multiput(73.1,63.98)(-.0324286,-.0328571){21}{\line(0,-1){.0328571}}
\multiput(71.738,62.6)(-.0324286,-.0328571){21}{\line(0,-1){.0328571}}
\multiput(70.376,61.22)(-.0324286,-.0328571){21}{\line(0,-1){.0328571}}
\multiput(69.014,59.84)(-.0324286,-.0328571){21}{\line(0,-1){.0328571}}
\multiput(67.652,58.46)(-.0324286,-.0328571){21}{\line(0,-1){.0328571}}
\multiput(66.29,57.08)(-.0324286,-.0328571){21}{\line(0,-1){.0328571}}
\multiput(64.928,55.7)(-.0324286,-.0328571){21}{\line(0,-1){.0328571}}
\multiput(63.566,54.32)(-.0324286,-.0328571){21}{\line(0,-1){.0328571}}
\multiput(62.204,52.94)(-.0324286,-.0328571){21}{\line(0,-1){.0328571}}
\multiput(60.842,51.56)(-.0324286,-.0328571){21}{\line(0,-1){.0328571}}
\multiput(59.48,50.18)(-.0324286,-.0328571){21}{\line(0,-1){.0328571}}
\multiput(58.118,48.8)(-.0324286,-.0328571){21}{\line(0,-1){.0328571}}
\multiput(56.756,47.42)(-.0324286,-.0328571){21}{\line(0,-1){.0328571}}
\multiput(55.394,46.04)(-.0324286,-.0328571){21}{\line(0,-1){.0328571}}
\multiput(54.032,44.66)(-.0324286,-.0328571){21}{\line(0,-1){.0328571}}
\multiput(52.67,43.28)(-.0324286,-.0328571){21}{\line(0,-1){.0328571}}
\multiput(51.308,41.9)(-.0324286,-.0328571){21}{\line(0,-1){.0328571}}
\multiput(49.946,40.52)(-.0324286,-.0328571){21}{\line(0,-1){.0328571}}
\multiput(48.584,39.14)(-.0324286,-.0328571){21}{\line(0,-1){.0328571}}
\multiput(47.222,37.76)(-.0324286,-.0328571){21}{\line(0,-1){.0328571}}
\multiput(45.86,36.38)(-.0324286,-.0328571){21}{\line(0,-1){.0328571}}
\multiput(44.498,35)(-.0324286,-.0328571){21}{\line(0,-1){.0328571}}
\multiput(43.136,33.62)(-.0324286,-.0328571){21}{\line(0,-1){.0328571}}
\multiput(41.774,32.24)(-.0324286,-.0328571){21}{\line(0,-1){.0328571}}
\multiput(40.412,30.86)(-.0324286,-.0328571){21}{\line(0,-1){.0328571}}
\multiput(39.05,29.48)(-.0324286,-.0328571){21}{\line(0,-1){.0328571}}
\multiput(37.688,28.1)(-.0324286,-.0328571){21}{\line(0,-1){.0328571}}
\multiput(36.326,26.72)(-.0324286,-.0328571){21}{\line(0,-1){.0328571}}
\multiput(34.964,25.34)(-.0324286,-.0328571){21}{\line(0,-1){.0328571}}
\multiput(33.602,23.96)(-.0324286,-.0328571){21}{\line(0,-1){.0328571}}
\multiput(32.24,22.58)(-.0324286,-.0328571){21}{\line(0,-1){.0328571}}
\multiput(30.878,21.2)(-.0324286,-.0328571){21}{\line(0,-1){.0328571}}
\multiput(29.516,19.82)(-.0324286,-.0328571){21}{\line(0,-1){.0328571}}
\multiput(28.154,18.44)(-.0324286,-.0328571){21}{\line(0,-1){.0328571}}
\multiput(26.792,17.06)(-.0324286,-.0328571){21}{\line(0,-1){.0328571}}
%\end
\put(172,198.75){\circle*{1.581}}
\put(207.75,199.25){\circle*{1.581}}
\put(25,53){\circle*{1.581}}
\put(93.75,171.25){\circle*{1.581}}
\put(24.5,74){\circle*{1.803}}
\put(209,49.5){\makebox(0,0)[cc]{$x_1$}}
\put(209.5,161.75){\makebox(0,0)[cc]{$x_2$}}
\put(109,207.25){\makebox(0,0)[cc]{$x_3$}}
\put(105,108.5){\makebox(0,0)[cc]{$0$}}
\put(11.5,78){\makebox(0,0)[cc]{$\gamma_3$}}
\put(180,207.75){\makebox(0,0)[cc]{$\gamma_2$}}
\put(86.5,177.25){\makebox(0,0)[cc]{$\gamma_1$}}
\put(202,206.25){\makebox(0,0)[cc]{$x$}}
\put(14.5,52){\makebox(0,0)[cc]{$y$}}
\put(0,10.5){\makebox(0,0)[cc]{$\mathcal D^\infty(x,y)$}}
\put(90.25,0){\makebox(0,0)[cc]
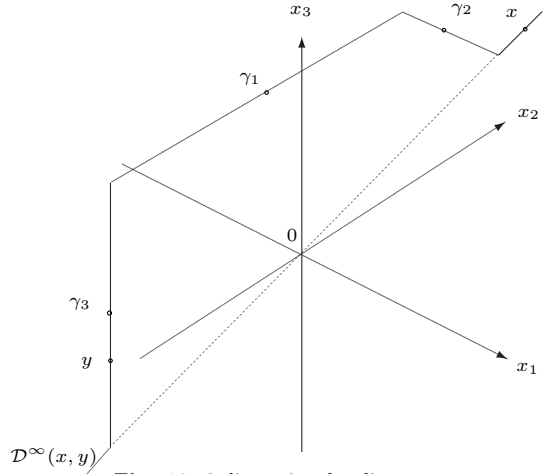
{{\bf Fig. 10:}  3-dimensional $\digamma$-line.}}
\end{picture}

}
\end{center}
 Figure 10 depicts the case of a 3-dimensional $\digamma$-line passing trough $x$ and $y$. For $i=1,2,3$, $\gamma_i$ is the intersection between the $\digamma$-line and the 2-dimensions plane with equation $x_i=0.$
\begin{cor}\label{Portion} For all $x,y\in \Real^n$ with $x\not=y$,
 $$Co^\infty(x,y)\subset\mathcal D^\infty(x,y).  $$
 
 \end{cor}
 {\bf Proof:} From \cite{b15}\begin{equation*}Co^\infty(x,y):=\Big \{t x\boxplus r x\boxplus s y\boxplus w y:
t\boxplus r\boxplus s\boxplus w=1,t,r,s,w\in \Real_+\Big\}.\end{equation*}The result is then immediate, since $Co^\infty(x,y)$ is the set of all the $\mathbb B$-convex combinations such that $t,r,s,w\geq 0.$ $\Box$\\

Corollary \ref{Portion} does not imply that $\mathcal D^\infty$ is {idempotent symmetric convex}. In general, this is not true. This can be seen in  Figure 10, where $x\boxminus y$ and $y\boxminus x$ belong to $\mathcal D^\infty(x,y)$. However,  $(x\boxminus y)\boxplus( y\boxminus x)=0\notin \mathcal D^\infty{(x,y)}.$

\begin{prop} Let $Y$ be an idempotent symmetric subspace of $\Real^n$. Then, for all $x,y\in \Real^n$ with $x\not=y$ we have $\mathcal D^\infty(x,y)\subset Y$.
\end{prop}
{\bf Proof:} We have shown that
$$ \mathcal D^\infty(x,y)=\Big \{t x\boxplus r x\boxplus s y\boxplus w y:
 t\boxplus r\boxplus s\boxplus w=1,t,r,s,w\in \Real\Big\}.$$
 From Proposition \ref{comb} the result is then immediate. $\Box$\\
 
 Notice that since $(\Real^n, \boxplus , \cdot)$ and $ ( \widetilde{\mathbb M}^n, \widetilde{\boxplus} , \otimes)$ are isomorphic one can deduce the algebraic form of a line passing trough $z$ and $u$. The multiplicative neutral element of $\cdot$ is replaced with $0=1_{ \widetilde{\mathbb M}}$ and we have 
 \begin{equation}
{\widetilde{ \mathcal D}}^\infty (z,u)=\Big \{t z \, \widetilde \boxplus \, r z\,\widetilde\boxplus \,s u\, \widetilde \boxplus \, w u:
 t\,\widetilde \boxplus\, r\widetilde\boxplus \,  s\, \widetilde\boxplus \,w=0,t,r,s,w\in \widetilde{\mathbb M}\Big\}. 
 \end{equation}

\subsection{Half-Lines and Inclusion Properties}

In the next statement, it is shown that, for all $x,y\in \Real^n$ with $x\not=y$, the line $\mathcal D^\infty(x,y)$ contains two half lines spanned from two specific points. $\mathcal D^\infty(x,y)$ is called the $\digamma$-line spanned by $x$ and $y$.
\begin{prop}\label{linex}Let $x,y\in \Real^n$ and suppose that $x\not=y$. Let $I=\{i\in [n]: x_i\not=y_i\}$ and $J=[n]\backslash I$. We have the following inclusions:\\
$(a)$ If $J=\emptyset$ then $$\{t  (y\boxminus x): t\in [1,\infty[\,\} \subset \mathcal D^\infty(x,y)\quad \text{ and }\quad
   \{t  (x\boxminus  y): t\in [1,\infty[\, \}\subset \mathcal D^\infty(x,y).$$
 $(b)$ If $J\not=\emptyset$ then $$x_{[J]}+\Big\{t  (y \boxminus x ): t\in [1,\infty[\,\Big\} \subset \mathcal D^\infty(x,y)\quad \text{ and }\quad x_{[J]}+\Big \{t  (x \boxminus y ): t\in [1,\infty[\,\Big \} \subset \mathcal D^\infty(x,y),$$
 where $x_{[J]}=\sum\limits_{k\in J}x_ke_k=\sum\limits_{k\in J}y_ke_k=y_{[J]}.$
\end{prop}
{\bf Proof:} $(a)$ We first note that since $J=\emptyset$ and $I=[n]$
$$x\boxminus y=x\boxplus x\boxplus 0 y\boxplus (-y)\in \Big \{t x\boxplus r x\boxplus s y\boxplus w y:
 t\boxplus r\boxplus s\boxplus w=1,t,r,s,w\in \Real\Big\}.$$
 For all $t\geq 1$, let us consider the vector
$$z=tx   \boxplus (-t)y.$$
Since {$t\geq 1$}, we have
$$z=tx \boxplus x\boxplus 0y  \boxplus (-t)y.$$
Moreover $ t  \boxplus 1\boxplus 0  \boxplus (-t)=1$. Clearly $z\in  \Big \{t x\boxplus r x\boxplus s y\boxplus w y:
 t\boxplus r\boxplus s\boxplus w=1,t,r,s,w\in \Real\Big\} $, which proves the first part of $(a)$. The proof of the second part is symmetrical swapping $x$ and $y$. $(b)$ Let us denote  $\Real_{[I]}=\{\sum_{k\in I}x_ke_k: x\in \Real^n\}$. By hypothesis $x,y\in x_{[J]}+\Real_{[I]}$. Hence 
 $\mathcal D^\infty(x,y)\subset  x_{[J]}+\Real_{[I]}.$ Let $n_I=\mathrm{Card}(I)$ and let $\mathcal D_{[I]}^\infty(x,y)$ be the canonical decomposition of $\mathcal D^\infty(x,y)$ over $\Real_{[I]}$. By definition for all $i\in I$ we have $x_i\not=y_i$. Since $\Real_{[I]}$ is isomorphic to its restriction projection onto $\Real^{n_I}$ we deduce from $(a)$ that 
   $$\big\{t  (x_{[I]}\boxminus  y_{[I]}): t\in [1,\infty[\,\big \}\subset \mathcal D_{[I]}^\infty(x,y).$$
However $\mathcal D^\infty(x,y)= x_{[J]}+\mathcal D_{[I]}^\infty(x,y)$. Moreover $x_{[J]}\boxminus  y_{[J]}=0.$ Hence
   \begin{align*}x_{[J]}+\big\{t  (x  \boxminus  y): t\in [1,\infty[\,\big \}&=x_{[J]}+\big\{t  (x_{[I]}\boxminus  y_{[I]}): t\in [1,\infty[\,\big \}\\&\subset x_{[J]}+ \mathcal D_{[I]}^\infty(x,y)\\&\subset\mathcal D^\infty(x,y).\end{align*}
The proof of the second part of the statement is symmetrical. $\Box$\\

Notice that Proposition \ref{linex}.(b)  can be deduced from  Proposition \ref{linex}.(a) by considering the canonical isomorphism $\psi_I$   defined from the set $\{z\in \Real^n: z_j=0,\; j\notin I\} $ to $\Real^m$ where $m=\mathrm{Card}(I)$.

\begin{center}{\scriptsize

\unitlength 0.4mm % = 1.138pt
\linethickness{0.4pt}
\ifx\plotpoint\undefined\newsavebox{\plotpoint}\fi % GNUPLOT compatibility

\unitlength 0.4mm % = 1.138pt
\linethickness{0.4pt}
\ifx\plotpoint\undefined\newsavebox{\plotpoint}\fi % GNUPLOT compatibility
\begin{picture}(234.412,154.35)(0,0)
\put(234.412,71.138){\makebox(0,0)[cc]{$x_1$}}
\put(116.813,154.35){\makebox(0,0)[cc]{$x_2$}}
\put(108.15,74.55){\makebox(0,0)[cc]{$0$}}
\put(117.6,0){\makebox(0,0)[cc]
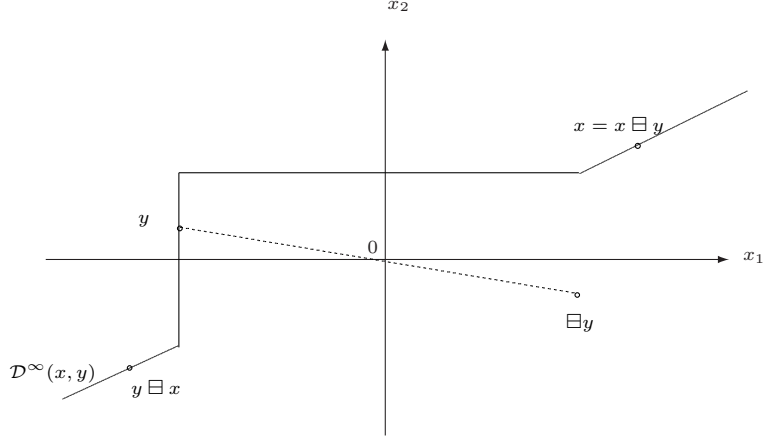
{{\bf Fig. 11:} Inclusion of half-lines.}}
\put(112.35,12.338){\vector(0,1){130.987}}
\put(0,70.613){\vector(1,0){226.012}}
\put(195.825,108.15){\circle*{1.66}}
\put(44.363,80.85){\circle*{1.575}}
\put(176.4,99.225){\line(-1,0){132.3}}
\put(44.1,99.225){\line(0,-1){57.75}}
\put(189.525,114.712){\makebox(0,0)[cc]{$x=x\boxminus y$}}
\put(36.525,27.712){\makebox(0,0)[cc]{$y\boxminus x$}}
\put(32.55,83.212){\makebox(0,0)[cc]{$y$}}
\put(2.363,32.975){\makebox(0,0)[cc]{$\mathcal D^\infty(x,y)$}}
%\emline(176.625,99)(232.125,126.375)
\multiput(176.625,99)(.1707692308,.0842307692){325}{\line(1,0){.1707692308}}
%\end
%\emline(43.875,42)(5.625,24.375)
\multiput(43.875,42)(-.183014354,-.084330144){209}{\line(-1,0){.183014354}}
%\end
%\dashline{1}(44.625,81.375)(175.125,59.625)
\multiput(44.449,81.199)(.4906,-.08177){2}{\line(1,0){.4906}}
\multiput(46.412,80.872)(.4906,-.08177){2}{\line(1,0){.4906}}
\multiput(48.374,80.545)(.4906,-.08177){2}{\line(1,0){.4906}}
\multiput(50.336,80.218)(.4906,-.08177){2}{\line(1,0){.4906}}
\multiput(52.299,79.891)(.4906,-.08177){2}{\line(1,0){.4906}}
\multiput(54.261,79.564)(.4906,-.08177){2}{\line(1,0){.4906}}
\multiput(56.224,79.237)(.4906,-.08177){2}{\line(1,0){.4906}}
\multiput(58.186,78.91)(.4906,-.08177){2}{\line(1,0){.4906}}
\multiput(60.149,78.583)(.4906,-.08177){2}{\line(1,0){.4906}}
\multiput(62.111,78.256)(.4906,-.08177){2}{\line(1,0){.4906}}
\multiput(64.073,77.929)(.4906,-.08177){2}{\line(1,0){.4906}}
\multiput(66.036,77.602)(.4906,-.08177){2}{\line(1,0){.4906}}
\multiput(67.998,77.274)(.4906,-.08177){2}{\line(1,0){.4906}}
\multiput(69.961,76.947)(.4906,-.08177){2}{\line(1,0){.4906}}
\multiput(71.923,76.62)(.4906,-.08177){2}{\line(1,0){.4906}}
\multiput(73.885,76.293)(.4906,-.08177){2}{\line(1,0){.4906}}
\multiput(75.848,75.966)(.4906,-.08177){2}{\line(1,0){.4906}}
\multiput(77.81,75.639)(.4906,-.08177){2}{\line(1,0){.4906}}
\multiput(79.773,75.312)(.4906,-.08177){2}{\line(1,0){.4906}}
\multiput(81.735,74.985)(.4906,-.08177){2}{\line(1,0){.4906}}
\multiput(83.697,74.658)(.4906,-.08177){2}{\line(1,0){.4906}}
\multiput(85.66,74.331)(.4906,-.08177){2}{\line(1,0){.4906}}
\multiput(87.622,74.004)(.4906,-.08177){2}{\line(1,0){.4906}}
\multiput(89.585,73.677)(.4906,-.08177){2}{\line(1,0){.4906}}
\multiput(91.547,73.35)(.4906,-.08177){2}{\line(1,0){.4906}}
\multiput(93.509,73.023)(.4906,-.08177){2}{\line(1,0){.4906}}
\multiput(95.472,72.696)(.4906,-.08177){2}{\line(1,0){.4906}}
\multiput(97.434,72.368)(.4906,-.08177){2}{\line(1,0){.4906}}
\multiput(99.397,72.041)(.4906,-.08177){2}{\line(1,0){.4906}}
\multiput(101.359,71.714)(.4906,-.08177){2}{\line(1,0){.4906}}
\multiput(103.321,71.387)(.4906,-.08177){2}{\line(1,0){.4906}}
\multiput(105.284,71.06)(.4906,-.08177){2}{\line(1,0){.4906}}
\multiput(107.246,70.733)(.4906,-.08177){2}{\line(1,0){.4906}}
\multiput(109.209,70.406)(.4906,-.08177){2}{\line(1,0){.4906}}
\multiput(111.171,70.079)(.4906,-.08177){2}{\line(1,0){.4906}}
\multiput(113.133,69.752)(.4906,-.08177){2}{\line(1,0){.4906}}
\multiput(115.096,69.425)(.4906,-.08177){2}{\line(1,0){.4906}}
\multiput(117.058,69.098)(.4906,-.08177){2}{\line(1,0){.4906}}
\multiput(119.021,68.771)(.4906,-.08177){2}{\line(1,0){.4906}}
\multiput(120.983,68.444)(.4906,-.08177){2}{\line(1,0){.4906}}
\multiput(122.946,68.117)(.4906,-.08177){2}{\line(1,0){.4906}}
\multiput(124.908,67.789)(.4906,-.08177){2}{\line(1,0){.4906}}
\multiput(126.87,67.462)(.4906,-.08177){2}{\line(1,0){.4906}}
\multiput(128.833,67.135)(.4906,-.08177){2}{\line(1,0){.4906}}
\multiput(130.795,66.808)(.4906,-.08177){2}{\line(1,0){.4906}}
\multiput(132.758,66.481)(.4906,-.08177){2}{\line(1,0){.4906}}
\multiput(134.72,66.154)(.4906,-.08177){2}{\line(1,0){.4906}}
\multiput(136.682,65.827)(.4906,-.08177){2}{\line(1,0){.4906}}
\multiput(138.645,65.5)(.4906,-.08177){2}{\line(1,0){.4906}}
\multiput(140.607,65.173)(.4906,-.08177){2}{\line(1,0){.4906}}
\multiput(142.57,64.846)(.4906,-.08177){2}{\line(1,0){.4906}}
\multiput(144.532,64.519)(.4906,-.08177){2}{\line(1,0){.4906}}
\multiput(146.494,64.192)(.4906,-.08177){2}{\line(1,0){.4906}}
\multiput(148.457,63.865)(.4906,-.08177){2}{\line(1,0){.4906}}
\multiput(150.419,63.538)(.4906,-.08177){2}{\line(1,0){.4906}}
\multiput(152.382,63.211)(.4906,-.08177){2}{\line(1,0){.4906}}
\multiput(154.344,62.883)(.4906,-.08177){2}{\line(1,0){.4906}}
\multiput(156.306,62.556)(.4906,-.08177){2}{\line(1,0){.4906}}
\multiput(158.269,62.229)(.4906,-.08177){2}{\line(1,0){.4906}}
\multiput(160.231,61.902)(.4906,-.08177){2}{\line(1,0){.4906}}
\multiput(162.194,61.575)(.4906,-.08177){2}{\line(1,0){.4906}}
\multiput(164.156,61.248)(.4906,-.08177){2}{\line(1,0){.4906}}
\multiput(166.118,60.921)(.4906,-.08177){2}{\line(1,0){.4906}}
\multiput(168.081,60.594)(.4906,-.08177){2}{\line(1,0){.4906}}
\multiput(170.043,60.267)(.4906,-.08177){2}{\line(1,0){.4906}}
\multiput(172.006,59.94)(.4906,-.08177){2}{\line(1,0){.4906}}
\multiput(173.968,59.613)(.4906,-.08177){2}{\line(1,0){.4906}}
%\end
\put(175.875,58.875){\circle*{1.5}}
\put(176.625,49.5){\makebox(0,0)[cc]{$\boxminus y$}}
\put(27.75,34.70){\circle*{1.5}}
\end{picture}

}

\end{center}

Let $x,y\in \Real^n$. The half lines $\mathcal D_+(x,y)=\big\{t  (x \boxminus y ): t\in [1,\infty[\,\big \}$ and 
$\mathcal D_-(x,y)=\big\{t  (y \boxminus x ): t\in [1,\infty[\,\big \}$ are respectively called {\bf the  upper and lower half-lines components} of the $\digamma$-line  $\mathcal D^\infty(x,y)$. 

  \begin{expl}Let us consider the case proposed in \ref{construct} with $x=(3,-2, 1 )$ and $y=(1,-1, 1 )$. In such a case $I=\{1,2\}$ and $J=\{3\}$.   We have:

   $$x\boxminus y= \begin{pmatrix}
   3\\-2\\ 1
   \end{pmatrix}\boxminus   \begin{pmatrix}
   1\\-1\\ 1
   \end{pmatrix}=\begin{pmatrix} 
   3 \\-2 \\ 0
   \end{pmatrix} $$
   and 
  $$x_{\{3\}}=y_{\{3\}}=\begin{pmatrix}
   0\\0\\ 1
   \end{pmatrix}.$$ It follows that $$D^\infty(x,y)\supset \begin{pmatrix}
   0\\0\\ 1
   \end{pmatrix}+\left\{t\begin{pmatrix} 
   3 \\-2 \\ 0
   \end{pmatrix}, t\in \Real\right\}=\left\{ \begin{pmatrix} 
   3 t\\-2 t\\ 1
   \end{pmatrix}: t\in \Real\right\} = \mathcal D_+(x,y)$$
   and
   $$D^\infty(x,y)\supset \begin{pmatrix}
   0\\0\\ 1
   \end{pmatrix}+\left\{t\begin{pmatrix} 
   -3 \\ 2 \\ 0
   \end{pmatrix}, t\in \Real\right\}=\left\{ \begin{pmatrix} 
   -3 t\\2 t\\ 1
   \end{pmatrix}: t\in \Real\right\}=  \mathcal D_-(x,y). $$
   
   \end{expl}

\subsection{Parallel Lines }

From Proposition \ref{linex}, any $\digamma$-lines contains two half-lines components. Let $\mathcal D^\infty(x,y)$ and $\mathcal E^\infty(u,v)$ be two $\digamma$-lines respectively spanned from $x$ and $y$ and  $u$ and $v$. We say that  $\mathcal E^\infty (u,v)$  and $\mathcal D^\infty (x,y)$ are {\bf $\digamma$-parallel} 
  if $x\boxminus y \propto u\boxminus v$. Equivalently, this means that $x\boxminus y$ and $ u\boxminus v$ are collinear, that is there is some $\alpha\in \Real\backslash\{0\}$ such that $(x\boxminus y )=\alpha (u\boxminus v)$. From Proposition \ref{linex}, there exist two half lines components $\mathcal E_-(u,v)$ and $\mathcal E_+(u,v)$  of  $\mathcal E^\infty(u,v)$ and
 two half lines components $\mathcal D_-(x,y)$ and $\mathcal D_+(x,y)$ of $\mathcal D^\infty(x,y)$ such that  $\mathcal E_+(u,v)\supset (\subset) \;\mathcal D_+(x,y)$ and $\mathcal E_-(u,v)\supset (\subset) \; \mathcal D_-(x,y)$.

 \begin{center}{\scriptsize

\unitlength 0.4mm % = 1.138pt
\linethickness{0.4pt}
\ifx\plotpoint\undefined\newsavebox{\plotpoint}\fi % GNUPLOT compatibility
\begin{picture}(244.162,154.35)(0,0)
\put(244.162,71.138){\makebox(0,0)[cc]{$x_1$}}
\put(126.563,154.35){\makebox(0,0)[cc]{$x_2$}}
\put(117.9,74.55){\makebox(0,0)[cc]{$0$}}
\put(127.35,0){\makebox(0,0)[cc]
{{\bf Fig. 12:} Parallel $\digamma$-lines.}}
\put(122.1,12.338){\vector(0,1){130.987}}
\put(9.75,70.613){\vector(1,0){226.012}}
\put(186.15,99.225){\line(-1,0){132.3}}
\put(53.85,99.225){\line(0,-1){57.75}}
\put(23.863,36.975){\makebox(0,0)[cc]{$\mathcal D^\infty$}}
\put(70.863,91.225){\makebox(0,0)[cc]{$\mathcal E_\infty$}}
\put(47.113,28.475){\makebox(0,0)[cc]{$\mathcal E_-$}}
\put(227.613,105.975){\makebox(0,0)[cc]{$\mathcal D_+$}}
\put(25.863,19.475){\makebox(0,0)[cc]{$\mathcal D_-$}}
\put(206.363,96.975){\makebox(0,0)[cc]{$\mathcal E_+$}}
\put(206.25,108.25){\line(-1,0){169.5}}
\put(36.75,108.25){\line(0,-1){73.75}}
%\emline(53.75,41.75)(0,18.25)
\multiput(53.75,41.75)(-.07711621234,-.03371592539){697}{\line(-1,0){.07711621234}}
%\end
%\emline(185.5,99)(232,119)
\multiput(185.5,99)(.0784148398,.0337268128){593}{\line(1,0){.0784148398}}
%\end
%\dashline{1}(185.5,99)(53.25,41.25)
\multiput(185.43,98.93)(-.076006,-.03319){12}{\line(-1,0){.076006}}
\multiput(183.606,98.133)(-.076006,-.03319){12}{\line(-1,0){.076006}}
\multiput(181.781,97.337)(-.076006,-.03319){12}{\line(-1,0){.076006}}
\multiput(179.957,96.54)(-.076006,-.03319){12}{\line(-1,0){.076006}}
\multiput(178.133,95.744)(-.076006,-.03319){12}{\line(-1,0){.076006}}
\multiput(176.309,94.947)(-.076006,-.03319){12}{\line(-1,0){.076006}}
\multiput(174.485,94.15)(-.076006,-.03319){12}{\line(-1,0){.076006}}
\multiput(172.661,93.354)(-.076006,-.03319){12}{\line(-1,0){.076006}}
\multiput(170.837,92.557)(-.076006,-.03319){12}{\line(-1,0){.076006}}
\multiput(169.012,91.761)(-.076006,-.03319){12}{\line(-1,0){.076006}}
\multiput(167.188,90.964)(-.076006,-.03319){12}{\line(-1,0){.076006}}
\multiput(165.364,90.168)(-.076006,-.03319){12}{\line(-1,0){.076006}}
\multiput(163.54,89.371)(-.076006,-.03319){12}{\line(-1,0){.076006}}
\multiput(161.716,88.575)(-.076006,-.03319){12}{\line(-1,0){.076006}}
\multiput(159.892,87.778)(-.076006,-.03319){12}{\line(-1,0){.076006}}
\multiput(158.068,86.981)(-.076006,-.03319){12}{\line(-1,0){.076006}}
\multiput(156.244,86.185)(-.076006,-.03319){12}{\line(-1,0){.076006}}
\multiput(154.419,85.388)(-.076006,-.03319){12}{\line(-1,0){.076006}}
\multiput(152.595,84.592)(-.076006,-.03319){12}{\line(-1,0){.076006}}
\multiput(150.771,83.795)(-.076006,-.03319){12}{\line(-1,0){.076006}}
\multiput(148.947,82.999)(-.076006,-.03319){12}{\line(-1,0){.076006}}
\multiput(147.123,82.202)(-.076006,-.03319){12}{\line(-1,0){.076006}}
\multiput(145.299,81.406)(-.076006,-.03319){12}{\line(-1,0){.076006}}
\multiput(143.475,80.609)(-.076006,-.03319){12}{\line(-1,0){.076006}}
\multiput(141.65,79.812)(-.076006,-.03319){12}{\line(-1,0){.076006}}
\multiput(139.826,79.016)(-.076006,-.03319){12}{\line(-1,0){.076006}}
\multiput(138.002,78.219)(-.076006,-.03319){12}{\line(-1,0){.076006}}
\multiput(136.178,77.423)(-.076006,-.03319){12}{\line(-1,0){.076006}}
\multiput(134.354,76.626)(-.076006,-.03319){12}{\line(-1,0){.076006}}
\multiput(132.53,75.83)(-.076006,-.03319){12}{\line(-1,0){.076006}}
\multiput(130.706,75.033)(-.076006,-.03319){12}{\line(-1,0){.076006}}
\multiput(128.881,74.237)(-.076006,-.03319){12}{\line(-1,0){.076006}}
\multiput(127.057,73.44)(-.076006,-.03319){12}{\line(-1,0){.076006}}
\multiput(125.233,72.644)(-.076006,-.03319){12}{\line(-1,0){.076006}}
\multiput(123.409,71.847)(-.076006,-.03319){12}{\line(-1,0){.076006}}
\multiput(121.585,71.05)(-.076006,-.03319){12}{\line(-1,0){.076006}}
\multiput(119.761,70.254)(-.076006,-.03319){12}{\line(-1,0){.076006}}
\multiput(117.937,69.457)(-.076006,-.03319){12}{\line(-1,0){.076006}}
\multiput(116.112,68.661)(-.076006,-.03319){12}{\line(-1,0){.076006}}
\multiput(114.288,67.864)(-.076006,-.03319){12}{\line(-1,0){.076006}}
\multiput(112.464,67.068)(-.076006,-.03319){12}{\line(-1,0){.076006}}
\multiput(110.64,66.271)(-.076006,-.03319){12}{\line(-1,0){.076006}}
\multiput(108.816,65.475)(-.076006,-.03319){12}{\line(-1,0){.076006}}
\multiput(106.992,64.678)(-.076006,-.03319){12}{\line(-1,0){.076006}}
\multiput(105.168,63.881)(-.076006,-.03319){12}{\line(-1,0){.076006}}
\multiput(103.344,63.085)(-.076006,-.03319){12}{\line(-1,0){.076006}}
\multiput(101.519,62.288)(-.076006,-.03319){12}{\line(-1,0){.076006}}
\multiput(99.695,61.492)(-.076006,-.03319){12}{\line(-1,0){.076006}}
\multiput(97.871,60.695)(-.076006,-.03319){12}{\line(-1,0){.076006}}
\multiput(96.047,59.899)(-.076006,-.03319){12}{\line(-1,0){.076006}}
\multiput(94.223,59.102)(-.076006,-.03319){12}{\line(-1,0){.076006}}
\multiput(92.399,58.306)(-.076006,-.03319){12}{\line(-1,0){.076006}}
\multiput(90.575,57.509)(-.076006,-.03319){12}{\line(-1,0){.076006}}
\multiput(88.75,56.712)(-.076006,-.03319){12}{\line(-1,0){.076006}}
\multiput(86.926,55.916)(-.076006,-.03319){12}{\line(-1,0){.076006}}
\multiput(85.102,55.119)(-.076006,-.03319){12}{\line(-1,0){.076006}}
\multiput(83.278,54.323)(-.076006,-.03319){12}{\line(-1,0){.076006}}
\multiput(81.454,53.526)(-.076006,-.03319){12}{\line(-1,0){.076006}}
\multiput(79.63,52.73)(-.076006,-.03319){12}{\line(-1,0){.076006}}
\multiput(77.806,51.933)(-.076006,-.03319){12}{\line(-1,0){.076006}}
\multiput(75.981,51.137)(-.076006,-.03319){12}{\line(-1,0){.076006}}
\multiput(74.157,50.34)(-.076006,-.03319){12}{\line(-1,0){.076006}}
\multiput(72.333,49.544)(-.076006,-.03319){12}{\line(-1,0){.076006}}
\multiput(70.509,48.747)(-.076006,-.03319){12}{\line(-1,0){.076006}}
\multiput(68.685,47.95)(-.076006,-.03319){12}{\line(-1,0){.076006}}
\multiput(66.861,47.154)(-.076006,-.03319){12}{\line(-1,0){.076006}}
\multiput(65.037,46.357)(-.076006,-.03319){12}{\line(-1,0){.076006}}
\multiput(63.212,45.561)(-.076006,-.03319){12}{\line(-1,0){.076006}}
\multiput(61.388,44.764)(-.076006,-.03319){12}{\line(-1,0){.076006}}
\multiput(59.564,43.968)(-.076006,-.03319){12}{\line(-1,0){.076006}}
\multiput(57.74,43.171)(-.076006,-.03319){12}{\line(-1,0){.076006}}
\multiput(55.916,42.375)(-.076006,-.03319){12}{\line(-1,0){.076006}}
\multiput(54.092,41.578)(-.076006,-.03319){12}{\line(-1,0){.076006}}
%\end
\put(196.25,103.5){\circle*{1.581}}
\put(221.25,114){\circle*{1.581}}
\put(53.5,78.5){\circle*{1.581}}
\put(36.75,89.75){\circle*{1.581}}
\put(219.75,120.25){\makebox(0,0)[cc]{$u$}}
\put(189,104.25){\makebox(0,0)[cc]{$x$}}
\put(27.5,91.25){\makebox(0,0)[cc]{$v$}}
\put(61.75,78){\makebox(0,0)[cc]{$y$}}
\end{picture}

}

 \end{center}

Figure 12 depicts the case of two parallel lines. In Figure 13, we consider the case of two secant lines.

\begin{center}
{\scriptsize

\unitlength 0.4mm % = 2.845pt
\linethickness{0.4pt}
\ifx\plotpoint\undefined\newsavebox{\plotpoint}\fi % GNUPLOT compatibility
\begin{picture}(245.912,154.35)(0,0)
\put(245.912,71.138){\makebox(0,0)[cc]{$x_1$}}
\put(128.313,154.35){\makebox(0,0)[cc]{$x_2$}}
\put(119.65,74.55){\makebox(0,0)[cc]{$0$}}
\put(129.1,0){\makebox(0,0)[cc]
{{\bf Fig. 13:} Secant $\digamma$-lines.}}
\put(123.85,12.338){\vector(0,1){130.987}}
\put(11.5,70.613){\vector(1,0){226.012}}
\put(187.9,99.225){\line(-1,0){132.3}}
\put(55.6,99.225){\line(0,-1){57.75}}
\put(25.613,36.975){\makebox(0,0)[cc]{$\mathcal D^\infty$}}
\put(18.863,91.725){\makebox(0,0)[cc]{$\mathcal E_\infty$}}
%\emline(55.5,41.75)(1.75,18.25)
\multiput(55.5,41.75)(-.07711621234,-.03371592539){697}{\line(-1,0){.07711621234}}
%\end
%\emline(187.25,99)(233.75,119)
\multiput(187.25,99)(.0784148398,.0337268128){593}{\line(1,0){.0784148398}}
%\end
\put(0,99){\line(1,0){242.25}}
\put(211.5,109.5){\circle*{1.5}}
\put(55.25,55.5){\circle*{1.5}}
\put(86.5,98.75){\circle*{1.5}}
\put(152.75,98.75){\circle*{1.5}}
\put(212,118){\makebox(0,0)[cc]{$x$}}
\put(45.25,57){\makebox(0,0)[cc]{$y$}}
\put(85.75,106.5){\makebox(0,0)[cc]{$u$}}
\put(152.5,106){\makebox(0,0)[cc]{$v$}}
\end{picture}

}

\end{center}

Traditional axioms of planar Euclidean geometry postulate that: $(a_1)$ a straight line may be drawn between any two points; $(a_2)$ any terminated straight line may be extended indefinitely; $(a_3)$ a circle may be drawn with any given point as center and any given radius, $(a_4)$ all right angles are equal; $(v)$ for any given point not on a given line, there is exactly one line through the point that does not meet the given line. 

It is well known that the fifth axiom does not follow from the first four. Clearly, the geometric implications of the algebraic structure proposed in this paper do not conflict with the first 4 axioms. Not surprisingly, this is not the case for the fifth. In particular:\\

 $(i)$ two distinct parallel lines may have an infinity of common points (Figure 12);
 
  $(ii)$ there may exist  an infinity of $\digamma$-lines passing by two distinct points $x$ and $y$ (Figure 12);

 $(iii)$ two  secant lines may have an infinity of common points (Figure 13);

  $(iv)$  a $\digamma$-line may not be spanned by two distinct points  it contains (Figure 13).\\

In the following   a system of inequations  characterizing a $\digamma$-line that passes from two distinct points is provided.
This we do by considering the suitable notion of determinant introduced in \cite{b19}. The next intermediary result establishes an algebraic description of the limit of a sequence of $\varphi_p$-hyperplanes.  The next result was established in \cite{b20}.

 \begin{prop} \label{Hyperplane}{Let $V$ be the $n\times n$ matrix with  $v_i$ as $i$-th column for each $i$}. Let $V_{(i)}$ be the matrix obtained from $V$ by replacing line $i$ with the transpose of the unit vector $1\!\!1_n$. Suppose that $|V|_\infty\not=0$. {Let $\{H_p(V)\}_{p\in \mathbb N}$ be the sequence of $\varphi_p$-hyperplanes passing trough each point $v_i$}. Then
$$\Lim_{p\longrightarrow \infty}H_p(V)=\Big\{x\in \Real^n: \bigsmileminus_{i\in [n]} |V_{(i)}|_\infty x_i\leq |V|_\infty\leq \bigsmileplus_{i\in [n]} |V_{(i)}|_\infty x_i\Big\}.$$
\end{prop}
\medskip
Given two distinct points $x,y\in \Real^2$, any point  $z\in \mathcal D(x,y)$ satisfies the relation:

\begin{equation}
|z^{}  {-}x, z^{}  {-}y| =\left |\begin{matrix}z_1^{}-x_1&z_1-y_1\\z_2 -x_2&z_2-y_2\end{matrix}\right|=0.
\end{equation}
Equivalently, we have
\begin{equation}
 (x_2-y_2)z_1^{}+(y_1-x_1)z_2^{}=|x,y|.
\end{equation}
The following proposition   establishes an analogous result in the case of the $\digamma$-line $\mathcal D^\infty$. This algebraic characterization is an immediate consequence of Proposition  \ref{Hyperplane} (see \cite{b20}). Let $x,y\in \Real^2$ with $x\not=y$. Then:
\begin{equation} \mathcal D^\infty(x,y)=
\Big\{z\in \Real^2: (  x_2\boxminus y_2)z_1\stackrel{-}{\smile}( y_1\boxminus x_1 ) z_2\leq |x,y|_\infty \leq (  x_2\boxminus y_2)z_1\stackrel{+}{\smile}( y_1\boxminus x_1 ) z_2\Big\}.\end{equation}

In the following, we show that the equations of two parallel lines in the plane obey algebraic rules analogous to those of the linear case.

\begin{lem} Let $\mathcal D^\infty$ and $\mathcal E^\infty$ be two parallel $\digamma$-lines of $\Real^2$. Then there exists $a\in \Real^2$, $c,d\in \Real$ such that:
$$\mathcal D^\infty=\{z\in \Real^2: \langle  a, z\rangle_\infty^-\leq c\leq \langle  a, z\rangle_\infty^+\leq c  \}$$
and
$$\mathcal E^\infty=\{z\in \Real^2: \langle  a, z\rangle_\infty^-\leq d\leq \langle  a, z\rangle_\infty^+\leq c  \}.$$

\end{lem}
{\bf Proof:} Suppose now that $\mathcal D^\infty $ and $\mathcal E^\infty $ are  two $\digamma$-line respectively spanned from $x,y$ and $u,v$, where $x,y,u,v\in \Real^2$. By definition we have:

$$ \mathcal E^\infty=
\Big\{z\in \Real^2: (  u_2\boxminus v_2)z_1\stackrel{-}{\smile}( v_1\boxminus u_1 ) z_2\leq |u,v|_\infty \leq (  u_2\boxminus v_2)z_1\stackrel{+}{\smile}( v_1\boxminus u_1 ) z_2\Big\}.$$
If $\mathcal D^\infty$ and $\mathcal E^\infty$ are parallel then $(u\boxminus v)\propto (x\boxminus y)$. Therefore,  there is a real number $\alpha\not=0$ such that $(u\boxminus v)=\alpha (x\boxminus y)$. Hence:
$$ \mathcal E^\infty=
\Big\{z\in \Real^2: \alpha (  x_2\boxminus y_2)z_1\stackrel{-}{\smile}\alpha( y_1\boxminus x_1 ) z_2\leq   |u,v|_\infty \leq \alpha (  x_2\boxminus y_2)z_1\stackrel{+}{\smile}\alpha( y_1\boxminus x_1 ) z_2\Big\}.$$
Suppose that $\alpha>0$, then:
$$ \mathcal E^\infty=
\Big\{z\in \Real^2:  (  x_2\boxminus y_2)z_1\stackrel{-}{\smile} ( y_1\boxminus x_1 ) z_2\leq  \alpha^{-1} |u,v|_\infty \leq   (  x_2\boxminus y_2)z_1\stackrel{+}{\smile} ( y_1\boxminus x_1 ) z_2\Big\}.$$
We deduce the result taking $a=x\boxminus y$, $c=|x,y|_\infty$, and $d=\alpha^{-1} |u,v|_\infty$. If $\alpha<0$, then:
$$ \mathcal E^\infty=
\Big\{z\in \Real^2:  (  x_2\boxminus y_2)z_1\stackrel{-}{\smile} ( y_1\boxminus x_1 ) z_2\geq  \alpha^{-1} |u,v|_\infty \geq   (  x_2\boxminus y_2)z_1\stackrel{+}{\smile} ( y_1\boxminus x_1 ) z_2\Big\}.$$
Since $(  x_2\boxminus y_2)z_1\stackrel{-}{\smile} ( y_1\boxminus x_1 )=-\big((  y_2\boxminus x_2)z_1\stackrel{+}{\smile} ( x_1\boxminus y_1 )\big)$ and
 $(  x_2\boxminus y_2)z_1\stackrel{+}{\smile} ( y_1\boxminus x_1 )=-\big((  y_2\boxminus x_2)z_1\stackrel{-}{\smile} ( x_1\boxminus y_1 )\big)$ we deduce that:
 $$ \mathcal E^\infty=
\Big\{z\in \Real^2:  (  y_2\boxminus x_2)z_1\stackrel{-}{\smile} ( x_1\boxminus y_1 ) z_2\leq |\alpha|^{-1} |u,v|_\infty \leq   (  y_2\boxminus x_2)z_1\stackrel{+}{\smile} ( x_1\boxminus y_1 ) z_2\Big\}.$$
Now  note that, since $\mathcal D^\infty(x,y)=\mathcal D^\infty(y,x)$, we also have 
 $$ \mathcal D^\infty=
\Big\{z\in \Real^2:  (  y_2\boxminus x_2)z_1\stackrel{-}{\smile} ( x_1\boxminus y_1 ) z_2\leq   |y,x|_\infty \leq   (  y_2\boxminus x_2)z_1\stackrel{+}{\smile} ( x_1\boxminus y_1 ) z_2\Big\}.$$
Hence, taking $a=y\boxminus x$, $c=|y,x|_\infty$ and $d= |\alpha|^{-1} |u,v|_\infty$ yields the result. $\Box$\\

\begin{expl}Suppose that $x=(3,1)$, $y=(1,-2)$, $u=(-2, 4)$ and $v=(-6, 1)$. We have:
$$x\boxminus y=\begin{pmatrix}
3\boxminus 1\\1\boxminus (-2)
\end{pmatrix}=\begin{pmatrix}
3 \\  2
\end{pmatrix}\quad \text{and}\quad u\boxminus v=\begin{pmatrix}
-2\boxminus (-6)\\4\boxminus 1
\end{pmatrix}=\begin{pmatrix}
6 \\ 4
\end{pmatrix}.$$
It follows that $u\boxminus v=2(x\boxminus y).$
Moreover,
$$|x,y|_\infty=3\cdot(-2)\boxminus 1\cdot 1=-6\quad\text{and}\quad |u,v|_\infty=(-2)\cdot 1 \boxminus (-6)\cdot 4= 24.$$
Therefore:
$$ \mathcal D^\infty(x,y)=
\Big\{z\in \Real^2:  \big(  1\boxminus (-2)\big)z_1\stackrel{-}{\smile} ( 1\boxminus 3 ) z_2\leq   -6 \leq   \big(  1\boxminus (-2)\big)z_1\stackrel{+}{\smile} ( 1\boxminus 3 ) z_2\Big\}.$$
Hence:
$$ \mathcal D^\infty(x,y)=
\Big\{z\in \Real^2:  2z_1\stackrel{-}{\smile}(- 3) z_2\leq   -6 \leq   2 z_1\stackrel{+}{\smile}(- 3)z_2\Big\}.$$
In addition we have:
 
$$ \mathcal D^\infty(u,v)=
\Big\{z\in \Real^2:  \big(  4\boxminus (-1)\big)z_1\stackrel{-}{\smile} \big( (-6)\boxminus (-2) \big) z_2\leq   -6 \leq   \big(  4\boxminus (-1)\big)z_1\stackrel{+}{\smile} \big( (-6)\boxminus (-2)  \big) z_2\Big\}.$$
Hence:
$$ \mathcal D^\infty(u,v)=
\Big\{z\in \Real^2:  (  -4)z_1\stackrel{-}{\smile}  6 z_2\leq   24 \leq   (-4) z_1\stackrel{+}{\smile} 6 z_2\Big\}.$$
Equivalently, by changing the sign and simplifying by 2:
$$ \mathcal D^\infty(u,v)=
\Big\{z\in \Real^2:  2z_1\stackrel{-}{\smile}(- 3) z_2\leq   -12 \leq   2 z_1\stackrel{+}{\smile}(- 3)z_2\Big\}.$$
We can notice that the coefficients multiplying $z_1$ and $z_2$ are the same the difference relating only to the constant term.

\end{expl}

\section*{Conclusion}

In this article, we have explored certain geometric properties emerging from the deformation of a vector space structure defined over a scalar field. In particular, we have focused on notions of lines that can be seen as special limits of linear varieties. For future research, it would be valuable to investigate the connections between these notions and those typically derived from Max-Plus algebras, which are traditionally used in the formalism of tropical mathematics, as initiated in \cite{s88} and \cite{p94}.
\\

\noindent {\bf Compliance with Ethical Standards:} Not applicable.\\ 

\noindent {\bf Data Availability Statement (DAS):}  No new data were created or analyzed in this study. Data sharing is not applicable to this article.

\end{document}